\documentclass{amsart}
\usepackage[english]{babel}
\usepackage{amsmath, amssymb, amsfonts, amsthm, mathtools, enumerate}
\usepackage{extpfeil}

\usepackage[
colorlinks=true,
linkcolor=blue!70!black,
citecolor=blue!70!black,
urlcolor=blue!70!black
]{hyperref}
\usepackage{hhline}

\usepackage{pgfplots}
\pgfplotsset{compat=1.18}
\usepackage{tikz}
\usepackage{graphics}
\usepackage{xcolor}

\newtheorem{theorem}{Theorem}[section]
\newtheorem{definition}[theorem]{Definition}
\newtheorem{lemma}[theorem]{Lemma}
\newtheorem{prop}[theorem]{Proposition}
\newtheorem{corollary}[theorem]{Corollary}
\newtheorem{remark}[theorem]{Remark}

\newtheorem{innercustomthm}{}
\newenvironment{customthm}[1]
  {\renewcommand\theinnercustomthm{#1}\innercustomthm}
  {\endinnercustomthm}

\newcommand{\inv}{^{-1}}

\newcommand{\card}{\mathop{\mathbf{card}}}
\newcommand{\BR}{\mathbb{R}}

\newcommand{\BN}{\mathbb{N}}

\newcommand{\R}{\BR}
\newcommand{\e}{\varepsilon}
\newcommand{\id}{\mathrm{id}}
\newcommand{\mcL}{\mathcal L}

\newcommand{\sgn}{\mathop{\mathrm{sgn}}}
\newcommand{\abs}[1]{\left|#1\right|}
\newcommand{\Abs}[1]{\Big|#1\Big|}

\newcommand{\derder}[2][\tau]{\frac{\partial^2 #2}{\partial #1^2}}

\newcommand{\const}{\mathcal{O}(1)}

\newcommand{\ind}[1]{\mathbf{1}_{#1}}

\newcommand{\restr}[1]{\big|_{#1}}
\newcommand{\Restr}[1]{\Big|_{#1}}

\newcommand{\TV}{\mathop{\mathrm{TV}}\nolimits}
\newcommand{\cont}{\mathrm{cont}}
\newcommand{\atom}{\mathrm{atom}}
\newcommand{\Loneloc}{L^1_{\mathrm{loc}}}
\newcommand{\norm}[1]{\left\|#1\right\|}
\newcommand{\BV}{\mathop{\mathrm{BV}}}

\newcommand{\RcOp}{ \mathbf{Cut}}
\newcommand{\Rc}[2]{\RcOp_{#1, #2}}
\newcommand{\Cr}[1]{\mathbf{Crop}_{#1}}

\newcommand{\conv}{\mathop{\mathbf{conv}}}
\newcommand{\conc}{\mathop{\mathbf{conc}}}
\newcommand{\Height}{\mathbf{H}}
\newcommand{\cost}{\mathop{\mathbf{cost}}\nolimits}

\newcommand{\scal}[2]{\left\langle#1\ ,\ #2\right\rangle}

\newcommand{\bleft}{ b^{\rm left}}
\newcommand{\bright}{b^{\rm right}}
\newcommand{\bridges}{\mathop{\mathbf{Br}}\nolimits}
\newcommand{\vspan}{\mathop{\mathbf{span}}}

\newcommand{\xor}{\mathop{\resizebox{1.7ex}{1.4ex}{$\triangle$}}}

\newcommand{\udhp}{\mathop{%
				\,\tikz[baseline =0.2ex, x=2ex, y=2ex]{
					\draw[line width = 0.1ex, line cap=round, rounded corners = 0.025ex, fill=black!15]
					(0,0)--(1,1)--(1,0) -- cycle;
			}\,}}
\newcommand{\smalludhp}{\mathop{%
				\,\tikz[baseline =0.15ex, x=1.5ex, y=1.5ex]{
					\draw[line width = 0.1ex, line cap=round, rounded corners = 0.025ex, fill=black!15]
					(0,0)--(1,1)--(1,0) -- cycle;
			}\,}}

\newcommand{\mcA}{\mathcal A}
\newcommand{\mcF}{\mathcal F}
\newcommand{\mcS}{\mathcal S}
\newcommand{\mcB}{\mathcal B}
\newcommand{\mcG}{\mathcal G}
\newcommand{\mcD}{\mathcal D}
\newcommand{\mcU}{\mathcal U}

\newcommand{\mcP}{\mathcal P}

\newcommand{\Aquadr}{\mathtt{A}^{\mathrm{quadr}}}
\newcommand{\Atrans}{\mathtt{A}^{\mathrm{trans}}}
\newcommand{\Acubic}{\mathtt{A}^{\mathrm{cubic}}}
\newcommand{\Acanc}{\mathtt{A}^{\mathrm{canc}}}
\newcommand \mfQ {\mathfrak{Q}}
\newcommand{\Qknown}{\mathcal{Q}^{\mathrm{GBB}}}
\newcommand{\Qtrans}{\mathcal{Q}^{\mathrm{trans}}}
\newcommand{\Qcubic}{\mathcal{Q}^{\mathrm{cubic}}}
\newcommand{\Qcost}{\mathcal{Q}^{\mathrm{cost}}}

\newcommand{\indfixbeta}[1]{\colorbox{cyan!30}{#1}}
\newcommand{\indfixC}[1]{\colorbox{green!30}{#1}}

\newcommand{\drp}{\delta_{\mathrm{RP}}}
\newcommand{\Beone}{\overline B_{\drp}(0)}
\newcommand{\Beonehalf}{\overline B_{\drp/2}(0)}

\newcommand{\I}{\mathbf{I}}

\newcommand{\Jleft}{J_{\mathrm{left}}}
\newcommand{\Jright}{J_{\mathrm{right}}}

\author{Egor Baitenov}
\author{Stefano Modena}
\address{GSSI, I-67100 L'Aquila, Italy}
\email{egor.baitenov@gssi.it, stefano.modena@gssi.it}
\date{\today}

\begin{document}

\title[ 
  Cost functional for systems of conservation laws]
 {The cost functional  for hyperbolic systems \\ of conservation laws}

\begin{abstract}
    For hyperbolic systems of conservation laws (without any genuine nonlinearity of linear degeneracy assumption on the flux) we construct a quadratic interaction functional, which controls the change of wavefront speed across interactions and it is local in time. This extends the analysis of \cite{BaitMod1} from scalar conservation laws to $n \times n$ systems. 
\end{abstract}

\maketitle

\section{Introduction}

We consider the $n \times n$ system of conservation laws
\begin{equation}
    \label{eq:syst-cons-law}
    u_t + f(u)_x  = 0, 
\end{equation}
with initial datum
\begin{equation}
    \label{eq:initial-datum}
    u|_{t=0} = \bar u,
\end{equation}
for the unknown $u \colon [0, + \infty) \times \R \to \R^n$. Here $f\colon \overline{\Omega} \subseteq \R^n \to \R^n$ is a given \emph{flux}. We assume $\Omega$ is open, $f$ is smooth, with uniformly bounded second derivative, and \emph{strictly hyperbolic}, i.e. for every $u \in \Omega$, $\nabla f(u)$  has $n$ distinct real eigenvalues
\begin{equation*}
    \lambda_1(u) < \ldots < \lambda_n(u).
\end{equation*}
No genuine nonlinearity (GNL) or linear degeneracy (LD) assumption is made on~$f$. 

We work in the vanishing viscosity setting of Bianchini and Bressan \cite{BianchiniBressan2005}, namely there is a subset $\mcD \subseteq \Loneloc(\R; \R^n)$, containing all $\BV$ maps with sufficiently small total variation and left state $u(-\infty)$ belonging to a fixed compact set $K \subset \Omega$, and a Lipschitz semigroup
\begin{equation*}
    S \colon [0, +\infty) \times \mcD \to \mcD, \qquad (t, \bar u) \mapsto S_t \bar u
\end{equation*}
such that, for every $\bar u \in \mcD$, $u(t,x) := S_t \bar u$ is the (vanishing viscosity) solution to the Cauchy problem~\eqref{eq:syst-cons-law}--\eqref{eq:initial-datum}. 

We prove the following two main theorems, stated for now in a quite informal way.
They claim the existence of a \emph{functional} $\Upsilon$ with some properties. The first one contains the \emph{fundamental} properties of $\Upsilon$, the second one contains a list of \emph{additional} properties. 

\begin{customthm}{Theorem~A}
\label{thm-main-a}
    There is a \emph{local-in-time} functional
    \begin{equation*}
        \Upsilon \colon \mathcal{D} \to \R_{\geqslant 0},
    \end{equation*}
    such that
    \begin{enumerate}[(i)]
        \item\label{thm-main:qb} {\underline{Quadratic bound}:}
        $\Upsilon(u) = \const \TV(u)^2$ 
        (i.e. $\Upsilon$ is quadratic w.r.t. the  total variation of its argument),\footnote{The (standard) notations $A = \const B$, $A \leqslant \const B$ mean that $A \leqslant CB$ for some constant $C>0$ depending only the flux function $f$ in~\eqref{eq:syst-cons-law}.}
        \item\label{thm-main:ie} {\underline{Interaction estimate}:}
        in each \emph{interaction between wavefronts} $\Upsilon$ decreases, and its decrease bounds  the \emph{change of speed} of the wavefronts involved in the interaction,
        \item\label{thm-main:pers} {\underline{Persistency}:}  $\Upsilon(S_t u) = \Upsilon(u)$, if $u$ is piecewise constant (with uniform spatial separation of jumps) and $t$ is small enough.
    \end{enumerate}
\end{customthm}

\begin{customthm}{Theorem~B}
\label{thm-main-b}
The functional $\Upsilon$ satisfies the following additional properties:
    \begin{enumerate}[(i)]
        \item \label{thm-main:rcut} {\underline{Riemann cuts}:} $\Upsilon$ is non-increasing across Riemann cuts.
        \item \label{thm-main:appr} {\underline{Approximation}:}  $\Upsilon(u) = \sup_v \Upsilon(v)$, where the supremum is taken over all piecewise constant samplings $v$ of $u$.
        \item \label{thm-main:stab} {\underline{Stability}:}
             $\abs{\Upsilon(u)-\Upsilon(\Phi(u))}= 
            \const(1+M^2)M\norm{\Phi-\id}_{W^{1,\infty}}$, for every\linebreak $\Phi \colon \Omega \to \R^n$ Lipschitz, where
            $M := \max\{\TV(u),\TV(\Phi\circ u)\}$ and $\id(p) \equiv p$.
        \item \label{thm-main:lsc} {\underline{Lower semicontinuity}:} If $u_m \to u$ in $\Loneloc$, then $\Upsilon(u) \leqslant \liminf_m \Upsilon(u_m)$.
        \item \label{thm-main:lyap} {\underline{Lyapunov property}:} For every $\bar u \in \mcD$, $t \mapsto \Upsilon(S_t \bar u)$ is non-increasing.
        \item \label{thm-main:invar} {\underline{Invariance}:} $\Upsilon$ is invariant under reparametrization of independent variable.
    \end{enumerate}
\end{customthm}

\ref{thm-main-a} was proven in the scalar case $n=1$ in \cite{BaitMod1} (with $\mcD = \BV(\R; \overline \Omega)$). In the present paper we extend the scalar theory to the case of systems. The core ideas concerning the construction of $\Upsilon$ and the proof of~\ref{thm-main-a} are essentially scalar, and they are already contained in~\cite{BaitMod1}. Extension to systems is mainly a matter of (often nontrivial) technicalities. 

\ref{thm-main-b} was as well stated  in~\cite{BaitMod1} for the scalar case, but no proof was given there. In the present paper we present a complete proof also of \ref{thm-main-b}, which  holds both for the scalar case $n=1$ and for the general system case $n \geqslant 1$.

The rigorous definition of $\Upsilon$ and the precise formulation of the statements in \ref{thm-main-a} and \ref{thm-main-b} require some notations and they are therefore postponed to Sections~\ref{s:rigorous-upsilon}--\ref{s:further-properties}. 

For the rest of this introduction, we assume  the reader has a vague idea about the notions of \emph{wavefronts}, \emph{strength} and \emph{speed} of wavefronts, and \emph{interaction} between wavefronts, and we briefly discuss the importance of controlling the \emph{change of speed} of wavefronts involved in an interaction by means of \emph{local-in-time} functionals (\ref{thm-main-a}). We also comment on the property listed in \ref{thm-main-b}. For a more extensive introduction to these topics, we refer the reader to the introduction of~\cite{BaitMod1}.

\subsection{On the change of speed across interactions}

Understanding how wavefronts interact is one of the central question in the analysis of the system~\eqref{eq:syst-cons-law}. When two (or more) wavefronts \emph{interact}, they undergo a change of their strength and of their speed. 

A classical approach, due to Glimm \cite{Glimm1965}, in the theory of system of conservation laws consists in constructing \emph{interaction functionals} which decrease across interactions, and whose decrease bounds the change in strength and/or in speed of wavefronts involved in the interaction. 

    \ref{thm-main-a} provides one such functional, 
     bounding the change in speed, see item~\eqref{thm-main:ie}. Item~\eqref{thm-main:qb} ensures that $\Upsilon(u)$ is uniformly bounded by the square of $\TV(u)$, while item~\eqref{thm-main:pers} guarantees that $t \mapsto \Upsilon(S_t \bar u)$ is constant as long as no interactions take place. 

Combining items~\eqref{thm-main:qb}--\eqref{thm-main:pers}, one can show, by standard techniques, that for every piecewise constant approximate solution $u$ (constructed e.g. by Glimm scheme\footnote{For wavefront tracking solutions \cite{Bressan2000, AM-WFT}, in order to ensure decay in time, the definition of $\Upsilon$ (see Section~\ref{s:rigorous-upsilon}) must be slightly modified including non-physical waves.}) $t \mapsto \Upsilon(u(t,\cdot))$ is non-increasing, and thus the global change of speed of wavefronts on the whole domain $[0, +\infty) \times \R$ is uniformly bounded by $\const \TV(u(0, \cdot))^2$.

We refer the reader to~\cite[Sections~1.1--1.3]{BaitMod1} for an extensive discussion about Glimm's strategy and the importance of controlling the global change of speed of wavefronts in approximate solutions. 

\subsection{On the locality in time}
\label{ss:intro-locality-time}

Functionals bounding the global change of speed of wavefronts in approximate solutions were already present in the literature  (\cite{BM-CMP, M-Forward} and also~\cite{AM2011, BM-Glimm, BM-CaseStudy, BM-scalar}). 

The main features of~$\Upsilon$ in~\ref{thm-main-a}, which make it substantially different than the functionals in~\cite{BM-CMP, M-Forward} (and the other references mentioned above) are, first, that~$\Upsilon(u)$ is defined for every $u \in \mcD$ (and not only for piecewise constant maps, for which wavefront tracking or Glimm approximate solutions can be constructed). 

Second, and more important, $\Upsilon$ in \ref{thm-main-a} is \emph{local-in-time}, i.e., for any given $u \in \mcD$, the value $\Upsilon(u)$ depends only on the profile $u$, and it does not depend on the (approximate) solution to~\eqref{eq:syst-cons-law} starting from $u$. 

This was not the case in~\cite{BM-CMP, M-Forward} (and the other references mentioned above), where wave tracing algorithms needed to be developed in order to trace the evolution of wavefronts in $u$ along the (approximate) solution, thus making the functionals in \cite{BM-CMP, M-Forward} highly nonlocal in time.  

We refer the reader to \cite[Section 1.4]{BaitMod1} for extensive discussion about the importance of having a \emph{local-in-time} $\Upsilon$, defined for every $u \in \mcD$. 

\subsection{On the additional properties of \texorpdfstring{$\Upsilon$}{Upsilon}}
\label{s:on-add-prop}
While \ref{thm-main-a} contains the fundamental properties of $\Upsilon$, as discussed before, we list  in \ref{thm-main-b} some \emph{additional} properties, which are useful when one wants to use~$\Upsilon$ in the analysis of~\eqref{eq:syst-cons-law}. 

We observe first of all that, for  the nonlocal functionals bounding the change in speed constructed in \cite{BM-CMP, M-Forward}, the additional properties listed in \ref{thm-main-b} do not hold,  mostly because they do not even make sense, being the functionals in \cite{BM-CMP, M-Forward} nonlocal in time and defined only for piecewise constant maps. 

On the other side, there are other well-known Lyapunov functionals in the literature, bounding  the change in \emph{strength} of wavefronts across interactions (and not the change in speed), which are defined for all $u \in \mcD$ and are local in time: e.g. the sum $\Qknown = \Qtrans + \tfrac{1}{4} \Qcubic$ of the Glimm transversal functional~$\Qtrans$~(\cite{Glimm1965}, see also~\cite{Bressan1995}) and the cubic functional~$\Qcubic$ by Bianchini and Bressan (\cite{Bianchini2002, Bianchini2003}), whose definitions will be recalled later (see Section~\ref{ss:def-V-k-Qknown}). For these Lyapunov functionals, properties analog to those in \ref{thm-main-b} hold true (see also \cite{Baiti1997}).

The fact that our functional~$\Upsilon$ is defined for all~$u \in \mcD$ and it is local-in-time suggests that~$\Upsilon$ may share the same ``additional properties'' as~$\Qknown$. This is indeed true, and such additional properties are listed in \ref{thm-main-b}. Despite being natural and somehow expected, the rigorous proof of \ref{thm-main-b} require some effort. 

We comment now briefly on items~\eqref{thm-main:rcut}--\eqref{thm-main:appr}, while items \eqref{thm-main:stab}--\eqref{thm-main:invar} should already be clear from the statement. 

\begin{figure}[h]
    \centering
    \includegraphics{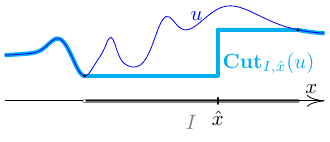}
    \caption{Riemann cut.}
    \label{fig:riemann_cut}
\end{figure}

Given $u \in \mcD$, a \emph{Riemann cut} $\Rc I{\hat x}$ is an operation that consists of collapsing the variation of $u$ on a given interval $I$ into a single jump located at some $\hat x \in I$ (see Fig.~\ref{fig:riemann_cut}). A Riemann cut can be interpreted as a simultaneous interaction of all wavefronts in $I$. Item~\eqref{thm-main:rcut} claims that~$\Upsilon$ decreases across Riemann cuts, consistently with \ref{thm-main-a}, item~\eqref{thm-main:ie}. This property will be made rigorous in Corollary~\ref{cor:riemann-cut}. 

The approximation property in item~\eqref{thm-main:appr} will be made rigorous in Theorems~\ref{prop:sup-prop-crop} and~\ref{thm:sup-prop-cut}. Theorem~\ref{prop:sup-prop-crop} states that $\Upsilon(u) = \sup_v \Upsilon(v)$, where the supremum is taken over all maps $v$ obtained by \emph{cropping} $u$ in such a way that the $v$'s are constant  outside arbitrarily large compact sets. Theorem~\ref{thm:sup-prop-cut}, on the other side, states   that $\Upsilon(u) = \sup_v \Upsilon(v)$, where the supremum is taken over all piecewise constant maps $v$ obtained from $u$ by means of a finite number or Riemann cuts.

We stress once again that  \ref{thm-main-b} do not hold for nonlocal functionals like the ones in \cite{BM-CMP, M-Forward}. Take, for instance, the most simple property listed above, the invariance w.r.t. reparametrization of independent variable (item~\eqref{thm-main:invar}), whose proof for our $\Upsilon$ is very natural. 
Nevertheless, for the nonlocal functional in \cite{BM-CMP, M-Forward} the invariance property is false: indeed, take a piecewise constant map and slightly  move the position of one of its jump. This is essentially a reparametrization of the independent variable, but it  might lead to complete different approximate evolutions (in wavefront tracking or Glimm scheme), and therefore it might drastically affect functionals whose value depend on the global evolution of wavefronts on $[0, +\infty)$.     

\subsection{Structure of the paper}

As in the scalar case~\cite{BaitMod1}, the most important object in the construction of~$\Upsilon$ is the $\cost$ function. We recall the definition of $\cost$ and its main properties (already proven in~\cite{BaitMod1}) in Section~\ref{s:cost-func-and-prop-recall}. 

Section~\ref{s:rigorous-upsilon} is devoted to the rigorous definition of~$\Upsilon$, for which we need also to recall the construction of vanishing viscosity solutions to Riemann problems, and the definition of the well-known functionals~$\Qtrans$ (introduced in~\cite{Glimm1965}) and~$\Qcubic$ (introduced in \cite{Bianchini2002, Bianchini2003}).  

In Section~\ref{s:system-quadratic-bound} we prove the quadratic bound (Theorem~\ref{thm:quadr-bd}, corresponding to~\ref{thm-main-a}, item~\eqref{thm-main:qb}).  In Section~\ref{s:interaction-estimate} we give a precise definition of \emph{interaction between wavefronts} and of the corresponding \emph{change of speed}, so that we can rigorously state and prove the main interaction estimate (Theorem~\ref{thm:interaction-estimate}, corresponding to~\ref{thm-main-a}, item~\eqref{thm-main:ie}). Prop.~\ref{prop:var-cost} contains the core of the proof. In Section~\ref{s:persistency} we prove the persistence property, namely that, for $\bar u$ piecewise constant, $\Upsilon(S_t \bar u)$ is constant for short times (precisely as long as no interaction takes place) (Theorem~\ref{thm:const-if-no-interaction}, corresponding to~\ref{thm-main-a}, item~\eqref{thm-main:pers}).

Sections~\ref{s:crop},~\ref{s:cut} are devoted to rigorously state and prove the approximation property (\ref{thm-main-b}, item~\eqref{thm-main:appr}). More precisely, in Def.~\ref{def:crop} we introduce the \emph{crop} of a map~$u$ and in Theorem~\ref{prop:sup-prop-crop}  we prove the approximation property w.r.t. cropping. In Def.~\ref{def:riemann-cut} we give a rigorous definition of \emph{Riemann cut} and in Theorem~\ref{thm:sup-prop-cut} we prove the approximation  property w.r.t. cutting. As a corollary we get also that~$\Upsilon$ is decreasing across Riemann cuts (Corollary~\ref{cor:riemann-cut}  corresponding to \ref{thm-main-b}, item~\eqref{thm-main:rcut}). 

In Section~\ref{s:cont-wrt-states} we prove the stability of~$\Upsilon$ w.r.t. perturbation of states (Theorem~\ref{thm:cont-wrt-states} corresponding to \ref{thm-main-b}, item~\eqref{thm-main:stab}).  

Finally in Section~\ref{s:further-properties} we prove lower semicontinuity, Lyapunov property and invariance (respectively in Theorems~\ref{thm:lsc},~\ref{thm:decr-along-semigroup},~\ref{thm:reparametr-state-curve}, corresponding to \ref{thm-main-b}, items~\eqref{thm-main:lsc}, \eqref{thm-main:lyap},~\eqref{thm-main:invar}).

\subsection{Notation}
\label{ss:notation}
We use the following notation:
\begin{itemize}
	\item $\delta_{kh} $ for Kronecker delta;
    \item $\lambda_k(u)$ ($u\in\Omega$) the $k$-th eigenvalue of the matrix $A=Df(u)$ with the ordering $\lambda_1(u)<\ldots<\lambda_n(u)$;
    \item $r_k(u)$ and $\ell_k(u)$ ($u\in\Omega$) the $k$-th right and left eigenvectors of the matrix $A=Df(u)$, respectively, with the normalization $\scal{r_k}{r_k} = 1$, $\scal{\ell_{k'}}{r_{k''}} = \delta_{k'k''}$;
    
    \item $\conv_{[a,c]} g$ (resp. $\conc_{[a,c]} g$) for the convex (res. convave) envelope of the function~$g$ over the segment $[a,c]$;
    \item $\scal{\ell}{r}$ for the scalar product in $\R^n$;
    \item $u_x $ for the derivative of a function $u\in\BV(\R;\R^n)$ treated as measure;
    \item $C_c(\R)$ for compactly supported continuous functions on $\R$;
    
    \item $\I(s) = \begin{cases}
        [0,s],&s\geqslant 0,\\
        [s,0],& s<0;
    \end{cases}$

    \item for $pq\geqslant 0$ and $g\colon\I(p)\to \R$, $h\colon p+\I(q)\to\R$, $g(p) = h(p)$ we denote
    \begin{align}
    \label{eq:cup-def}
        g\cup h\colon \I(p+q)\to \R,\quad (g\cup h)(y) =\begin{cases}
          g(y), &y\in\I(p),\\
          h(y), &y\in p+\I(q);
        \end{cases} 
    \end{align}
    \item for $pq< 0$ and $g\colon\I(p)\to \R$, $h\colon p+\I(q)\to\R$ we denote
    \begin{align}
    \label{eq:xor-def}
        g\xor h\colon \I(p+q)\to \R,\quad g\xor h =\begin{cases}
          g\restr{\I(p+q)}, &\abs{p}\geqslant\abs{q},\\
          h\restr{\I(p+q)},& \abs{p}<\abs{q};
        \end{cases} 
    \end{align}    
    \item $\mu\restr{A}\colon B\mapsto \mu(A\cap B)$ the restriction of Radon measure $\mu$ onto a Borel set $A$;
    \item $\sgn \mu $ the Radon--Nikodym derivative of measure $\mu$ w.r.t. $\abs\mu$;
    \item $\vspan \{v_1,\ldots,v_m\}$ the convex hull of vectors $\{v_i\}_{i=1}^m$;
    \item $\ind A$ for the indicator function of $A\subset\R$;

    \item $\bigsqcup$ for disjoint union;
    \item for any interval $A$,  $A_+ : = \{x\colon x> a\text{ for every $a\in A$}\}$, $A_- :=\R\setminus A_+\setminus A $;
    \item $\mcL(A)$ and $|A|$ for the Lebesgue measure of $A\subset\R$;
    \item $\psi_\sharp\mu$ for the pushforward of the measure $\mu$ through the map $\psi$.
    
\end{itemize}

\section{\texorpdfstring{$\cost$}{Cost} function and its properties}
\label{s:cost-func-and-prop-recall}

As in the scalar case, the main object in the definition of $\Upsilon$ is the $\cost$ function, introduced
in~\cite[Section 2]{BaitMod1}. We recall now the definition of $\cost$ and its main properties, that will be used throughout the paper. For the proofs, we refer to \cite{BaitMod1}.

\subsection{Definition of \texorpdfstring{$\cost$}{cost}}
\label{ss:cost-def-recall} 
\begin{figure}[h]
			\centering
			\includegraphics{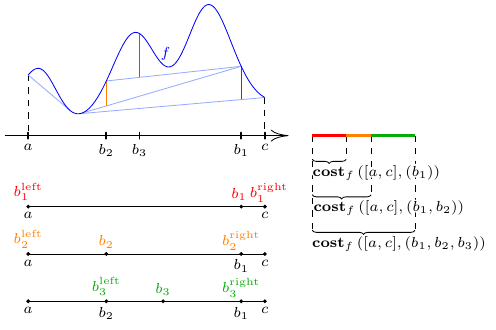}
		\caption{The definition of merge cost.}
		\label{fig:mergercostdef-recall}
\end{figure}

      Let $a\leqslant b \leqslant c$ and $f\in C^{1,1}([a,c])$. Put   \begin{align*}
        \Height_f(a,b,c) =
        f(b) -\Big(\conv_{[a,c]} f\Big)(b).
    \end{align*}
 Suppose now that $b_1,\ldots,b_m \in [a,c]$. Define the \emph{merge cost} of the interval $[a,c]$ w.r.t. to the \emph{mergers} $(b_i)_{i=1}^m$ recursively by
     \begin{align*}
         \cost_f([a,c],(b_i)_{i=1}^m) = \begin{cases}
             \Height_f(a,b_1,c),& m = 1,\\
             \cost_f([a,c],(b_i)_{i=1}^{m-1}) +\Height_f(\bleft_m, b_m, \bright_m),& m>1,
         \end{cases}
     \end{align*}
     where
     \begin{align*}
        \bleft_j& = \max\Big(\{a\}\cup\left( \{b_i\}_{i=1}^{j-1}\cap[a,b_j]\right)\Big),\\  
        \bright_j & = \min\Big(\{c\}\cup\left( \{b_i\}_{i=1}^{j-1}\cap[b_j,c]\right)\Big). 
     \end{align*}
     Alternatively, one easily derives
$
     \cost_f([a,c],(b_i)_{i=1}^m) = \sum_{i=1}^m\Height_f(\bleft_i, b_i, \bright_i).
$
 Notice that, in general, $\cost_f([a,c],(b_i)_{i=1}^m)$ depends on the order of the points $b_i$. An example is shown in Fig.~\ref{fig:mergercostdef-recall}.

We can finally introduce the \emph{cost function} $\cost$.

\begin{definition}\label{def:cost}
        Let $a\leqslant c$ and $f\in C^{1,1}([a,c])$. We define the cost function as 
        \begin{align*}
            \cost_f(a,c) =\sup\Big\{ \cost_f([a,c],(b_i)_{i=1}^m) \, \colon \,  m\in \BN,  \text{ $b_1,\ldots, b_m\in [a,c]$}\Big\}.
        \end{align*}
        Define as well $
            \cost_f(c,a) = \cost_{-f}(a,c)$.
    \end{definition}

\subsection{Quadratic bound on \texorpdfstring{$\cost$}{cost} and interaction-like estimates}
\label{ss:cost-QB-ILE-recall}
First, we recall an upper bound for $\cost$ function. We will refer to it as \emph{quadratic bound} meaning the dependence on the length of the interval.
 \begin{prop}\label{prop:cost_Linfty_est-recall}
        Let $a<c$ and $f\in C^{1,1}([a,c])$, then
        $$
            \cost_f(a,c)\leqslant\frac14(c-a)^2\norm{\left(\derder[x]{f}\right)^-}_{L^\infty}.
        $$
\end{prop}
Secondly, we restate the inequalities relating the values of $\cost$ on subsequent intervals and their union. In view of conservation laws, treating intervals as an interacting jumps, we will refer to these inequalities as \emph{interaction-like estimates}.
\begin{prop}\label{prop:ILE-recall}
		If $a\leqslant b\leqslant c$ and $f\in C^{1,1}([a,c])$, then
        \begin{enumerate}[(i)]
            \item \label{prop:ILE-recall:superadd}$\cost_f(a,c)\geqslant\cost_f(a,b)+\cost_f(b,c)+\Height_f(a,b,c)$,
            \item \label{prop:ILE-recall:subadd}$\cost_f(a,c)\leqslant\cost_f(a,b)+\cost_f(b,c)+
        \frac12
            (b-a)(c-b)\norm{\left(\derder[x]f\right)^-}_{L^\infty}
        .$
        \end{enumerate}
	\end{prop}

\subsection{Persistence property}
\label{ss:cost-persistency}
\begin{figure}[h]
    \centering
    \includegraphics{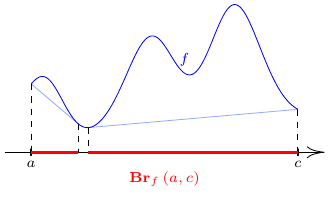}
    \caption{Definition of $\bridges_f(a,c)$}
    \label{fig:bridges}
\end{figure}

We now recall a useful proposition stating that $\cost$ is in fact localized in the gaps between flux function and its convex envelope. For the rigorous statement, we need the following definition.
\begin{definition}\label{def:bridges}
    Let $f\in C([a,c])$. We define the family of \emph{bridges of $f$ above $[a,c]$} denoted by $\bridges_f(a,c)$ as the countable family of maximal open intervals on which $f>\conv_{[a,c]}f$, see Fig.~\ref{fig:bridges}.
\end{definition}
\begin{prop}\label{prop:persistence}
    Let $a \leqslant c$ and $f\in C^{1,1}([a,c])$. Then
    \begin{align*}
        \cost_f(a,c) = \sum_{(p,q)\in\bridges(a,c)}\cost_f(p,q).
    \end{align*}
\end{prop}

\subsection{Stability estimate}
\label{ss:cost-SE-recall}
We now reproduce the estimate quantifying the change of $\cost$ if both function and interval are perturbed.
\begin{prop}\label{prop:stability-recall}
    Let $a,b_1,b_2\in\BR$. Denote $I_i =[\min\{a,b_i\},\max\{a,b_i\}] $, $i = 1,2$. Suppose that $f_i\in C^{1,1}(I_i)$, $i = 1,2$. Then
    \begin{gather*}
        \Abs{\cost_{f_1}(a,b_1) - \cost_{f_2}(a,b_2)}\leqslant\\
        \leqslant\frac12\abs{I_1\cap I_2}\cdot\norm{\frac{\partial ^2f_1}{\partial\tau^2}-\frac{\partial ^2f_2}{\partial\tau^2}}_{L^1(I_1\cap I_2)} +\\
        +\frac12\max\left\{\norm{\frac{\partial ^2f_1}{\partial\tau^2}}_{L^\infty},\norm{\frac{\partial ^2f_2}{\partial \tau^2}}_{L^\infty}\right\}\cdot\frac{\abs{I_1}+\abs{I_2}}2\cdot \abs{b_2-b_1}.
    \end{gather*}
\end{prop}

\section{Rigorous definition of \texorpdfstring{$\Upsilon$}{Upsilon}}
\label{s:rigorous-upsilon}

In this section we give the rigorous definition of the functional $\Upsilon$. In the scalar case (see~\cite[Def.~2.6]{BaitMod1}), $\Upsilon(u)$ was the difference between $C\TV^2(u)$ and sum of $\cost(u_-, u_+)$ over jumps of $u$. In the case of systems, 
the definition is more technical. 

Indeed, $\TV$ is replaced by the \emph{total variation  along curves $V$}, $\cost$ must be computed using \emph{reduced fluxes} and, in addition, the well known Glimm--Bianchini--Bressan  functional $\Qknown = \Qtrans + \tfrac{1}{4}\Qcubic$ (see below in Section~\ref{ss:def-V-k-Qknown}) must be added to the definition of~$\Upsilon$, in order to take care of the presence of~$n$ different families.

In order to introduce all these objects, we need to recall the core building block that is \emph{Riemann problem}. First, we remind how Riemann problems are solved in the general case, where no GNL/LD assumption is made (Section~\ref{ss:sln-rp}). We can then define, for each family $k=1,\dots, n$,  the $k$-\emph{cost} (Section~\ref{ss:def-cost-syst}) and the functionals~$V$ and~$\Qknown$ (Section~\ref{ss:def-V-k-Qknown}). The functional $\Upsilon$ is defined in Section~\ref{ss:def-Upsilon}. Even though $\Upsilon$ is well defined also for (some) $\BV$ maps with large total variation, most of the results in \ref{thm-main-a} and \ref{thm-main-b} hold true only if one assume smallness of total variation. This assumption is conveniently phrased using a \emph{semigroup domain} approach, in the spirit of Bressan (see~\cite{Bressan2000} and the reference therein). We  comment on this approach in Section~\ref{ss:domain-def}. Finally, we recall some basic estimates about solutions to Riemann problems in Section~\ref{ss:basic-rp}.

\subsection{Solution to Riemann problem}
\label{ss:sln-rp}

For given $u^L, u^R \in \Omega$ close enough to zero, we describe now how one constructs~\cite{BianchiniBressan2005} the vanishing viscosity solution to the Riemann problem $(u^L, u^R)$, namely the Cauchy problem to~\eqref{eq:syst-cons-law} with initial datum
$\bar u(x) = u^L \ind{(-\infty, 0)} + u^R \ind{(0, + \infty)}$. We assume, by the smallness of $\abs{u^L}$, $\abs{u^R}$, that $\bar u$ fits the setting of~\cite{BianchiniBressan2005}
and thus admits a unique vanishing viscosity solution.

First of all, one constructs, for each family $k=1,\dots, n$,  so called \emph{generalized eigenvectors} and \emph{generalized eigenvalues}
\begin{equation*}
\begin{aligned}
    \tilde r_k & : \mathcal{U}_k \subseteq \R^n \times \R \times \R \to \R^n, \\
    \tilde \lambda_k & : \mathcal{U}_k \subseteq \R^n \times \R \times \R \to \R,    
\end{aligned}
\end{equation*}
where $\mathcal{U}_k$ is a neighborhood of $(0, 0, \lambda_k(0))$. Recall that $\lambda_k(0)$ is the $k$-th eigenvalue of the matrix $A(0) = Df(0)$. The generalized eigenvalues are strictly separated, i.e. for some $\hat \lambda_0 < \cdots < \hat \lambda_n$, 
\begin{equation}
    \label{eq:separation-eigenvalues}
    \hat \lambda_{k-1} < \tilde \lambda_k(u,v,\sigma) < \hat \lambda_k \text{ for all } (u,v,\sigma) \in \mcU_k.
\end{equation}
Moreover, the generalized eigenvectors and eigenvalues satisfy:
\begin{equation}
\label{eq:unit-length}
    |\tilde r_k (u,v,\sigma) | = 1, \qquad \tilde r_k(u,0,\sigma) = r_k(u), \qquad \tilde \lambda_k(u,0,\sigma) = \lambda_k(u)
\end{equation}
for all $(u,v,\sigma) \in \mathcal{U}_k$. 

\medskip

Let now $\mcU \subset \R^n$ be a bounded neighborhood of~$0$ such that for every $ u^L\in\overline\mcU $ it follows $ (u^L, 0, \lambda_k(u^L))\in\mcU_k$, $k = 1,\ldots, n$.
Given a left state $u^L\in\mcU$  and $k  \in \{1, \dots, n\}$, one shows that there is a curve $s \mapsto T^k_s u^L \in \R^n$, with $T^k_0 u^L=u^L$, defined for $|s| \ll 1$, such that, for each fixed $s$, the Riemann problem $(u^L, T^k_s u^L)$ admits a unique vanishing viscosity solution, that is constructed as follows. There is another curve, denoted by $\gamma_k(u^L, s)$ or simply by $\gamma_k$, called \emph{exact curve} of the~$k$-th family, starting from~$u^L$ with length~$s$,
\begin{equation*}
    \gamma_k(u^L, s) = \gamma_k =(u_k, v_k, \sigma_k) : \I(s) \to \mathcal{U}_k 
\end{equation*}
with $u_k(0) = u^L$, $u_k(s) = T^k_s u^L$, solving, if $s \geqslant 0$, the fixed point problem
\begin{equation}
\label{eq:exact-curve}
\left\{
\begin{aligned}
    u_k(\tau) & = u^L + \int_0^\tau \tilde r_k (u_k(\varsigma), v_k(\varsigma), \sigma_k(\varsigma)) d\varsigma, \\
    v_k(\tau) & = f_k(\tau) - \conv_{[0, s]} f_k(\tau), \\
    \sigma_k(\tau) & = \frac{d}{d\tau} \conv_{[0,s]} f_k(\tau),
\end{aligned}
\right.
\end{equation}
where
\begin{equation*}
    f_k(\tau) = \int_0^\tau \tilde \lambda_k (u_k(\varsigma), v_k(\varsigma), \sigma_k(\varsigma)) d \varsigma,
\end{equation*}
and the same fixed point problem with $\conc_{[s,0]}$ instead of $\conv_{[0,s]}$ if $s \leqslant 0$. The map $f_k$ is called $k$-th \emph{reduced flux} associated to the curve $\gamma_k(u^L,s)$. 

\begin{remark}
    \label{rmk:shifts}
It will sometimes be convenient to work with a \emph{shifted} version of $\gamma_k$, meaning that $\gamma_k$ has domain $a + \I(s)$ instead of $\I(s)$. We may also shift $f_k$ vertically, i.e. set $f(a)$ equal to a prescribed value, possibly different than zero. Notice that~\eqref{eq:exact-curve} is robust under these horizontal and vertical shifts. 
\end{remark}

The solution to the Riemann problem $(u^L, T^k_s u^L)$ is then given by 
\begin{equation}
\label{eq:rp-sln-k}
    \omega_k(t,x) := 
    \begin{cases}
        u^L, & \text{if }  x/t \leqslant \sigma_k(0), \\
        u_k(\tau), & \text{where } \tau = \max\{ \varsigma \in \I(s) \, : \, x/t = \sigma_k(\varsigma) \}, \\
        T^k_s u^L, & \text{if } x/t \geqslant \sigma_k(s).
    \end{cases}
\end{equation}
\begin{figure}[h]
    \centering
    \includegraphics{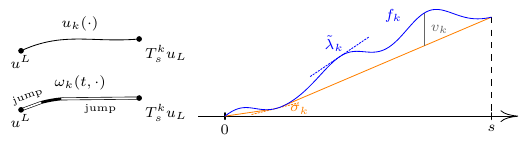}
    \caption{$k$-th exact curve $(u_k, v_k, \sigma_k)$ and the solution~$\omega_k$ to the Riemann problem $(u^L, T_s^ku^L)$ (case $s>0$). Observe the correspondence between the jumps of~$\omega_k$ and the flat parts of~$\sigma_k$, according to~\eqref{eq:rp-sln-k}.}
    \label{fig:sol-curve}
\end{figure}
\begin{remark}
    It is useful to rewrite~\eqref{eq:rp-sln-k} in the form
    \begin{align}\label{eq:rp-sln-k-pf}
        \omega_k(-\infty) = u^L, \qquad \omega_k(t,\cdot)_x = (t\sigma_k)_{\sharp}\Big( (\tilde r_k\circ \gamma_k)\cdot\mcL\restr {\I(s)} \Big)
    \end{align}
    (up to a.e. equality).
    Indeed, let w.l.o.g. $s\geqslant 0$, then, being $\sigma_k$  non-decreasing and $\omega_k(t,\cdot)$ right continuous, we deduce from~\eqref{eq:exact-curve}--\eqref{eq:rp-sln-k} that for every $x$ it holds
    \begin{align*}
         \Big(\omega_k(t,\cdot)_x\Big)((-\infty, x])&= \omega_k(t,x)-u^L = \int_{
        {\{\tau\in[0,s]\colon\sigma(\tau) \leqslant x/t\}}}(\tilde r_k\circ\gamma_k)(\varsigma)\, d\varsigma \\
        &= \Big((\tilde r_k\circ \gamma_k)\cdot\mcL\restr {[0,s]}\Big)\Big((t\sigma_k)\inv((-\infty, x])\Big),
    \end{align*}
    that justifies~\eqref{eq:rp-sln-k-pf}.
\end{remark}

\medskip
 
\begin{figure}[h]
    \centering
    \begin{minipage}{0.47\textwidth}
        \centering
        \includegraphics{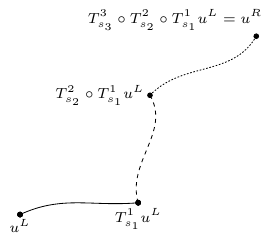}\\
        (a)
    \end{minipage}\hfill
    \begin{minipage}{0.47\textwidth}
        \centering
        \includegraphics{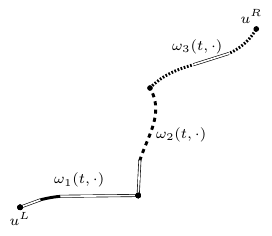}\\
        (b)
    \end{minipage}
    \caption{(a) Joining two states $u^L$, $u^R$ by the solution curves (the $u$-components are drawn) ordered according to the number of family; (b) the corresponding solution~$\omega_k$ to the Riemann problem $(u^L, u^R)$.}
    \label{fig:RP}
\end{figure}

Finally, using Lipschitz version of the Inverse Function Theorem, one proves that
there exists $\drp>0$ such that $\Beone\subset \mcU$ and, 
given $u^L, u^R\in \Beone$, there are unique numbers $s_1, \dots, s_n \in \R$ satisfying 
\begin{equation}
\label{eq:sln-rp}
    u^R = T_{s_n}^n \circ \cdots \circ T_{s_1}^1 u^L
\end{equation}
and
\begin{align}
    \label{eq:tv-along-curves}
    \sum_{k=1}^n |s_k| = \const |u^R-u^L|, \qquad |u^R - u^L| = \const \sum_{k=1}^n |s_k|,
\end{align}
the map 
\begin{equation}
\label{eq:s-to-ur}
    (u^L, u^R) \mapsto (s_1, \dots, s_n)
\end{equation}
being Lipschitz continuous, with $\const$ Lipschitz constant, see Fig.~\ref{fig:RP}(a).

We will use also the notation
\begin{equation*}
    s_k(u^L, u^R) := s_k.
\end{equation*}

By construction, the Riemann problem between ${k}$-th and $(k+1)$-th elements of the sequence
    \begin{align}\label{eq:sequence-Tksk}
        u^L,\qquad T^1_{s_1}u^L,\qquad  T^2_{s_2}\circ T^1_{s_1}u^L,\qquad \ldots,\qquad  \smash{\overbrace{T^n_{s_n}\circ\ldots\circ T^1_{s_1}u^L}^{u^R}}
    \end{align}
is solved by a the map $\omega_k(t,x)$ as in~\eqref{eq:rp-sln-k} (with $k$-th element of the sequence~\eqref{eq:sequence-Tksk} instead of $u^L$).

The solution to $(u^L, u^R)$ is obtained by patching together the~$n$ different solutions~$\omega_k$:
\begin{equation}
\label{eq:rp-patch}
    u(t,x) := \omega_k(t,x) \text{ for } \hat \lambda_{k-1} < \frac{x}{t} < \hat \lambda_k,
\end{equation}
see Fig.~\ref{fig:RP}(b). 
This is possible in view of the separation of characteristic speeds in~\eqref{eq:separation-eigenvalues}. We summarize the above discussion in the following proposition.
\begin{prop}\label{prop:rp}
    There exists a number $\drp>0$ such that the following holds. For every $u^L, u^R \in \Beone$, there exist unique numbers $s_1,\ldots, s_n$ such that~\eqref{eq:sln-rp} holds, and a unique vanishing viscosity solution $u$ to the Riemann problem $(u^L, u^R)$. This solution is given by~\eqref{eq:rp-patch}, where $\omega_k(t,x)$ is a solution~\eqref{eq:rp-sln-k} to the Riemann problem between ${k}$-th and $(k+1)$-th elements of the sequence~\eqref{eq:sequence-Tksk}. Moreover,~\eqref{eq:tv-along-curves} holds, and the map $(u^L, u^R)\mapsto s_k$ is Lipschitz with $\const$ Lipschitz constant.
\end{prop}

Note that although $u_0=u^L\in\Beone$ and $u_n=u^R\in \Beone$, the solution $u$ may take values in $\mcU \setminus \Beone$ as well.

\subsection{Definition of \texorpdfstring{$\Qcost$}{Qcost}}
\label{ss:def-cost-syst}

Consider now two states $u^L, u^R\in\Beone$, so that
\begin{equation*}
    u^R = T_{s_n}^n \circ \cdots \circ T_{s_1}^1 u^L.
\end{equation*}
We can now define the $k$-th cost $\cost_k(u^L, u^R)$ of the jump $(u^L, u^R)$ by 
\begin{equation*}
    \cost_k(u^L, u^R) := \cost_{f_k} (0, s_k),
\end{equation*}
where the $\cost_{f_k}$ on the r.h.s. is the (scalar) cost function introduced in Section~\ref{s:cost-func-and-prop-recall}, applied to the reduced flux $f_k$ (see previous section) on the interval $\I(s_k)$. In view of Remark~\ref{rmk:shifts}, we notice that the definition of $\cost$ is robust we shift the domain of $f_k$ from $\I(s_k)$ to $a+ \I(s_k)$ and if we substitute $f_k$ with $b  + f_k$.
We now define $\Qcost_k,\Qcost\colon\BV(\R,\Beone)\to\R_{\geqslant 0}$ by
\begin{align*}
\Qcost_k(u) := \sum_{\abs{u_x}(\{x\})>0}\cost_k(u_-(x),u_+(x)), \qquad \Qcost(u) :=\sum_{k=1}^n\Qcost_k(u).
\end{align*}

\subsection{Definition of \texorpdfstring{$V_k$}{Vk} and \texorpdfstring{$\Qknown$}{QGBB}}
\label{ss:def-V-k-Qknown}

We recall  now the rigorous definition of the \emph{$k$-total variation along curves} $V_k$ and the \emph{Glimm--Bianchini--Bressan} functional $\Qknown$.

Let $u\in\BV(\R;\Beone)$. Consider first of all the derivative of~$u$ as a vector valued measure $u_x= u_x^\cont + u_x^\atom$, decomposed into its continuous and atomic part. The~$k$-th wave measure $\mu_k = \mu_k^\cont + \mu_k^\atom$ is defined as follows. The continuous part $\mu_k^\cont$ is the Radon measure such that
\begin{equation*}
    \int \varphi(x) \mu_k^\cont(dx) = \int \varphi(x) \scal{\ell_k(u(x))}{u_x^\cont(dx)}
\end{equation*}
for every $\varphi \in C_c(\R;\R)$. The atomic part $\mu_k^\atom$ is concentrated on the (at most countable)  set $\{x_\alpha\}_{\alpha}$ of jumps of~$u$ and it is defined by
\begin{equation*}
    \mu_k^\atom(x_\alpha) := s_k(u_-(x_\alpha), u_+(x_\alpha)).
\end{equation*}
Note that the r.h.s. is well-defined due to $u_{\pm}(x^\alpha)\in \Beone$.
We can now define, for each $k=1,\dots, n$, $V_k,V \colon \BV(\R; \Beone) \to \R_{\geqslant 0}$ by
\begin{equation*}
    V_k(u) :=  |\mu_k|(\R),\qquad V(u):=\sum_{k=1}^n V_k(u).
\end{equation*}
It is not difficult to see
\begin{align}\label{eq:V-equiv-TV}
    V(u) \leqslant\const \TV(u),\qquad \TV(u)\leqslant\const V(u).
\end{align}

We also define $\Qknown \colon \BV(\R; \Beone) \to \R_{\geqslant 0}$ by
\begin{equation*}
    \Qknown(u) := \Qtrans(u) + \frac{1}{4} \sum_{k=1}^n  \Qcubic_k(u). 
\end{equation*}
Here the \emph{transversal} and \emph{cubic} functionals $\Qtrans$,  $\Qcubic_k \colon\allowbreak \BV(\R; \Beone) \to \R_{\geqslant 0}$ are constructed as follows.
\begin{itemize}
    \item  $\Qtrans$ measures the potential transversal  interaction that could take place along the semigroup evolution starting from initial datum $u$:
    \begin{equation*}
    \Qtrans(u):= \sum_{h<k} \big( |\mu_h| \times |\mu_k| \big) (\udhp),\quad \udhp:= \{ (x,y) \in \R^2 \, : \, x>y \}.
\end{equation*}
    \item $\Qcubic_k$ measures the potential interaction between the wavefronts of $k$-th family weighted with their relative speed:
    \begin{equation}
    \label{eq:qcubic}
        \Qcubic_k(u) := \int_0^{|\mu_k|(\R)} \int_{0}^{|\mu_k|(\R)} |\sigma_k(\tau) - \sigma_k(\tau')| d\tau d\tau'.
    \end{equation}
\end{itemize}
Here the \emph{speed} $\sigma_k \colon (0, |\mu_k(\R)) \to \R$ is defined as follows. Fix $\tau \in (0, |\mu_k(\R)|)$, let $x\in\R$ be maximal number satisfying $\abs{\mu_k}((-\infty, x))\leqslant\tau$. Now put
\begin{align}
    \sigma_k(\tau) = 
    \begin{cases}
        \lambda_k(u_-(x)),&\mu_k(\{x\}) = 0,\\
        \frac{d}{d \varsigma}\big[\conv\limits_{[0, \mu_k(\{x\})]}f_k\big](\varsigma)\Restr {\varsigma = \tau - \abs{\mu_k}((-\infty, x))},& \mu_k(\{x\}) > 0,\\
        \frac{d}{d \varsigma}\big[\conc\limits_{[ \mu_k(\{x\}),0]}f_k\big](\varsigma)\Restr {\varsigma = \abs{\mu_k}((-\infty, x))-\tau},& \mu_k(\{x\}) < 0,\\
    \end{cases}\label{eq:abs-sigma-map}
\end{align}
where $f_k$ in the last two cases is the $k$-th reduced flux corresponding to the Riemann problem $(u_-(x),u_+(x))$.

There is a slight abuse of notation, because we use the symbol $\sigma_k$ both for the map defined in~\eqref{eq:abs-sigma-map} and for a component of the curves solving Riemann problems in~\eqref{eq:exact-curve}. No confusion should occur.

Finally observe that
\begin{align}
    \label{eq:qtrans-qcubic-bound-tv}
    \Qtrans(u)  \leqslant \const \TV(u)^2, \qquad   \Qcubic  \leqslant \const \TV(u)^3.
\end{align}

\subsection{Definition of \texorpdfstring{$\Upsilon$}{Upsilon}}
\label{ss:def-Upsilon}

\begin{definition}\label{def:Upsilon}
    We define $\Upsilon \colon \BV(\R; \Beone) \to \R_{\geqslant 0}$ by
    \begin{align}
    \label{eq:upsilon}
        \Upsilon(u) := C_1 V(u)^2 + C_2 \Qknown(u) - \Qcost(u).
    \end{align}
\end{definition}
Being the difference of nonnegative quantities, it is not clear apriori that ${\Upsilon\geqslant0}$. We will prove this fact in Theorem~\ref{thm:quadr-bd}. Constants $C_1, C_2$ will be fixed  in Theorems~\ref{thm:quadr-bd} and~\ref{thm:interaction-estimate}, with $C_1 \gg 1$ and $C_2 - C_1 \gg 1$, in order for the proofs to work. We indicate every new requirement on $C_1$ and $C_2$ by a \indfixC{green box.} In the scalar case it was enough to choose $C_1 = \frac{1}{4} \norm{\derder{f}}_{L^\infty}$ and $C_2 = 0$.

\subsection{A remark on semigroup domain \texorpdfstring{$\mcD$}{D}}
\label{ss:domain-def}
The functional~$\Upsilon$, as well as~$V$ and~$\Qknown$, is defined for every $u \in \BV(\R; \Beone)$. Indeed, in order to define~$\Upsilon(u)$, one just need to know that $u$ takes values in~$\Omega$, it has (at most) countably many jumps with finite total size and that every Riemann problem $(u_-(x), u_+(x))$ can be solved. These three requirements are satisfied if $u \in \BV(\R; \Beone)$. 

The domain $\BV(\R; \Beone)$, however, has also three disadvantages:
\begin{enumerate}[(a)]
    \item \label{pt:large-tv} maps in $\BV(\R; \Beone)$ may have large $\TV$, while, for proving some statements in \ref{thm-main-a} and \ref{thm-main-b}, we need sufficiently small $\TV$;
    \item \label{pt:no-solution} for $\bar u \in \BV(\R; \Beone)$, there might be no solution to \eqref{eq:syst-cons-law}--\eqref{eq:initial-datum}, while, for some statements in \ref{thm-main-a} and \ref{thm-main-b}, we need to know that $S_t \bar u$ exists;
    \item  \label{pt:no-invariance} in addition, even if $S_t \bar u$ exists (e.g. if $\bar u$ is a Riemann problem), it might not belong to $\BV(\R; \Beone)$. 
\end{enumerate}

These issues are conveniently solved by a \emph{semigroup approach} in the spirit of Bressan (see~\cite{Bressan2000} and the reference therein). 

\begin{definition}
\label{def:semigroup}
    A \emph{semigroup domain} is a subset $\mcD \subseteq \BV(\R; \Beone)$ such that
    \begin{enumerate}[(i)]
        \item \label{item:semigroup-delta} there are $\delta_0, \delta_1 \in (0, +\infty]$ with
        $\{u \in \BV(\R; \Beonehalf) \colon \TV(u) < \delta_0 \}    \subseteq \mcD \subseteq  \{u \in \BV(\R; \Beone) \colon \TV(u) < \delta_1 \}$;
        \item \label{item:semigroup-sln}  for every $\bar u \in \mcD$ there is a unique vanishing viscosity solution $u$ to \eqref{eq:syst-cons-law}--\eqref{eq:initial-datum} with $u(t,\cdot) \in \mcD$ for all $t \in [0,+\infty)$. 
    \end{enumerate}
\end{definition}

We define now a family $\{\mcD_\beta\}_{\beta>0}$ of subsets of $\BV(\R; \Beone)$, as follows:
\begin{align}\label{eq:def-D-beta}
    \mcD_\beta := \left\{
        u\in\BV(\R;\Beone)\colon \abs{u(-\infty)}\leqslant\frac\drp2,\ V(u)+C\Qknown(u)\leqslant\beta\right\}.
\end{align}
It turns out that, for \indfixbeta{$\beta \ll 1$ (depending only on the flux function $f$),}  $\mcD_\beta$ is a semigroup domain.  Indeed, it is not difficult to see, by~\eqref{eq:V-equiv-TV} and~\eqref{eq:qtrans-qcubic-bound-tv}, that for every $\beta$ there are $\delta_0(\beta)$, $\delta_1(\beta)$ such that~\eqref{item:semigroup-delta} in Def.~\ref{def:semigroup} holds, and $\delta_0, \delta_1 \to 0$ as $\beta \to 0$.
In particular, by \cite{BianchiniBressan2005}, if $\beta \ll 1$, every $\bar u \in \mcD_\beta$ admits a unique vanishing viscosity solution~$u$. The fact that $u(t,\cdot) \in \mcD_\beta$ for all $t$ follows from the analysis in \cite{Bianchini2003}, and in particular from the fact that $V + C\Qknown$ is a Lyapunov functional, i.e. $t \mapsto V+C\Qknown(u(t,\cdot))$ is non-increasing, if $C \gg 1$.\footnote{The Lyapunov property of $V + C\Qknown$ was proven in \cite{Bianchini2003} for piecewise constant functions. See Remark~\ref{rmk:qknown-lyap}.  For general semigroup trajectories, we present a sketch of the proof in Section~\ref{ss:decr-along-semigroup}. } 

Given that, for $\beta \ll 1$, $\mcD_\beta$ is a semigroup domain and it contains only maps with  total variation of order $\beta$, $\mcD_\beta$ does not have the issues \eqref{pt:large-tv}--\eqref{pt:no-invariance} listed above. In the rest of the paper, we will specify in each statement if it holds for every $u \in \BV(\R; \Beone)$ or only for $u \in \mcD_\beta$, with $\beta$ ``small enough''. We indicate every new smallness requirement on $\beta$ by a \indfixbeta{blue box.} Note that, assuming $\beta \leqslant \const$, the choice of $C_1$, $C_2$ in~\eqref{eq:upsilon} is independent of the choice of $\beta$.

\subsection{Basic estimates on Riemann problems}
\label{ss:basic-rp}

\label{a:int-est}
We conclude this section with a list of some elementary properties of Riemann problems, being used as tools in the following. We begin with a standard linear approximation lemma.
\begin{lemma}
\label{l:length-vs-scalarprod}
    Let $u^L, u^R\in\Beone$, then it holds
    \begin{equation*}
        s_k(u^L, u^R) = \scal{\ell_k(u^L)}{ u^R - u^L} + \const |u^R - u^L|^2. 
    \end{equation*}
\end{lemma}
\begin{proof}
    We have
    \begin{equation*}
        u^R - u^L = \smash{\sum_{h=1}^n} \int_0^{s_h} \tilde r_h (\gamma_h(\varsigma)) d\varsigma.
    \end{equation*}
    Hence
    \begin{equation*}
        \begin{aligned}
        &\scal{\ell_k(u^L)}  {u^R - u^L}
             = \sum_{h=1}^n \int_0^{s_h} \scal{\ell_k(u^L)}{ \tilde r_h (\gamma_h(\varsigma))} d\varsigma \\
             &= \sum_{h=1}^n \Big(\int_0^{s_h} \smash{\underbrace{\scal{\ell_k(u^L)}{r_h(u^L)}}_{\delta_{kh}}} \, d\varsigma + \int_0^{s_h} \scal{\ell_k(u^L)}{ \tilde r_h (\gamma_h(\varsigma))- \smash{\underbrace{r_h(u^L)}_{\mathclap{=\tilde r_h(u^L, 0, \sigma_h(\varsigma)) \text{ by } \eqref{eq:unit-length}}}}  } d\varsigma\Big)\\
             &
             = s_k + \const \sum_{h=1}^n \int_0^{s_h}\Big(\underbrace{|u_h(\varsigma) - u^L| + |v_h(\varsigma)|}_{\const\abs{s_h}\text{ by  \eqref{eq:exact-curve}}} \Big)
             d\varsigma \stackrel{\eqref{eq:tv-along-curves}}= s_k + \const |u^R - u^L|^2,
        \end{aligned}
    \end{equation*}
    which concludes the proof. 
\end{proof}

Then, we recall three results from~\cite{BM-CMP}.
\begin{lemma}[{\cite[Lemma~3.5]{BM-CMP}}, Translation of starting point]\label{lemma:3.5}
    Let $\gamma = \gamma_k(u_0,s) = (u, v_k, \sigma_k)$ and $\gamma' =\gamma_k'(u_0',s) = (u', v_k', \sigma_k') $ and $f_k, f_k'\colon\I(s)\to\R$ denote the corresponding reduced fluxes. Then $\norm{u-u'}_{L^\infty}\leqslant\const\abs{u_0-u_0'}$ and
    \begin{align*}
        \norm{v_k-v_k'}_{L^\infty},\       \norm{\sigma_k-\sigma_k'}_{L^1},
        \
        \norm{\derder {f_k}-\derder{f'_k}}_{L^1}\leqslant\const\abs s\cdot\abs{u_0-u_0'}.
    \end{align*}
\end{lemma}

\begin{lemma}[{\cite[Lemma~3.6]{BM-CMP}}, Change of length]\label{lemma:3.6}
    Let $s's''\geqslant0$ and
    \begin{align*}
        &\gamma = \gamma_k(u_0,s'+s'') = (u, v_k, \sigma_k),\quad \gamma' =\gamma_k(u_0, s')= (u', v_k', \sigma_k'),\\
        &\gamma'' = \gamma_k(u(s'), s''),
    \end{align*}
    with $\gamma,\gamma',\gamma''$ being defined on $\I(s'+s'')$, $\I(s')$, $s'+\I(s'')$, respectively. Denote $f_k$, $f_k'$, $f_k''$ the corresponding reduced fluxes.
    Then the eight quantities
    \begin{align*}
        &\norm{u - u'}_{L^\infty(\I(s'))},\  \norm{v_k - v_k'}_{L^\infty(\I(s'))},\
        \norm{\sigma_k - \sigma_k'}_{L^1(\I(s'))} ,\\
        &\norm{u - u''}_{L^\infty({s'+\I(s'')})},\  \norm{v_k - v_k''}_{L^\infty({s'+\I(s'')})},\
        \norm{\sigma_k - \sigma_k''}_{L^1({s'+\I(s'')})}  ,\\
        &\norm{\derder{f_k}-\derder{f_k'}}_{L^1(\I(s'))},\
        \norm{\derder{f_k}-\derder{f_k''}}_{L^1(s'+\I(s''))}
    \end{align*}
    are bounded by $\const|s'||s''|$.
\end{lemma}

\begin{lemma}[{\cite[Lemma~3.9]{BM-CMP}}, Same family, same sign interaction]\label{lemma:3.9}
    Let $s's''\geqslant 0$, $\gamma' = (u', v_k', \sigma_k') = \gamma_k(u^L, s')$, $u^M = u'(s')$, $\gamma'' = (u'', v_k'', \sigma_k'')= \gamma_k(u^M, s'')$, $\gamma = (u, v_k, \sigma_k) = \gamma_k(u^L,s'+s'') $, with $\gamma', \gamma'', \gamma$ being defined on $\I(s')$, $s'+\I(s'')$, $\I(s'+s'')$, respectively. Denote by $f'_k$, $f''_k$, $f_k$ the reduced fluxes of $ \gamma_k', \gamma_k'',\gamma_k $, respectively, assuming $f''_k(s') = f'_k(s')$. Then it follows
    \begin{align*}
        &\norm{u'\cup u''- u}_{L^\infty(\I(s'+s''))}, \
        \norm{v_k'\cup v_k''- v_k}_{L^\infty(\I(s'+s''))}, \
        \norm{\sigma_k'\cup \sigma_k''- \sigma_k}_{L^1(\I(s'+s''))},\\
        &\norm{\derder{f_k'\cup f_k''} - \derder{f_k}}_{L^1(\I(s'+s''))}\leqslant\const
        \begin{cases}
        \Big[f'_k- \conv\limits_{[0,s'+s'']}(f_k'\cup f_k'')\Big](s'),& s', s''\geqslant 0,\\
        \Big[\conc\limits_{[s'+s'',0]}(f_k'\cup f_k'')-\tilde f'_k\Big](s'),& s', s''\leqslant 0.
        \end{cases}
    \end{align*}
\end{lemma}

\section{Nonnegativity and quadratic bound}
\label{s:system-quadratic-bound}

We now state  and prove rigorously \ref{thm-main-a}, item~\eqref{thm-main:qb}. 

\begin{theorem}[Nonnegativity and quadratic bound]
\label{thm:quadr-bd}
    For all $u \in \BV(\R; \Beone)$,
    \begin{equation*}
        0 \leqslant \Upsilon(u) \leqslant \const \TV(u)^2.
    \end{equation*}
\end{theorem}

\begin{proof}
The functional $V$ is known to be equivalent to total variation.
Then, since $\Qknown(u) = \const\TV(u)^2$ and the  $\cost$ terms are non-negative, we immediately get $\Upsilon(u) \leqslant \const\TV(u)^2 $.

Next, applying Prop.~\ref{prop:cost_Linfty_est-recall} to each $\cost$ term in~\eqref{eq:upsilon}, we obtain
    \begin{align*}
        \Upsilon(u) &\geqslant C_1 V(u)^2-\frac 14\sup_{k, f_k}\norm{\derder {f_k}}_{L^\infty} \sum_{k=1}^n \ \sum_{\substack{ |u_x|(\{x\})>0}\hspace{1mm}} s_k(u_-(x),u_{+}(x))^2\\
        &\geqslant C_1
        V(u)^2-\const\TV(u)^2,
    \end{align*}
    where the supremum is taken over all possible reduced fluxes  in the solutions of Riemann problems in our framework. Hence, by~\eqref{eq:V-equiv-TV}, we get $\Upsilon(u)\geqslant 0$ for all $u\in \BV(\R; \Beone)$, \indfixC{provided that $C_1$ is large enough.}\label{C1-C2-fix-1}
\end{proof}

\section{Interaction estimate}
\label{s:interaction-estimate}

We  now state and prove the interaction estimate claimed in \ref{thm-main-a}, item~\eqref{thm-main:ie}. First of all we give the rigorous definition of \emph{change of speed} of wavefronts involved in an interaction.

\begin{definition}[Change of speed]
\label{def:change-speed}
Let $u^L, u^M, u^R\in\Beone$. Put $s_k' = s_k(u^L, u^M)$, $s_k'' = s_k(u^M, u^R)$, $s_k = s_k(u^L, u^R)$, and let the families of curves
\begin{align*}
    &\{\gamma'_k = (u'_k, v'_k,\sigma'_k)\colon \I(s_k')\to\mcU_k\}_{k=1}^n,\quad
    \{\gamma''_k = (u''_k, v''_k,\sigma''_k)\colon s_k'+\I(s_k'')\to \mcU_k\}_{k=1}^n,\\
    &\{\gamma_k = (u_k, v_k,\sigma_k)\colon \I(s_k)\to \mcU_k\}_{k=1}^n
\end{align*}
solve the Riemann problems $(u^L, u^M)$, $(u^M, u^R)$, $(u^L, u^R)$ respectively.  Note that we assume $\gamma''_k$ to be defined on the shifted domain $s_k'+\I(s_k'')$.
The \emph{change of speed} is defined  as:
\begin{align*}
    \Delta\sigma(u^L, u^M, u^R) := \begin{cases}
    \norm{\big(\sigma_k'\cup\sigma_k''\big) - \sigma_k}_{L^1(\I(s_k'+s_k'')\cap \I(s_k))},& s_k's_k''\geqslant 0,\\
    \norm{\big(\sigma_k'\xor\sigma_k''\big) - \sigma_k}_{L^1(\I(s_k'+s_k'')\cap \I(s_k))},& s_k's_k''< 0.
    \end{cases}
\end{align*}
{(Recall the notation $\sigma_k'\cup\sigma_k''$,  $\sigma'\xor\sigma_k''$ from~\eqref{eq:cup-def}--\eqref{eq:xor-def}.)}
\end{definition}

\begin{remark}
Def.~\ref{def:change-speed} corresponds to the one in~\cite[Section 1.1]{BM-CMP}. This is relevant, because we will use many results from \cite{BM-CMP} as a black box.    
\end{remark}

\noindent The following theorem is the rigorous version of  \ref{thm-main-a}, item~\eqref{thm-main:ie}.

\begin{theorem}[Interaction estimate]
    \label{thm:interaction-estimate}
    Let $\beta \ll 1$. Let $u, \hat u \in \mcD_\beta$.  Let $I \subseteq \R$ be an open interval. Assume that
    \begin{enumerate}[(a)]
        \item $u(x) = \hat u(x)$ for $x \notin I$,
        \item on $I$ both $u, \hat u$ are piecewise constant,
        \item on $I$, $u$  has exactly two jumps from $u^L$ to $u^M$ and from $u^M$ to $u^R$ respectively, whilst $\hat u$  has exactly one jump from $u^L$ to $u^R$.
    \end{enumerate}
    Then 
    \begin{equation}
    \label{eq:int-est-main}
        \Delta \sigma (u^L, u^M, u^R) \leqslant \const\Big( \Upsilon (u) - \Upsilon(\hat u)\Big).
    \end{equation}
    \begin{figure}[h]
        \centering
        \includegraphics{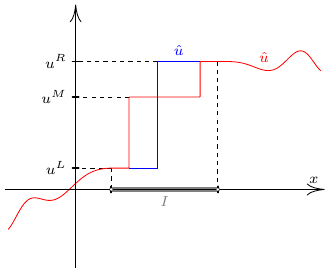}
        \caption{To Theorem~\ref{thm:interaction-estimate}.}
        \label{fig:thm:interaction-estimate}
    \end{figure}
\end{theorem}

To prove Theorem~\ref{thm:interaction-estimate} we follow the approach of~\cite{BM-CMP} (which was inspired by~\cite{Bianchini2003}). We decompose the \emph{interaction} between the two Riemann problems $(u^L, u^M)$,\linebreak $(u^M, u^R)$ into simpler steps, namely:
\begin{enumerate}
    \item \emph{transveral interactions} (between wavefronts belonging to different families),
    \item \emph{same family interactions} (these are analog to interactions in the scalar case, where only one family is present) 
    \item  \emph{adjustment of the length of the curves} (due to the fact that $s_k \not = s_{k'} + s_{k''}$ in general). 
\end{enumerate}
One can quantify how strong each of these simpler interactions is, by means of suitable \emph{amounts of interactions} (a.o.i.), see the specific definitions below in Section~\ref{ss:a-o-i}: transversal interactions are measured by the \emph{transversal}  a.o.i. $\Atrans(u^L, u^M, u^R)$; same family interactions are measured by the \emph{quadratic} a.o.i. $\Aquadr_k(u^L, u^M, u^R)$ or by the \emph{amount of cancellation} $\Acanc_k(u^L, u^M, u^R)$ (respectively if the two interacting wavefronts of the same family have the same sign or opposite sign); the adjustments of the length are measured by the  sum of \emph{transversal} and \emph{cubic} a.o.i., $\Atrans(u^L, u^M, u^R)+\sum_k\Acubic_k(u^L, u^M, u^R)$.

The proof of Theorem~\ref{thm:interaction-estimate} is split in two parts: first one bounds $\Delta \sigma$ by the amount of interactions.

\begin{prop}
\label{prop:ie-intro-delta-sigma}
For $u^L, u^M, u^R \in \Beone$, 
\begin{align}
\label{eq:ie-intro-delta-sigma}
    &\Delta \sigma(u^L, u^M, u^R) \nonumber\\
    &\leqslant \const \bigg[\Atrans 
              + \sum_{k=1}^n \Big( \Acubic_k + \Aquadr_k  +  (|s_k'|+|s_k''|) \Acanc_k \Big) \bigg](u^L, u^M, u^R).
\end{align}
\end{prop}

Then one proves that the amounts of interactions are bounded by the decrease of $\Upsilon$.

\begin{prop}
    \label{prop:ie-intro-delta-upsilon}
For $u, \hat u$ as in Theorem~\ref{thm:interaction-estimate}
\begin{align}
\label{eq:ie-intro-delta-upsilon}
    &\bigg[\Atrans 
              + \sum_{k=1}^n \Big( \Acubic_k + \Aquadr_k  +  (|s_k'|+|s_k''|) \Acanc_k \Big) \bigg](u^L, u^M, u^R)
    \nonumber\\
    &\leqslant 2\Big( \Upsilon(u) - \Upsilon(\hat u)\Big).
\end{align}
\end{prop}

Prop.~\ref{prop:ie-intro-delta-sigma} is essentially (up to minor modification) \cite[Theorem~3.3]{BM-CMP}. On the contrary, Prop.~\ref{prop:ie-intro-delta-upsilon}  is heavily based on the properties of $\Upsilon$ and $\cost$ (see in particular the proof of Prop.~\ref{prop:var-cost}).

\begin{proof}[Proof of Theorem~\ref{thm:interaction-estimate}]
    Combine Prop.~\ref{prop:ie-intro-delta-sigma}--\ref{prop:ie-intro-delta-upsilon}. 
\end{proof}

In Section~\ref{ss:system-two-rp} we rigorously describe the decompositions of the interaction into simpler steps. In Section~\ref{ss:a-o-i} we recall the definition of the various amounts of interactions. In Section~\ref{ss:analysis-element-inter}, we recall some estimates for elementary interaction from \cite{BM-CMP}. In Section~\ref{ss:delta-sigma-bound} we recall 
the proof of Prop.~\ref{prop:ie-intro-delta-sigma} from \cite{BM-CMP} (with a minor improvement). Finally, in Section~\ref{ss:delta-upsilon-bound} we prove
Prop.~\ref{prop:ie-intro-delta-upsilon}: this proof heavily relies on the definition of $\Upsilon$ and in particular on the properties of $\cost$, and  it is the most important contribution of this section (and of the whole paper). We will highlight in Section~\ref{ss:delta-upsilon-strategy-scalar}  how the scalar proof of the analog of Prop.~\ref{prop:ie-intro-delta-upsilon} in~\cite{BaitMod1}  can be adapted to the case of systems.

\subsection{Decomposition of interactions into elementary interactions}
\label{ss:system-two-rp}

We now describe how the \emph{interaction} between two Riemann problems can be decomposed in simpler steps. 

Consider two contiguous Riemann problem (situation \emph{before} the interaction)
\begin{equation}
\label{E_two_riemann_1}
u^M = T^{n}_{s_n'} \circ \dots \circ T^{1}_{s_1'} u^{L}, \qquad  
u^{R} = T^{n}_{s_n''} \circ \dots \circ T^{1}_{s_1''} u^{M},
\end{equation}
and the Riemann problem obtained joining them (situation \emph{after} the interaction)
\begin{equation*}
u^{R} = T^{n}_{s_n} \circ \dots \circ T^{1}_{s_1} u^{L}.
\end{equation*}
In particular the exact curve before the interaction are denoted by $\gamma_k', \gamma_k''$:
\begin{equation*}
    \begin{aligned}
        \gamma_1' & = (u_1', v_1', \sigma_1') := \gamma_1(u^L, s_1'), & \gamma_k' & = (u_k', v_k', \sigma_k') := \gamma_k\big(u_{k-1}'(s_{k-1}'), s_k'\big), & k & \geqslant 2, \\
        \gamma_1'' & =  (u_1'', v_1'', \sigma_1'') := \gamma_1(u^M, s_1''), & \gamma_k'' & = (u_k'', v_k'', \sigma_k'')  := \gamma_k\big(u_{k-1}''(s_{k-1}''), s_k''\big), & k &\geqslant 2,
    \end{aligned}
\end{equation*}
while the exact curve after the interaction are denoted by $\gamma_k$:
\begin{equation*}
    \begin{aligned}
        \gamma_1 & = (u_1, v_1, \sigma_1) := \gamma_1(u^L, s_1), & \gamma_k & = (u_k, v_k, \sigma_k) := \gamma_k\big(u_{k-1}(s_{k-1}), s_k\big), & k & \geqslant 2.
    \end{aligned}
\end{equation*}

We will similarly denote by $f_k', f_k'', f_k$ the corresponding reduced fluxes associated to $\gamma_k', \gamma_k'', \gamma_k$ respectively. It is convenient to assume (w.l.o.g.) that $\gamma_k''$ and $f_k''$ are defined on $s_k' + \I(s_k'')$, instead of $\I(s_k'')$ and $f_k''(s_k') = f_k'(s_k')$. Notice that, if $s_k's_k'' \geqslant 0$, $f_k' \cup f_k''$ is thus continuous (but, in general, it is not differentiable at point $s_k'$, see first and fourth graphs in Fig.~\ref{fig:cubic}). 

\begin{figure}[h]
\includegraphics{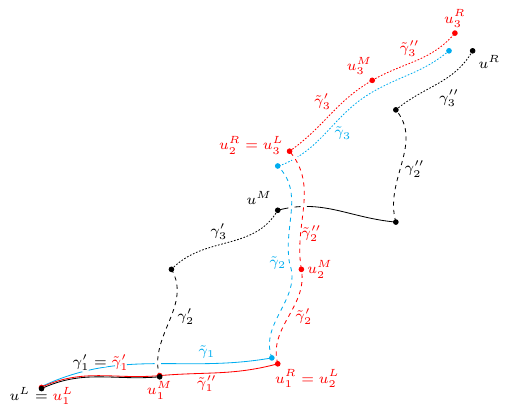}
\caption{Elementary curves: two Riemann problems before interactions (black); after transversal interactions (red); after same family interactions (cyan). Only $u$-components are drawn.}
\label{F_el_curves_after_trans}
\end{figure}
It is also useful to consider an intermediate situation where, roughly speaking, all \emph{transversal interactions} have already taken place. Specifically, fix a family $k \in \{1, \ldots, n\}$ and consider the points (see Fig.~\ref{F_el_curves_after_trans})
\begin{equation*}
    \begin{aligned}
        u^L_1 & := u^L,  &  u^L_k & := T^{k-1}_{s''_{k-1}} \circ T^{k-1}_{s'_{k-1}} \circ \dots \circ T^1_{s''_1} \circ T^1_{s'_1} u^L,  & k & = 2,\dots, n, \\
        u^M_k & := T^k_{s'_k} u^L_k,  & u^R_k & := T^k_{s''_k} u^M_k, & k & = 1,\dots,n.
    \end{aligned}
\end{equation*}
By definition, the Riemann problem between $u^L_k$ and $u^M_k$ is solved by an elementary wave of the $k$-th family with strength $s_k'$ and the Riemann problem between $u^M_k$ and $u^R_k$ is solved by an elementary wave of the $k$-th family with strength $s''_k$. 

We denote by $\tilde \gamma_k' = (\tilde u_k', \tilde v_k', \tilde \sigma_k')$ the curve which solves the Riemann problem $(u^L_k, u^M_k)$ and by $\tilde f_k'$ the associated reduced flux. Similarly, let $\tilde \gamma_k'' = (\tilde u_k'', \tilde v_k'', \tilde \sigma_k'')$ be the curve solving the Riemann problem $(u^M_k, u^R_k)$ and let $\tilde f_k''$ be the associated reduced flux.

Clearly, $\tilde \gamma_k'$, $\tilde f_k'$ are defined on $\I(s_k')$, whereas it is convenient, as before, to assume (w.l.o.g.) that $\tilde \gamma_k''$ and $\tilde f_k''$ are defined on $s_k' + \I(s_k'')$ (instead of $\I(s_k'')$) and that $\tilde f_k''(s_k') = \tilde f_k'(s_k')$. Notice that, if $s_k' s_k'' \geqslant 0$, $\tilde f_k' \cup \tilde f_k''$ is continuous and, in contrast with $f_k' \cup f_k''$, it is also differentiable at $s_k'$ (with derivative given by $\lambda_k(u_k^M)$), so that, overall, $\tilde f_k' \cup \tilde f_k'' \in C^{1,1}(\I(s_k' + s_k''))$.

Finally, we introduce the curves resulting from the \emph{same family interactions},
\begin{align*}
    &\tilde \gamma_1 = (\tilde u_1,\tilde v_1,\tilde\sigma_1) :=\gamma_1(u^L, s_1'+s_1''),\\
    &\tilde\gamma_k = (\tilde u_k,\tilde v_k,\tilde\sigma_k) := \gamma_k(\tilde u_{k-1}(s'_{k-1}+s''_{k-1}), s_k'+s_k''), k = 2,\ldots, n,
\end{align*}
see Fig.~\ref{F_el_curves_after_trans}.

\subsection{Amounts of interaction}
\label{ss:a-o-i}
\begin{figure}[h]
\includegraphics{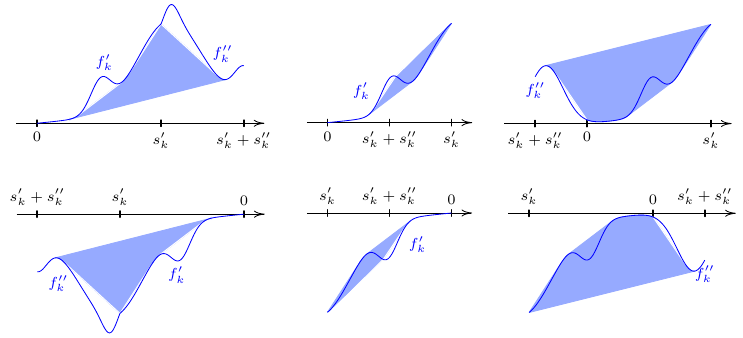}
\caption{Cubic amount of interaction is represented by the violet-blue area.}
\label{fig:cubic}
\end{figure}

We now recall the definition of four quantities, called \emph{amounts of interaction} introduced in~\cite{BM-CMP, Bianchini2003}. The \emph{transversal amount of interaction} and the \emph{amount of cancellation} are given by
\begin{equation}\label{eq:Atrans-Acanc-def}
\begin{aligned}
    \Atrans(u^L, u^M, u^R) & := \sum_{h<k} |s_k'||s_h''|, \\
    \Acanc_k(u^L, u^M, u^R) & := 
    \begin{cases}
        0, & \text{if } s_k' s_k'' \geqslant 0, \\ 
        \min \{ |s_k'|, |s_k''|\}, & \text{if } s_k' s_k'' \leqslant 0. \\
    \end{cases}
\end{aligned}
\end{equation}
Next, we reproduce the definition of \emph{cubic amount of interaction} from~\cite{Bianchini2003},
\begin{align*}
    &\Acubic_k(u^L, u^M, u^R) := \\
    &\begin{cases}\int_0^{s_k'+s_k''}  \Big[\conv\limits_{[0,s_k']} f_k' \cup \conv\limits_{[s_k', s_k'+s_k'']} f_k'' - \conv\limits_{[0, s_k'+s_k'']} \big(f_k' \cup f_k'' \big)\Big] \, d\tau,& s_k', s_k''  \geqslant 0, \\
         \int_0^{s_k'} \Big[\conv\limits_{[0,s_k'+s_k'']} f_k' \cup \conc\limits_{[s_k'+s_k'', s_k']} f_k'  - \conv\limits_{[0,s_k']} f_k'\Big]\, d \tau,& -s_k'\leqslant s_k''\leqslant 0,\\
        \int_{s_k'+s_k''}^{s_k'}\, 
         \Big[\conc\limits_{[s_k'+s_k'', s_k']}f''_k-
            \conc\limits_{[s_k'+s_k'', 0]} f_k''\cup \conv\limits_{[0,s_k']}f_k''\Big]\,d\tau, & s_k''\leqslant-s_k'\leqslant 0,\\
            \int_{s_k'+s_k''}^0  \Big[  \conc\limits_{[ s_k'+s_k'',0]} \big(f_k'' \cup f_k' \big)-\conc\limits_{[s_k'+s_k'',s_k']} f_k''\cup \conc\limits_{[s_k',0]} f_k' \Big] \, d\tau,& s_k', s_k''  \leqslant 0, \\
         \int_{s_k'}^0 \Big[\conc\limits_{[s_k',0]} f_k'-\conv\limits_{[s_k',s_k'+s_k'']} f_k'\cup \conc\limits_{[s_k'+s_k'',0]} f_k'  \Big]\, d \tau,& -s_k'\geqslant s_k''\geqslant 0,\\
        \int_{s_k'}^{s_k'+s_k''}\, 
         \Big[\conc\limits_{[s_k',0]}f_k''\cup\conv\limits_{[0,s_k'+s_k'']} f_k''-\conv\limits_{[s_k'+s_k'', s_k']}f''_k\Big]\,d\tau, & s_k''\geqslant-s_k'\geqslant 0,
    \end{cases}
\end{align*}
see Fig.~\ref{fig:cubic}.

The \emph{quadratic amount of interaction} is given by
\begin{equation*}
    \Aquadr_k(u^L, u^M, u^R)  := \begin{cases}
        \tilde f'_k(s_k') - \Big[\conv\limits_{[0,s'_k+s''_k]}(\tilde f_k'\cup \tilde f_k'')\Big](s_k'),& s_k', s_k''\geqslant 0,\\
        \Big[\conc\limits_{[s'_k+s''_k,0]}(\tilde f_k'\cup \tilde f_k'')\Big](s_k')-\tilde f'_k(s_k'),& s_k', s_k''\leqslant 0,\\
        0, & s_k's_k''\leqslant 0.
    \end{cases}
\end{equation*}
Observe that for $s_k',s_k''\geqslant0$, $\Aquadr_k = \Height_{\tilde f_k'\cup \tilde f_k''} (0,s_k',s_k'+s_k'')$.

We conclude this section by recalling a seminal theorem from~\cite{Bianchini2003}. In view of the following usage, we also add a quadratic bound (the second and third inequalities).
\begin{prop}[{\cite[Theorem~3.7]{Bianchini2003}}]\label{prop:sdef-quadr-bound}
In the notation of Section~\ref{ss:a-o-i}, it holds
\begin{align*}
    \sum_{k=1}^n\abs{s_k - (s_k'+s_k'')}&\leqslant\const\Big[
    \Atrans(u^L, u^M , u^R)+\sum_{k=1}^n\Acubic_k(u^L, u^M , u^R)
    \Big]\\ &\leqslant\const\left[\sum_{k=1}^n|s_k'|\right]\cdot\left[\sum_{k=1}^n|s_k''|\right]\leqslant\const\abs{u^M-u^L}\abs{u^R-u^M}.
\end{align*}
\end{prop}
\begin{proof}
For the first inequality, we refer to~\cite[Theorem~3.7]{Bianchini2003}. The quadratic bound in the second line follows from the definitions of $\Atrans$ and $\Acubic_k$.
\end{proof}

\subsection{Analysis of elementary interactions}
\label{ss:analysis-element-inter}

Before proving Prop.~\ref{prop:ie-intro-delta-sigma}--\ref{prop:ie-intro-delta-upsilon}, we recall some elementary estimates from \cite{BM-CMP}. 
Lemmas~\ref{lemma:310} and~\ref{lemma:313} are copied and pasted from \cite[Lemma~3.10, Corollary~3.13]{BM-CMP}, respectively. Lemma~\ref{l:improved-311} is a slightly improved version  of \cite[Lemma~3.11]{BM-CMP}, namely one should add the multiplier $|s_k'| + |s_k''|$ to the r.h.s. in~\cite[Lemma~3.11]{BM-CMP}. We will prove only Lemma~\ref{l:improved-311}. 

\begin{lemma}[{\cite[Lemma~3.10]{BM-CMP}}]
    \label{lemma:310}

    For any $k=1,\ldots, n$, the quantities
    \begin{align*}
        &\norm{u_k'-\tilde u_k'}_{L^\infty(I(s_k'))},\ 
        \norm{u_k''-\tilde u_k''}_{L^\infty(s_k'+I(s_k''))},\
        \norm{\sigma_k'-\tilde \sigma_k'}_{L^1(I(s_k'))},\\ 
        &\norm{\sigma_k''-\tilde \sigma_k''}_{L^1(s_k'+I(s_k''))},\
        \norm{\derder{f_k'} - \derder{\tilde f_k'}}_{L^1(\I(s_k'))}, \
        \norm{\derder{f_k''} - \derder{\tilde f_k''}}_{L^1(s_k'+\I(s_k''))}
    \end{align*}
    are bounded by $\const\Atrans(u^L, u^M, u^R)$.
\end{lemma}
\begin{lemma}[{Refined version of \cite[Lemma~3.11]{BM-CMP}}]
    \label{l:improved-311}
    For any $k =1,\ldots, n$,
    \begin{itemize}
        \item if $s_k's_k''\geqslant 0$, then
        \begin{align*}
        &\begin{rcases}
            \norm{(\tilde u_k'\cup\tilde u_k'')-\tilde u_k}_{L^\infty(\I(s_k'+s_k''))}\\
            \norm{(\tilde \sigma_k'\cup\tilde \sigma_k'')-\tilde \sigma_k}_{L^1(\I(s_k'+s_k''))}\\
            \norm{\derder{\tilde f'_k} \cup\derder{\tilde f''_k} - \derder{\tilde f_k}}_{L^1(\I(s_k'+s_k''))}
        \end{rcases}  
        \\
        &\qquad\leqslant
        \const \sum_{h=1}^k\Big[\Aquadr_h(u^L, u^M , u^R)
        +
        (|s_h'|+|s_h''|)
        \Acanc_h(u^L, u^M , u^R)
        \Big];
        \end{align*}
        \item if $s_k's_k''< 0$, then
        \begin{align*}
        &\begin{rcases}
            \norm{(\tilde u_k'\xor\tilde u_k'')-\tilde u_k}_{L^\infty(\I(s_k'+s_k''))}\\
            \norm{(\tilde \sigma_k'\xor\tilde \sigma_k'')-\tilde \sigma_k}_{L^1(\I(s_k'+s_k''))}\\
            \norm{\derder{\tilde f'_k} \xor\derder{\tilde f''_k} - \derder{\tilde f_k}}_{L^1(\I(s_k'+s_k''))}
        \end{rcases}  
        \\
        &\qquad\leqslant
        \const \sum_{h=1}^k\Big[\Aquadr_h(u^L, u^M , u^R)
        +
        (|s_h'|+|s_h''|)
        \Acanc_h(u^L, u^M , u^R)
        \Big].
        \end{align*}
    \end{itemize}
\end{lemma}

\begin{lemma}[{\cite[Corollary~3.13]{BM-CMP}}]
    \label{lemma:313}
    For any $k =1,\ldots, n$, it holds
    \begin{align*}
        &\norm{\tilde u_k - u_k}_{L^\infty(\I(s_k'+s_k'')\cap\I(s_k))},\
        \norm{\tilde \sigma_k - \sigma_k}_{L^1(\I(s_k'+s_k'')\cap\I(s_k))},\\
        &
        \hspace{2.2cm}\mathclap{\norm{\derder{\tilde f_k} - \derder{f_k}}_{L^1(\I(s_k'+s_k'')\cap\I(s_k))}}\hspace{1.8cm}\leqslant\const\Big[\Atrans(u^L, u^M, u^R)+\sum_{h=1}^n\Acubic_h(u^L, u^M, u^R)\Big].
    \end{align*}
\end{lemma}

\begin{proof}[Proof of Lemma~\ref{l:improved-311}]
The proof below is reminiscent to the one of~\cite[Lemma~3.11]{BM-CMP}, the only difference being in the analysis of the interaction between $k$-th waves having opposite sign, $s_k' s_k'' \leqslant 0$, and precisely in estimates \eqref{l:improved-311:eq1}--\eqref{l:improved-311:eq2}, \eqref{l:improved-311:eq4}--\eqref{l:improved-311:eq5}. Following~\cite{BM-CMP}, we divide the proof into two steps.
Throughout the proof, we omit the arguments $(u^L, u^M, u^R)$ of the amounts of interaction.
   
    \textit{Step 1.} We claim that for each $k = 1,\ldots, n$ it holds
    \begin{align}
        \abs{\tilde u''_k(s_k'+s_k'')-\tilde u_k(s_k'+s_k'')}\leqslant\const \sum_{h=1}^k\Big[&\Aquadr_h+(|s_h'|+|s_h''|)
        \Acanc_h\Big].\label{l:improved-311:eq0.1}
    \end{align}
    We proceed by induction on $k$.
    Both the base and the inductive step reduce to proving that for every $k$,
    \begin{align}
        &\abs{\tilde u''_k(s_k'+s_k'')-\tilde u_k(s_k'+s_k'')} \nonumber\\
        &\qquad\leqslant\const\Big[
        \Aquadr_k+(|s_k'|+|s_k''|)
        \Acanc_k +\abs{\tilde u_k'(0)-\tilde u_k(0)}
        \Big],\label{l:improved-311:eq0.2}
    \end{align}
    being l.h.s. equal to the third term in the r.h.s. for the subsequent $k$.
    We distinguish three cases.
    \begin{itemize}
        \item \textit{Case $s_k's_k''\geqslant 0$.} Set $\hat \gamma_k = (\hat u_k,\hat v_k,\hat\sigma_k):= \gamma_k(\tilde u_k'(0), s_k'+s_k'')$, see Fig.~\ref{fig:improved-311}(a). We have
        \begin{align}
            \abs{\tilde u''_k(s_k'+s_k'')-\hat u_k(s_k'+s_k'')}\leqslant\const\Aquadr_k\label{l:improved-311:eq0.4}
        \end{align}
        by Lemma~\ref{lemma:3.9}, and
        \begin{align}
            \abs{\hat u_k(s_k'+s_k'')-\tilde u_k(s_k'+s_k'')}\leqslant \const\abs{\tilde u_k'(0)-\tilde u_k(0)} \label{l:improved-311:eq0.5}
        \end{align}
    by Lemma~\ref{lemma:3.5}. From~\eqref{l:improved-311:eq0.4} and~\eqref{l:improved-311:eq0.5} we get~\eqref{l:improved-311:eq0.2}.

\begin{figure}[h]
    \begin{minipage}{0.32\textwidth}
    \centering
        \includegraphics{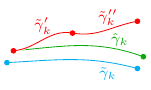}
        \\
        (a)
    \end{minipage}
    \hfill
    \begin{minipage}{0.32\textwidth}
    \centering
        \includegraphics{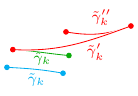}
        \\
        (b)
    \end{minipage}
    \hfill
    \begin{minipage}{0.32\textwidth}
    \centering
        \includegraphics{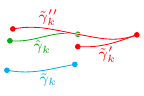}
        \\
        (c)
    \end{minipage}
    \caption{To the proof of Lemma~\ref{l:improved-311}: (a) case $s_k's_k''\geqslant 0$; (b) case $s_k's_k''< 0$, $\abs{s_k'}\geqslant\abs{s_k''}$;  (c) case $s_k's_k''< 0$, $\abs{s_k'}<\abs{s_k''}$; $u$-components of curves are drawn.}
    \label{fig:improved-311}
    \end{figure}
        
        \item \textit{Case $s_k's_k''<0$, $\abs{s_k'}\geqslant\abs{s_k''}$.} Set $\hat \gamma_k := (\hat u_k,\hat v_k,\hat\sigma_k)= \gamma_k(\tilde u_k'(0), s_k'+s_k'')$, see Fig.~\ref{fig:improved-311}(b).
        We estimate
        \begin{align}
            \abs{\tilde u''_k(s_k'+s_k'')-\tilde u'_k(s_k'+s_k'')}&\leqslant \const\abs{s_k''}^2 = \const\abs{s_k''}\Acanc_k\label{l:improved-311:eq1}
        \end{align}
    by the tangency of the $u$-components of the curves $\tilde\gamma'_k$ and $\tilde\gamma''_k$ at $\tilde u_k'(s_k')$ and the fact that such two curves have uniformly bounded second derivative,
    \begin{align}
        \abs{\tilde u'_k(s_k'+s_k'') - \hat u_k(s_k'+s_k'')}\leqslant \const\abs{s_k'+s_k''}\cdot\abs{s_k''} \leqslant\const\abs{s_k'}\cdot\Acanc_k\label{l:improved-311:eq2}
    \end{align}
    by Lemma~\ref{lemma:3.6}, and
    \begin{align}
        \abs{\hat u_k(s_k'+s_k'')-\tilde u_k(s_k'+s_k'')} \leqslant \const\abs{\tilde u_k'(0)-\tilde u_k(0)}
        \label{l:improved-311:eq3}
    \end{align} by Lemma~\ref{lemma:3.5}.
    From~\eqref{l:improved-311:eq1}--\eqref{l:improved-311:eq3} we get~\eqref{l:improved-311:eq0.2}.
    \item \textit{Case $s_k's_k''<0$, $\abs{s_k'}<\abs{s_k''}$.} Set $\hat \gamma_k = (\hat u_k,\hat v_k,\hat\sigma_k):= \gamma_k(\tilde u_k''(0), s_k'+s_k'')$, see Fig.~\ref{fig:improved-311}(c). We estimate
    \begin{align}
        \abs{\tilde u''_k(s_k'+s_k'') - \hat u_k(s_k'+s_k'')}\leqslant \const\abs{s_k'+s_k''}\cdot\abs{s_k'} \leqslant\const\abs{s_k''}\Acanc_k\label{l:improved-311:eq4}
    \end{align}
    by Lemma~\ref{lemma:3.6}, and
    \begin{align}
        &\abs{\hat u_k(s_k'+s_k'')-\tilde u_k(s_k'+s_k'')} \leqslant \const\Big(\abs{\hat u_k(0)-\tilde u_k'(0)}+\abs{\tilde u_k'(0) - \tilde u_k(0)}\Big)\nonumber \\
        &\leqslant\const\Big(\abs{s_k'}^2+\abs{\tilde u_k'(0) - \tilde u_k(0)}\Big)=\const\Big(\abs{s_k'}\Acanc_k+\abs{\tilde u_k'(0) - \tilde u_k(0)}\Big)
        \label{l:improved-311:eq5}
    \end{align} by Lemma~\ref{lemma:3.5} and using tangency of the $u$-components of the curves $\tilde\gamma'_k$, $\tilde\gamma''_k$ at $\tilde u _k'(s_k')$. From~\eqref{l:improved-311:eq4} and~\eqref{l:improved-311:eq5} we get~\eqref{l:improved-311:eq0.2}.
    \end{itemize}

    \textit{Step 2} is completely analogous to the one in~\cite{BM-CMP}, with the same modification as in Step 1.
       \end{proof}
       
\subsection{Bound on \texorpdfstring{$\Delta \sigma$}{delta sigma}}
\label{ss:delta-sigma-bound}
We are now ready to prove Prop.~\ref{prop:ie-intro-delta-sigma}. 
We have already observed that  this proposition is a refined version of~\cite[Theorem~3.3]{BM-CMP}. Indeed, comparing to~\cite[Theorem~3.3]{BM-CMP}, an additional multiplier $|s_k'|+|s_k''|$ appears before $\Acanc_k(u^L,\allowbreak u^M,\allowbreak u^R)$ in the r.h.s. The presence of this multiplier shows in particular that $\Delta \sigma(u^L, u^M, u^R) = \const\abs{u^R-u^M}\cdot\abs{u^M-u^L} $ i.e. it is indeed a ``quadratic quantity'', as one expects from being bounded by the decrease of the quadratic functional $\Upsilon$.

\begin{proof}[Proof of Prop.~\ref{prop:ie-intro-delta-sigma}]
    Combine Lemmas~\ref{lemma:310},~\ref{l:improved-311} and~\ref{lemma:313}.
\end{proof}
\begin{remark}
     The derivation of Prop.~\ref{prop:ie-intro-delta-sigma} from  Lemmas~\ref{lemma:310},~\ref{l:improved-311},~\ref{lemma:313} is reminiscent to the one of \cite[Theorem~3.3]{BM-CMP} from \cite[Lemmas~3.10--3.11, Corollary~3.13]{BM-CMP} respectively. The difference is only in the additional multiplier $|s_k'|+|s_k''|$ in Lemma~\ref{l:improved-311} (w.r.t. \cite[Lemma~3.11]{BM-CMP}).
\end{remark}

\subsection{The strategy for \texorpdfstring{$\Delta \Upsilon$}{Delta Upsilon}  in the scalar case}
\label{ss:delta-upsilon-strategy-scalar}

Before proving Prop.~\ref{prop:ie-intro-delta-upsilon}, for which we will exploit the properties of $\Upsilon$ and $\cost$, it is instructive to make some analogy to the scalar case in \cite[Section 4]{BaitMod1}. 

First of all, observe that the amounts of interactions (in particular $\Aquadr(u^L,\allowbreak u^M,\allowbreak u^R)$ and $\Acanc(u^L,u^M, u^R)$) are well defined also in the scalar case (though this notation is not used in \cite{BaitMod1}). Moreover, $\Aquadr, \Acanc$ are mutually singular (i.e. if one amount is nonzero, the other one must be zero). In the two prototypical examples of \emph{same-sign} ($u^L \leqslant u^M \leqslant u^R$)  and \emph{opposite sign} ($u^L \leqslant u^R \leqslant u^M$) interactions, it holds 
$$
\begin{array}{|c|c|c|}
\hline
\multicolumn{3}{|c|}{\text{\it scalar case}}\\
\hline
    &\rule{0pt}{10pt} u^L \leqslant u^M \leqslant u^R & u^L \leqslant u^R \leqslant u^M \\ \hline 
   \rule{0pt}{10pt} \Aquadr(u^L, u^M, u^R)  & \Height(u^L, u^M, u^R) & 0 \\ \hline 
   \rule{0pt}{10pt} \Acanc(u^L, u^M, u^R) & 0 & |u^M - u^R| \\ \hline
\end{array}
$$

Though it is not done explicitly in \cite{BaitMod1}, also in the scalar case the proof of Theorem~\ref{thm:interaction-estimate} can be split into analogs of Prop.~\ref{prop:ie-intro-delta-sigma}--\ref{prop:ie-intro-delta-upsilon}. For $\Delta \sigma$, it holds 
(omitting, here and in this whole section, the arguments $(u^L, u^M, u^R)$ of the amounts of interaction and the change of speed):
\begin{align}
\label{eq:scalar-Dsigma}
    \Delta \sigma \leqslant 2 \Aquadr + \const \big( |u^R - u^M| + |u^M - u^L| \big) \Acanc,\quad n=1.
\end{align}
Indeed,~\eqref{eq:scalar-Dsigma} follows from \cite[(32)]{BaitMod1}  in the case of same-sign interactions,  while it follows from \cite[(34)]{BaitMod1} in the case of opposite-sign interactions. Notice that~\eqref{eq:scalar-Dsigma} is essentially Prop.~\ref{prop:ie-intro-delta-sigma} (without $\Atrans$ and $\sum_k \Acubic_k$, which we can treat as additional ``errors'' coming from system case).

Concerning Prop.~\ref{prop:ie-intro-delta-upsilon}, the scalar analog is
\begin{align}
    \label{eq:scalar-Dupsilon}
    \Aquadr + \left( |u^R - u^M| + |u^M - u^L| \right) \Acanc \leqslant \const \big[ \Upsilon(u) - \Upsilon(\hat u) \big], \quad n = 1.
\end{align}
Recall that, in the scalar case, $\Upsilon = C \TV^2 - \Qcost$, for some $C>0$ large enough, to be chosen.\footnote{In \cite{BaitMod1}, it is shown that it is enough to pick $C = \frac{1}{4} \norm{\derder{f}}_{L^\infty}$.}
Estimate~\eqref{eq:scalar-Dupsilon} is proven in \cite{BaitMod1} 
by analyzing separately $\TV^2(u) - \TV^2(\hat u)$ and $\Qcost(u) -\Qcost(\hat u)$. 
The two estimates are\footnote{These estimates are not written exactly in this form in \cite{BaitMod1} (in particular, the notation with amounts of interaction was not used), but we point out in the next paragraph the precise reference to estimates in \cite{BaitMod1}.}
\newcommand{\mclrbx}[2][0pt]{\setlength{\fboxsep}{#1}\colorbox{white}{$#2$}}
$$
\begin{array}[b]{p{2mm}|c|c|r|c}
\hhline{~|---|~}
&\multicolumn{3}{c|}{\text{\it scalar case}}&\\
\hhline{~|-|-|-|~}
&& \rule{0pt}{.7cm}\begin{minipage}[b]{1.8cm}
\centering
    same-sign interaction
\end{minipage}& \text{opposite-sign interaction\hspace{1cm}}
\\
\hhline{~|-|-|-|~}
    &\rule{0pt}{10pt}\Qcost(\hat u) - \Qcost(u) & \mathllap{\mclrbx{\geqslant}\hspace{3mm}}\hspace{1.8mm}\Aquadr\hspace{1.8mm} & \mathllap{\mclrbx{-\hspace{.4mm}}} \hspace{1mm}\const\left( |u^R - u^M| + |u^M - u^L| \right) \Acanc&{\refstepcounter{equation}(\theequation)\label{eq:scalar-Qcost}}\\
\hhline{~|-|-|-|~}
    &\rule{0pt}{10pt}\TV(u)^2 - \TV(\hat u)^2 & \mathllap{\mclrbx{\geqslant}\hspace{3mm}} 
    \hspace{5.5mm}0\hspace{5.5mm}
    & \mathllap{\mclrbx[2pt]{+}\hspace{6.35mm}}   2\left( |u^R - u^M| + |u^M - u^L| \right) \Acanc &{\refstepcounter{equation}(\theequation)\label{eq:scalar-TV2}}\\
\hhline{~|-|-|-|~}
\end{array}
$$

In the r.h.s. of  the above inequalities  there are always two terms: the first one is obtained from the analysis of \emph{same sign} interactions, the second one from \emph{opposite sign} interactions. 
For $\Delta \Qcost$ see \cite[(31)]{BaitMod1} and \cite[(35)]{BaitMod1} respectively. For $\Delta \TV^2$, the case of same sign-interactions is trivial, see \cite[(33)]{BaitMod1} for opposite-sign interactions.

These inequalities, together with the choice of a sufficiently large constant $C$ in the definition $\Upsilon$, are enough to conclude the proof of~\eqref{eq:scalar-Dupsilon} in the scalar case.  
Indeed, denote by $B>0$ the constants $\const$ appearing in~\eqref{eq:scalar-Qcost}; now $\eqref{eq:scalar-Qcost}+C\eqref{eq:scalar-TV2}$ reads 
\begin{align*}
    \Upsilon(u) - \Upsilon(\hat u) & \geqslant \Aquadr + (2C -  B) \left( |u^R - u^M| + |u^M - u^L| \right) \Acanc \\
    & \overbrace{\geqslant}^{(C \geqslant B)} \Aquadr + B \left( |u^R - u^M| + |u^M - u^L| \right) \Acanc 
\end{align*}
and this gives~\eqref{eq:scalar-Dupsilon}. In a nutshell, $-\Qcost$ can increase in case of opposite sign interactions, but its increase is controlled by the decrease of $\TV^2$. 

\subsection{Bound on \texorpdfstring{$\Delta \Upsilon$}{Delta Upsilon}}
\label{ss:delta-upsilon-bound}

We finally prove Prop.~\ref{prop:ie-intro-delta-upsilon}, adopting a similar strategy to the one described in the previous section for scalar case. We therefore analyze separately the variation of the various components  of $\Upsilon = C_1 V^2 + C_2\Qknown - \Qcost$. 
The constants $C_1, C_2$, with $C_1, C_2 - C_1 \gg 1$, will be fixed at the end of this section.
\begin{prop}[Variation of $V^2$]
\label{prop:var-v2}
    Let $u$, $\hat u$ be as in Theorem~\ref{thm:interaction-estimate}. It holds
    \begin{equation*}
        V^2(u) - V^2(\hat u) \geqslant 2 \sum_{k=1}^n (|s_k'| + |s_k''|) \sum_{h=1}^n \Acanc_h - \const \big[\Atrans + \sum_{k=1}^n \Acubic_k \big].
    \end{equation*} 
\end{prop}
\begin{proof}
We have
\begin{align*}
    V^2(u) -& V^2(\hat u)
     = (V(u) + V(\hat u)) (V(u) - V(\hat u)) \\
    & = (V(u) + V(\hat u)) \Big( \sum_h |s_h'| + |s_h''| - |s_h| \Big) \\
    & = (V(u) + V(\hat u)) \bigg[ \Big( \sum_h \underbrace{|s_h'| + |s_h''| - |s_h' + s_h''|}_{= 2\Acanc_h \geqslant 0} \Big) + \Big( \sum_h \underbrace{|s_h' + s_h''| - |s_h|}_{
    \geqslant -|s_h' + s_h'' - s_h|} \Big) \bigg] \\
    & 
    \geqslant 2 \sum_k \big( |s_k'| + |s_k''| \big) \sum_h \Acanc_h - \const \Big[\Atrans + \sum_k \Acubic_k \Big],
\end{align*}
where in the last inequality we used $\sum_k (|s_k'| + |s_k''| )\leqslant V(u) + V(\hat u) \leqslant 2 \beta \leqslant \const $, and Prop.~\ref{prop:sdef-quadr-bound}.
\end{proof}

We now recall the related Lyapunov property for $\Qknown$~\cite{Bianchini2003}.
\begin{prop}[{\cite[Prop.~4.1]{Bianchini2003}}]\label{prop:Qknown-Lyap2}
  Let $u$, $\hat u$, $u^L$, $u^M$, $u^R$ be as in Theorem~\ref{thm:interaction-estimate} and \indfixbeta{$\beta$ be small enough}. Then
    \begin{equation*}
        \Atrans(u^L, u^M , u^R)+\sum_{k=1}^n\Acubic_k(u^L, u^M , u^R)\leqslant\const\Big[\Qknown(u)-\Qknown(\hat u)\Big].
    \end{equation*}
\end{prop}

\begin{remark}
\label{rmk:qknown-lyap}
    Combining Prop.~\ref{prop:sdef-quadr-bound} (which is \cite[Theorem~3.7]{Bianchini2003}) with Prop.~\ref{prop:Qknown-Lyap2},  one can deduce (using the same notation as above) that
    \begin{equation*}
        (V + C\Qknown) (u) \geqslant (V+C\Qknown)(\hat u),
    \end{equation*}
    for $C \gg 1$.
\end{remark}

We now estimate the variation of $\Qcost$. This is the core proof of the whole paper. As in the scalar case, we prove that $-\Qcost$ has a decreasing part (corresponding to same-sign interactions) and an increasing part (corresponding to opposite-sign interactions, as well as to transversal ones and adjustments of the length of the curves). 
\begin{prop}[Variation of $\Qcost$]
\label{prop:var-cost}
     Let $u$, $\hat u$ be as in Theorem~\ref{thm:interaction-estimate}. If $\beta\ll 1$, then
    \begin{equation*}
        \Qcost(\hat u) - \Qcost(u) \geqslant \frac{1}{2} \sum_k \Aquadr_k - \const \sum_{k=1}^n\big(|s_k'|+|s_k''|\big) \Big( \Atrans + \sum_{h=1}^n \big( \Acubic_h  + \Acanc_h  \big) \Big).
    \end{equation*}
\end{prop}
\begin{proof}
First we observe that it is enough to show, for every $k$,
\begin{align}
    &\cost_k(u^L, u^R)-\cost_k(u^L, u^M) - \cost_k(u^M, u^R)\geqslant\Aquadr_k 
    \nonumber\\
    &- \const(\abs{s_k'}+\abs{s_k''}+\abs{s_k})\Big(\Atrans+\sum_{h=1}^n(\Acubic_h+\Aquadr_h+\Acanc_h)\Big).
    \label{eq:main-prop-k}
\end{align}
Indeed, assume that~\eqref{eq:main-prop-k} is proven. Adding for $k=1,\ldots, n$, the l.h.s. sums up to $\Qcost(\hat u) - \Qcost(u)$. Noticing that, due to Prop.~\ref{prop:sdef-quadr-bound}, 
\begin{align*}
    \sum_k \abs{s_k} \leqslant\sum_k\abs{s_k'+s_k''} +\const\Big[\sum_k|s_k'|\Big]\cdot{\overbrace{\Big[\sum_k|s_k''|\Big]}^{\leqslant\beta}}\leqslant\const(1+\beta)\sum_k(|s_k'|+|s_k''|)
\end{align*}
the sum of r.h.s. of~\eqref{eq:main-prop-k} for all $k$  is in turn bounded from below by
\begin{align*}
    \sum_k \Aquadr_k - \const&(1+\beta)\Big(\sum_k(|s_k'|+|s_k''|)\Big)\cdot \Big( \Atrans + \sum_h (\Acubic_h + \Aquadr_h + \Acanc_h )\Big) 
    \\
     =&\Big(1-\const(1 + \beta)\beta \Big)\sum_k \Aquadr_k \\
     &- \const(1+\beta)\Big(\sum_k(|s_k'|+|s_k''|)\Big)\cdot \Big( \Atrans + \sum_h (\Acubic_h  + \Acanc_h )\Big)
\end{align*}
which yields the conclusions, if \indfixbeta{$\beta$ is small enough}.

Let us thus prove~\eqref{eq:main-prop-k}. We consider separately the two cases $s_k', s_k''\geqslant 0$ (``same-sign interaction'') and $-s_k'\leqslant s_k''\leqslant 0$ (``opposite-sign interaction''). The remaining four cases are treated similarly.

\medskip

\textit{Same-sign interaction, $s_k', s_k''\geqslant 0$.} We adopt the notation $a = 0$, $b = s_k'$, $c = s_k'+s_k''$, $\bar c = s_k$. 
We have
\begin{align}
    &\cost_k(u^L, u^M) = \cost_{f_k'}(a,b),\quad\cost_k(u^M, u^R) = \cost_{f_k''}(b,c),\label{eq:costkuLuM}\\
    &\cost_k(u^L, u^R) = \cost_{f_k}(a,\bar c).\label{eq:costkuLuR}
\end{align}
{In order to apply the interaction-like estimate from Section~\ref{s:cost-func-and-prop-recall}, we will operate with the  ``fake flux'' $\tilde f_k'\cup\tilde f_k''$ defined on $[a,c]$, and then quantify the errors w.r.t.~\eqref{eq:costkuLuM} and~\eqref{eq:costkuLuR}}.\footnote{Recall that, as observed in Section~\ref{ss:system-two-rp}, $\tilde f_k' \cup \tilde f_k'' \in C^{1,1}(\I(s_k'+s_k''))$ and therefore we can apply all results from Section~\ref{s:cost-func-and-prop-recall}.}
Analogously to the scalar case, we exploit the super-additivity of $\cost$ (Prop.~\ref{prop:ILE-recall},~\eqref{prop:ILE-recall:superadd}) that gives
\begin{align}
    \cost_{ \hspace{.7mm}{\mathclap{\tilde f'_k}}\hspace{1mm}\cup\tilde f''_k}(a,c) -\cost_{ \hspace{.7mm}{\mathclap{\tilde f'_k}}\hspace{1mm}\cup\tilde f''_k}(a,b)- \cost_{ \hspace{.7mm}{\mathclap{\tilde f'_k}}\hspace{1mm}\cup\tilde f''_k}(b,c)\geqslant\Height_{ \hspace{.7mm}{\mathclap{\tilde f'_k}}\hspace{1mm}\cup\tilde f''_k}(a,b,c)=\Aquadr_k.
    \label{eq:SSI-cost-tilde-int-est}
\end{align}
Now we are to quantify the corrections due to change of fluxes and amounts of waves.
First, by Prop.~\ref{prop:stability-recall} with $b_1=b_2=b$,
\begin{align}
    &\abs{\cost_{\tilde f'_k\cup\tilde f''_k}(a,b) - \cost_k(u^L, u^M)} =  \abs{\cost_{\tilde f'_k}(a,b) - \cost_{f_k'}(a, b)}\nonumber \\
    &\leqslant\frac12(b-a)\norm{\derder{\tilde f'_k}  - \derder{f_k'}}_{L^1([a,b])}= \const\abs{s_k'}\Atrans,\label{eq:SSI-cost-uL-uM-error}
\end{align}
where the latter estimate holds by Lemma~\ref{lemma:310}. Similarly,
\begin{align}
    \abs{\cost_{\tilde f'_k\cup\tilde f''_k}(b,c) - \cost_k(u^M, u^R)} = \const\abs{s_k''}\Atrans.\label{eq:SSI-cost-uM-uR-error}
\end{align}
Secondly, again by Prop.~\ref{prop:stability-recall} with $b_1=c$, $b_2 = \bar c$, 
\begin{align}
     &\abs{\cost_{\tilde f'_k\cup\tilde f''_k}(a,c) - \cost_k(u^L, u^R)} = \abs{\cost_{\tilde f'_k\cup\tilde f''_k}(a,c) - \cost_{f_k}(a, \bar c)}\nonumber \\
     &\leqslant\frac12(s_k'+s_k'')\, \boxed{\norm{\derder{(\tilde f'_k\cup\tilde f''_k)} - \derder{f_k}}_{L^1(\I(s_k'+s_k'')\cap \I(s_k))}}\nonumber \\
     &+\frac12\max\left\{\norm{\derder{(\tilde f'_k\cup\tilde f''_k)}}_{L^\infty},\ \norm{\derder{f_k\vphantom{\tilde f}}}_{L^\infty}\right\}
     \cdot
     \frac{s'_k+s''_k+\abs{s_k}}2
     \cdot
     \underbrace{\abs{s'_k+s''_k-s_k}}_{\smash{\mathclap{\substack{\vspace{4mm}\\ \leqslant \const (\Atrans+\sum_{h=1}^n\Acubic_h)\\ \text{by Prop.~\ref{prop:sdef-quadr-bound}}}}}\hspace{6mm}}.\label{eq:SSI-cost-uL-uR-error-pre}
\end{align}

We estimate the framed $L^1$-norm in the r.h.s. as 
\begin{align*}
    &\overbrace{\norm{\derder{(\tilde f'_k\cup\tilde f''_k)}- \derder{\tilde f_k}}_{L^1(\I(s_k'+s_k''))}}^{{\substack{ \leqslant \const \sum_{h=1}^k (\Aquadr_h  + (|s_h'| + |s_h''|)\Acanc_h )     \\ \text{by Lemma~\ref{l:improved-311}}}}} + \overbrace{\norm{\derder{\tilde f_k}- \derder{f_k}}_{L^1(\I(s_k'+s_k'')\cap \I(s_k))}}^{\vphantom{\sum}\smash{\mathclap{\leqslant\const(\Atrans+
    \sum_{h=1}^n
    \Acubic_h)\text{ by~Lemma~\ref{lemma:313}}}}}
    \\
    &
    \leqslant \const\Big(\Atrans+\sum_{h=1}^n(\Acubic_h+\Aquadr_h + (|s_h'| + |s_h''|)\Acanc_h)\Big).
\end{align*}
Thus~\eqref{eq:SSI-cost-uL-uR-error-pre} gives, due to the maximum in the r.h.s. being uniformly bounded,
\begin{align}
    &\abs{\cost_{\tilde f'_k\cup\tilde f''_k}(a,c) - \cost_k(u^L, u^R)}\nonumber \\
    &\leqslant\rule{0pt}{.6cm}\const(s_k'+s_k''+\abs{s_k})\smash{\Big(\Atrans+\sum_{h=1}^n(\Acubic_h+\Aquadr_h + (|s_h'| + |s_h''|)\Acanc_h ) \Big).}\label{eq:SSI-cost-uL-uR-error}
\end{align}
By~\eqref{eq:SSI-cost-tilde-int-est}--\eqref{eq:SSI-cost-uM-uR-error} and~\eqref{eq:SSI-cost-uL-uR-error} we obtain
\begin{align*}
    &\cost_k(u^L, u^R)-\cost_k(u^L, u^M) - \cost_k(u^M, u^R)\geqslant\Aquadr_k 
    \nonumber\\
    &- \const (s_k'+s_k''+\abs{s_k})\rule{0pt}{.5cm}\smash{\Big(\Atrans+\sum_{h=1}^n(\Acubic_h+\Aquadr_h + (|s_h'| + |s_h''|)\Acanc_h  )\Big)}.
\end{align*}
This implies~\eqref{eq:main-prop-k}, in case $s_k', s_k'' \geqslant 0$.

\medskip
\textit{Opposite sign interaction, $-s_k'\leqslant s_k''\leqslant 0$}.
We adopt the notation $a = 0$, $b = s_k'+s_k''$, $c = s_k'$, $\bar b = s_k$. It holds
\begin{align}
    &\cost_k(u^L, u^M) = \cost_{f_k'}(a,c),\quad\cost_k(u^M, u^R) = \cost_{f_k''}(c,b),\label{eq:costkuLuM-os}\\
    &\cost_k(u^L, u^R) = \cost_{f_k}(a,\bar b).\label{eq:costkuLuR-os}
\end{align}

Using Prop.~\ref{prop:ILE-recall}~\eqref{prop:ILE-recall:subadd} and Prop.~\ref{prop:cost_Linfty_est-recall}, we obtain
\begin{align}
    &\cost_{\tilde f_k'}(a,b) - \cost_{\tilde f_k'}(a,c) \geqslant - \cost_{\tilde f_k'} (b,c)  -\frac12\norm{\derder{\tilde f_k'}}_{L^\infty} (c-b)(b-a) 
    \nonumber\\
    &\geqslant -\const \abs{s_k'}\cdot\abs{s_k''} = -\const\abs{s_k'}\Acanc_k.\label{eq:OSI-cost-tilde-int-est}
\end{align}
Analogously to~\eqref{eq:SSI-cost-uL-uM-error}, we have
\begin{align}
    \abs{\cost_{\tilde f_k'}(a,c) - \cost_{k}(u^L, u^M)} = \const\abs{s_k'}\Atrans.
    \label{eq:OSI-cost-uL-uM-error}
\end{align}
Also, by Prop.~\ref{prop:stability-recall} with $b_1= b$, $b_2= \bar b$
\begin{align}
    &\abs{\cost_{\tilde f_k'}(a,b) - \cost_k(u^L, u^R) } = \abs{\cost_{\tilde f_k'}(a,b) - \cost_{f_k}(a, \bar b) }
    \nonumber\\
    &\leqslant \frac12(s_k'+s_k'')\,\boxed{\norm{\derder{\tilde f_k'} - \derder{f_k}}_{L^1(\I(s_k'+s_k'')\cap \I(s_k))}}
    \nonumber\\
    &+\frac12\max\left\{\norm{\derder{\tilde f_k'}}_{L^\infty},\ \norm{\derder{\vphantom{\tilde f }f_k}}_{L^\infty}\right\}\cdot\frac{s_k'+s_k''+\abs{s_k}}2\cdot\underbrace{\abs{s_k'+s_k''-s_k}}_{\smash{\mathclap{\substack{\vspace{4mm}\\ \leqslant \const(\Atrans+\sum_{h=1}^n\Acubic_h)\\ \text{by Prop.~\ref{prop:sdef-quadr-bound}}}}}}.
    \label{eq:OSI-cost-uL-uR-error-pre}
\end{align}
We estimate the framed $L^1$-norm in the r.h.s. as
\begin{align*}
    \overbrace{\norm{\derder{\tilde f'_k}- \derder{\tilde f_k}}_{L^1(\I(s_k'+s_k''))}}^{{\substack{\vspace{4mm}\\ \leqslant \const\sum_{h=1}^k(\Aquadr_h+ (|s_h'|+|s_h''|)\Acanc_h)\\ \text{by Lemma~\ref{l:improved-311}}}}} \quad+\quad \overbrace{\norm{\derder{\tilde f_k}- \derder{f_k}}_{L^1(\I(s_k'+s_k'')\cap \I(s_k))}}^{\vphantom{\sum}\smash{\mathclap{\leqslant\const(\Atrans+
    \sum_{h=1}^n
    \Acubic_h)\text{ by Lemma~\ref{lemma:313}}}}}
    \hspace{2cm}
    \\
    \leqslant \const\Big(\Atrans+\sum_{h=1}^n(\Acubic_h+ \Aquadr_h + 
    (|s_h'|+|s_h''|)
    \Acanc_h)\Big).
\end{align*}
Thus~\eqref{eq:OSI-cost-uL-uR-error-pre} gives, due to the maximum in the r.h.s. being uniformly bounded,
\begin{align}
    &\abs{\cost_{\tilde f'_k}(a,b) - \cost_k(u^L, u^R)}\nonumber \\
    &\leqslant\const(s_k'+s_k''+\abs{s_k})\Big(\Atrans+\sum_{h=1}^n(\Acubic_h+ \Aquadr_h + (|s_h'|+|s_h''|)\Acanc_h)\Big).\label{eq:OSI-cost-uL-uR-error}
\end{align}
Finally, we roughly estimate $\cost_k(u^M, u^R)$ by quadratic bound (Prop.~\ref{prop:cost_Linfty_est-recall}),
\begin{align}
    \cost_k(u^M, u^R) \leqslant \const \abs{s_k''}^2\leqslant \const\abs{s_k'}\Acanc.
    \label{eq:OSI-cost-uM-uR}
\end{align}
By~\eqref{eq:OSI-cost-tilde-int-est}--\eqref{eq:OSI-cost-uL-uM-error} and \eqref{eq:OSI-cost-uL-uR-error}--\eqref{eq:OSI-cost-uM-uR} we obtain
\begin{align*}
    &\cost_k(u^L, u^R)-\cost_k(u^L, u^M) - \cost_k(u^M, u^R)\geqslant - \const\Big[\abs{s_k'}\Acanc_k
    \\
    &+(|s_k'|+|s_k''|+\abs{s_k})\Big(\Atrans+\sum_{h=1}^n(\Acubic_h+ \Aquadr_h + (|s_h'|+|s_h''|)\Acanc_h)\Big)\Big]
\end{align*}
This yields~\eqref{eq:main-prop-k}.
\end{proof}
\begin{remark}
    Observe that in~\eqref{eq:SSI-cost-uL-uR-error} and~\eqref{eq:OSI-cost-uL-uR-error} we in fact receive an additional prefactor $\abs{s_h'}+\abs{s_h''}$ before $\sum_h\Acanc_h$.
     This may be interpreted as follows: the contribution
     from $k$-th family to $\Delta\cost_k$ is $\Aquadr_k-\const(\abs{s_k'}+\abs{s_k''})\Acanc$ while the correction coming from $h$th family, $h\neq k$, has higher order prefactors before $\Aquadr_h$ and $\Acanc_h$. Accordingly, the r.h.s. of~\eqref{eq:main-prop-k} might be refined by adding this prefactors to all $\Acanc_h$, $h\neq k$. However, after adding~\eqref{eq:main-prop-k} for all $k$, the effect of refinement is not seen anymore, and so we keep~\eqref{eq:main-prop-k} with rougher r.h.s.
\end{remark}

\begin{proof}[Proof of Prop.~\ref{prop:ie-intro-delta-upsilon}]
    We now show that Prop.~\ref{prop:ie-intro-delta-upsilon} holds, provided we choose the constant $C_1, C_2-C_1$ in the definition of $\Upsilon$ large enough. Let us denote by $B>0$ (an upper bound for) the constants $\const$ appearing in Prop.~\ref{prop:var-v2}, \ref{prop:Qknown-Lyap2},~\ref{prop:var-cost} and let us  sum up the corresponding inequalities:
    \begin{align*}
        \Upsilon & (u)  - \Upsilon(\hat u) \\
        & = C_1 (V^2(u) - V^2(\hat u)) + C_2 (\Qknown(u) - \Qknown(\hat u)) - (\Qcost(u) - \Qcost(\hat u)) \\
        & \geqslant 2 C_1 \sum_{k} (|s_k'| + |s_k''|) \sum_h \Acanc_h - C_1 B \big[\Atrans + \sum_{k} \Acubic_k \big] \\
        & \qquad + C_2 B \big[ \Atrans + \sum_k \Acubic_k \big] \\
        & \qquad + \frac{1}{2} \sum_k \Aquadr_k - B \sum_k\big(|s_k'|+|s_k''|\big) \Big( \Atrans + \sum_{h} \big( \Acubic_h  + \Acanc_h  \big) \Big) \\
        & = \frac{1}{2} \sum_k \Aquadr_k + (2C_1 - B) \sum_k \big(|s_k'| + |s_k''| \big) \sum_h \Acanc_h \\ 
        & \qquad + (C_2 - C_1 - \beta) B \Big( \Atrans + \sum_k \Acubic_k \Big) \\
        & \geqslant \frac{1}{2} \Big( \Atrans + \sum_k  \Acubic_k + \Aquadr_k + \big( |s_k'|+ |s_k''| \big)\Acanc_k \Big),
    \end{align*}
   where the last inequality holds, provided that  \indfixC{$C_1$ and $C_2 - C_1$ are chosen large} \label{C1-C2-fix-2} \indfixC{enough}.
\end{proof}
 
\section{\texorpdfstring{$\Upsilon$}{Upsilon} is constant if no interaction takes place}
\label{s:persistency}

We  now state and prove the \emph{persistency} property, claimed in \ref{thm-main-a}, item~\eqref{thm-main:pers}.

\begin{theorem}[$\Upsilon$ is constant if no interactions]
\label{thm:const-if-no-interaction}
  Let $\beta \ll 1$. Let $u \in \mcD_\beta$ be piecewise constant with distance between adjacent jump points bounded from below by $L_0>0$. Then   
  for all $t \in [0,t_0)$, $S_t \bar u \in \mcD_\beta$ and $$\Upsilon(S_t u) = \Upsilon(u),$$
  with 
  \begin{equation*}
      t_0 := \frac{L_0}{\max
      \limits_{\xi, \eta\in\Beone}\abs{\lambda_n(\xi)-\lambda_1(\eta)}}
  \end{equation*}
\end{theorem}

\begin{remark}
\label{rmk:persistencty-bianchini}
    We will actually show (see the proof of Prop.~\ref{prop:syst-RP-consistency}) that each component of $\Upsilon$ (namely $V, \Qtrans, \Qcubic, \Qcost$) is constant up to time $t_0$. 
\end{remark}
We derive Theorem~\ref{thm:const-if-no-interaction} from the following proposition.
\begin{prop}\label{prop:syst-RP-consistency}
        Let $u(t,x)$ be the vanishing viscosity solution to  Riemann problem $(u^L,u^R)$. 
    Then for all $t>0$ we have $\Upsilon(u(t,\cdot)) = \Upsilon(\bar u)$.
\end{prop}

\begin{proof}[Proof of Prop.~\ref{prop:syst-RP-consistency}]
    For simplicity, consider the case $t=1$, the general one being completely analogous up to a rescaling of the spatial argument. Let $w = u(1,\cdot)$. We  prove separately that every component of $\Upsilon = C_1 V  + C_2(\Qtrans + \tfrac{1}{4} \Qcubic) - \Qcost$ takes the same values for $\bar u $ and $w$.

    According to the construction of the solution to the Riemann problems (Section~\ref{ss:sln-rp}), $w$ consists of $n$ parts patched together, $k$-th part corresponding to $k$-th family. Since the parts are spatially ordered by its number, $\Qtrans(w)=0$, and also, trivially, $\Qtrans(\bar u) = 0$. All the other components in $\Upsilon$, namely $V$, $\Qcubic$ and $\Qcost$, deal with all  families separately without measuring cross-family interaction. Having observed that, we will now assume, for the simplicity of notation, that only one family is present. The general case will follow by analyzing all families in the manner reproduced below and adding their contributions.

    So, we assume that $u^L$ and $u^R$ are joined by an exact curve $\gamma_k(u^L, s) \allowbreak = (u_k, v_k, \sigma_k)\colon\allowbreak \I(s)\to\R^n\times \R\times \R$ for some fixed $k$ and $s$. W.l.o.g., we also assume $s>0$ (and thus $\I(s) = [0,s]$).
    Denote by $f_k$ the reduced $k$-flux. Denote also by $\mu_k$ the $k$-th wave measure of $w$.

    Let now $\mcA$ be the (at most countable) collection of maximal intervals $A\subset \I(s)$ of positive measure on which $\sigma_k$ is constant. We prove first of all that
    \begin{align}
    \label{eq:w-pf-jump-cont}
        w_x^\atom & ={\sigma_k}_\sharp\Big((\tilde r_k\circ\gamma_k)\cdot\mcL\restr{\cup\mcA}\Big), & 
        w_x^\cont & ={\sigma_k}_\sharp\Big((\tilde r_k\circ\gamma_k)\cdot\mcL\restr{\I(s)\setminus (\cup\mcA)}\Big).
    \end{align}
    By the construction of the solution to the Riemann problems,
    \begin{align}
        w_x \stackrel{\eqref{eq:rp-sln-k-pf}}={\sigma_k}_\sharp\Big((\tilde r_k\circ\gamma_k)\cdot\mcL\restr{\I(s)}\Big).\label{eq:w-pf}
    \end{align}
    In particular, we have $w_x(\{x\})\neq 0\Longrightarrow\abs{\sigma_k\inv(\{x\})}>0\Longrightarrow\sigma_k\inv(\{x\})\in \mcA$. Therefore $ w_x^\atom = \sum_{A \in \mcA} w_x \big( \sigma_k (A) \big) \delta_{\sigma_k(A)}$.
    We deduce, for $\varphi \in C_c(\R)$,
    \begin{align*}
        \int \varphi(x) w_x^\atom = \sum_{A \in \mcA} \varphi(\sigma_k(A)) w_x(\sigma_k(A)) \overset{\eqref{eq:w-pf}}{=} \sum_{A \in \mcA} \int_A \varphi(\sigma_k(\varsigma)) \tilde r_k(\gamma_k(\varsigma)) d\varsigma,
    \end{align*}
    from which the  equality for $w_x^\atom$ in~\eqref{eq:w-pf-jump-cont} follows. The equality for $w_x^\cont$ follows now from the one for $w_x^\atom$ and~\eqref{eq:w-pf}. 

    Observe now that, by~\eqref{eq:w-pf-jump-cont}, for $A = [a_-, a_+] \in \mcA$ and $\sigma_k(A) = \{x\}$, the map $\gamma_k\restr{A}$ is an exact curve joining $w_-(x) = u_k(a_-)$ and $w_+(x) = u_k(a_+)$, with the corresponding flux $f_k\restr{A}$ (up to an appropriate shift of the domain). We infer
    \begin{align}\label{eq:muh-atom-and-costh}
        \mu_h(\{x\}) = \abs{A}\delta_{kh}, \qquad
        \cost_h(w_-(x), w_+(x)) = \cost_{f_k}(a_-,a_+)\delta_{kh}.
    \end{align}
    
    We now show that 
    \begin{align}
    \label{eq:mu-pf}
        \mu_h^\atom   = \delta_{hk} {\sigma_k}_\sharp (\mcL\big|_{\cup \mcA}) ,\qquad \mu_h^\cont   = \delta_{hk} {\sigma_k}_\sharp (\mcL\big|_{\I(s) \setminus \cup \mcA}).
    \end{align}
    The first equality follows immediately from \eqref{eq:w-pf-jump-cont},~\eqref{eq:muh-atom-and-costh}.
    Next, let $s\in\I(s)\setminus (\cup \mcA)$. We have $v_k(s) = 0$ and so $(\tilde r_k\circ\gamma_k)(s) = (r_k\circ u_k)(s) $; also $(w\circ\sigma_k)(s) \stackrel{\eqref{eq:rp-sln-k}}= u_k(s) $. Thus, given $\varphi\in C_c(\R)$ we obtain, by definition of $\mu_h$,
    \begin{align*}
        &\int\varphi\mu_h^\cont= \int\varphi(x)\scal{\ell_h(w(x))}{w_x^\cont(dx)} \\
        &\stackrel{\mathclap{\eqref{eq:w-pf-jump-cont}}}{=}
\int\limits_{\mathrlap{\I(s)\setminus (\cup\mcA)}}\hspace{3mm}(\varphi\circ\sigma_k)(s)\scal{\smash{\underbrace{(\ell_h\circ w\circ\sigma_k)(s)}_{=\ell_h(u_k(s))}}}{\smash{\underbrace{(\tilde r_k\circ\gamma_k)(s)}_{=r_k(u_k(s))}}}\, ds
\\
&= \int\limits_{\mathrlap{\I(s)\setminus{( \cup\mcA)}}}\hspace{3mm}(\varphi\circ\sigma_k)\delta_{kh}\, ds,
\end{align*}
that implies the second equality in~\eqref{eq:mu-pf}.

From~\eqref{eq:mu-pf}, we  
obtain $\mu_h = \delta_{hk} {\sigma_k}_\sharp(\mcL\restr{\I(s)})$, which implies, by monotonicity of $\sigma_k$, $\mu_k = \abs{\mu_{k}}$ and thus $V(w) = V(\bar u) = |s|$.

Let us now prove that $\Qcubic(w) = \Qcubic(\bar u)$. Recall that, as before in this proof, $\sigma_k$ is the third component of the exact curve solving the Riemann problem $(u^L, u^R)$. Clearly $\Qcubic(\bar u)$ can be computed in~\eqref{eq:qcubic} using $\sigma_k$. We now show that also $\Qcubic(w)$ can be computed  in~\eqref{eq:qcubic} using the same map $\sigma_k$. This amounts to showing that $\sigma_k$ satisfies~\eqref{eq:abs-sigma-map} for $w$. 

Fix thus $\tau\in(0,\abs{\mu_k}(\R)) = (0,s)$. Then the maximal number $x$ satisfying $\abs{\mu_k((-\infty,x))}\leqslant\tau$ is equal to $x=\sigma_k(\tau)$. Now
\begin{itemize}
    \item if $\mu_k(\{x\}) = 0$, then $\tau\in\I(s)\setminus(\cup\mcA)$ is a singleton $\sigma_k$-preimage of $x$, that means $w_-(x) = w(x) = u_k(\tau) $ and $\sigma_k(\tau) = \tilde\lambda_k(\gamma(\tau)) = \lambda_k(u_k(\tau))$;
    \item if $\mu_k(\{x\})>0$, then $A:=\sigma_k\inv(\{x\})$ belongs to $\mcA$, and
\begin{align*}
    \sigma_k(\tau) = \frac{d}{d\varsigma}(\conv_{\I(s)} f_k)(\tau) = \frac{d}{d\varsigma}(\conv_{A} f_k)(\tau),
\end{align*}
\end{itemize}
so that, overall, $\sigma_k$ satisfies~\eqref{eq:abs-sigma-map} for $w$.
This proves $\Qcubic(w) = \Qcubic(\bar u)$.

Finally, having atoms of $w_x$ identified above as  $\sigma_k$-images of elements of $\mcA$, one writes
\begin{align*}
    &\Qcost(w) = \sum_{x}\cost_k(w_-(x),w_+(x))\stackrel{\eqref{eq:muh-atom-and-costh}}{=}\sum_{[a_-,a_+]\in\mcA}\cost_{f_k}(a_-,a_+)\\
&\xlongequal{\text{Prop.~\ref{prop:persistence}}}\sum_{[a_-,a_+]\in\mcA}\sum_{\substack{(p,q)\in\\
\bridges_{f_k}(a_-,a_+)}}\cost_{f_k}(p,q) = \sum_{\substack{(p,q)\in\\
\bridges_{f_k}(0,s)}}\cost_{f_k}(p,q)=\Qcost(\bar u). 
\end{align*}
In the second last equality, we used the identity 
\begin{equation*}
    \bigsqcup_{[a_-,a_+]\in \mcA}\bridges_{f_k}(a_-,a_+) = \bridges_{f_k}(0,s).\qedhere
\end{equation*}
\end{proof}

\begin{proof}[Proof of Theorem~\ref{thm:const-if-no-interaction}]
    Consider two adjacent jumps of $u$. Due to Lax admissibility condition, any speed originating from the first jump differs from any speed originating from the second by at most
    \begin{equation*}
    \max
      \limits_{\xi, \eta\in\Beone}\abs{\lambda_n(\xi)-\lambda_1(\eta)}
    \end{equation*} 
    Also, the jumps are separated spatially by at least $L_0$.
    Thus, before time $t_0$, wavefronts starting (at time zero) from different points can not interact.
    Therefore,
    for all $t<t_0$,  by Prop.~\ref{prop:syst-RP-consistency}, $\Upsilon(S_t \bar u) = \Upsilon(\bar u)$.
\end{proof}

\section{The approximation property w.r.t. cropping}
\label{s:crop}

In this and the next section we state the two theorems (Theorems~\ref{prop:sup-prop-crop}--\ref{thm:sup-prop-cut}), corresponding to \ref{thm-main-b}, item~\eqref{thm-main:appr}. In this section we discuss approximations of a given map $u$ by means of a \emph{cropping} procedure, consisting in removing variations of $u$ outside compact sets.

\begin{definition}
    \label{def:crop}
Let $I \subseteq \R$ be a nonempty interval. 
For $u \in \BV(\R;\R^n)$, we define the \emph{cropping} $\Cr{I}(u)$ as the map $v \in \BV(\R;\R^n)$ such that
\begin{align}\label{eq:crop-def}
    v_x = v_x\restr{I} =  u_x\restr{I} \quad\text{ and }v =  u\text{ a.e. on } I.
\end{align}
\end{definition}
Note that
\begin{align}\label{eq:absorption-crop}
   \Cr{I}\circ\Cr{J} = \Cr {I\cap J} \text{ whenever } I \cap J\neq \varnothing.
\end{align}
\begin{theorem}[Approximation w.r.t. $\Cr{}{}$]
\label{prop:sup-prop-crop}
    Let $\Theta \subseteq \R$ be such that $\inf \Theta = - \infty$, $\sup \Theta = + \infty$. For every $u \in \BV(\R;\Beone)$,
    \begin{align*}
        \Upsilon(u) = \sup_{\substack{I \subseteq \R \\ \inf I,  \sup I \in \Theta}} \Upsilon(\Cr{I} (u)).
    \end{align*}
\end{theorem}

To prove Theorem~\ref{prop:sup-prop-crop} we first discuss in Section~\ref{ss:perturb-ins-interval} how $\Upsilon(u)$ (precisely: each component of $\Upsilon$) changes if we modify the values of $u$ on an interval $I \subseteq \R$. Then we show that the operation $\Cr{}$ reduces the value of $\Upsilon$ (Section~\ref{ss:crop-Upsilon-reduce}), and finally we prove that the decrease vanishes as the limits of $\Cr{}$ approach $\pm\infty$  (Section~\ref{ss:crop-Upsilon-reduce-e}).

\subsection{Perturbation inside an interval}
\label{ss:perturb-ins-interval}
\begin{prop}
\label{prop:u-v-I}
Assume that
\begin{align}\label{eq:u-v-I}
    u,\tilde u \in \BV(\R;\Beone),\quad u_x = \tilde u_x\text{ on }\R \setminus A,\quad\text{ and }u = \tilde u\text{ a.e. on }\R \setminus A,
\end{align}
where $A\subseteq\R$ is an interval. Then
\begin{align*}
    \abs{\Upsilon(\tilde u) - \Upsilon(u)}\leqslant\const\Big[&(\TV(u) + \TV(\tilde u))\sum_k\abs{V_k(\tilde  u)-V_k(u)} \\
    &+\abs{u_x}(A)^2+\abs{\tilde u_x}(A)^2\Big](1+\TV(u)+\TV(\tilde u)).
\end{align*}
\end{prop}
The proof is based on the following three lemmas assuming~\eqref{eq:u-v-I}. For the proofs, recall the notation $A_\pm$ (see Section~\ref{ss:notation}). We denote~$\mu_k$ (resp.~$\tilde\mu_k$) the $k$-th wave measure of $u$ (resp.~$\tilde u$).

\begin{lemma}\label{l:alt-appr-Qtrans}
 Let~\eqref{eq:u-v-I} holds. Then
\begin{align*}
\abs{\Qtrans(\tilde u) - \Qtrans(u)}\leqslant\const\Big[&
\abs{u_x}(\R\setminus A){\sum_{k=1}^n\abs{V_k(\tilde u)-V_k(u)}} \\
&+ \abs{u_x}(A)^2+\abs{\tilde u_x}(A)^2\Big] .  
\end{align*}
\end{lemma}
\begin{proof}
    We have
    \begin{align*}
       \Qtrans(u) = \sum_{k'<k''}\big(\abs{\mu_{k'}\times\mu_{k''}} \big)(\udhp),\quad \Qtrans(\tilde u) = \sum_{k'<k''}\big(\abs{\tilde \mu_{k'}\times\tilde\mu_{k''}} \big)(\udhp).
    \end{align*}
    and, due to~\eqref{eq:u-v-I}, 
    \begin{align*}
        \Qtrans(\tilde u) - \Qtrans(u) =
        \sum_{k'<k''}\smash{\big(\abs{\tilde\mu_{k'}\times\tilde\mu_{k''}}-\abs{\mu_{k'}\times\mu_{k''}}\big)(\overbrace{\udhp\cap[A\times \R\cup \R\times A]}^{\mathclap{(A\times A_-)\ \sqcup\ (A_+\times A)\ \sqcup\  (\smalludhp\cap A^2)}}).}
    \end{align*}
    Then,
    \begin{align*}
        &\big(\abs{\tilde\mu_{k'}\times\tilde\mu_{k''}}-\abs{\mu_{k'}\times\mu_{k''}}\big)
        \smash
        {
            ({A\times A_-}) = \overbrace{(\abs{\tilde\mu_{k'}} - \abs{\mu_{k'}})(A)}^{\smash{V_{k'}(\tilde u) - V_{k'}( u)}}\overbrace{\abs{\mu_{k''}}(A_-)}^{\smash{\const\abs{u_x}(\R\setminus A)}}
        },\\
            &\big(\abs{\tilde\mu_{k'}\times\tilde\mu_{k''}}-\abs{\mu_{k'}\times\mu_{k''}}\big)({A_+\times A}) =
            \smash{\underbrace{\abs{\mu_{k'}}(A_+)}_{\smash{\const\abs{u_x}(\R\setminus A)}}\underbrace{(\abs{\tilde\mu_{k''}} - \abs{\mu_{k''}})(A)}_{\smash{V_{k''}(\tilde u) - V_{k''}( u)}}.}
    \end{align*}
    and 
    \begin{align*}
    &\Big|\big(\abs{\tilde\mu _{k'}\times\tilde\mu _{k''}}-\abs{\mu_{k'}\times\mu_{k''}}\big)(\udhp\cap A^2)\Big|\leqslant\big(\abs{\tilde\mu _{k'}\times\tilde\mu _{k''}}+\abs{\mu_{k'}\times\mu_{k''}}\big)(\udhp\cap A^2)\\
    &\leqslant\const(\abs{u_x}(A)^2+\abs{\tilde u_x}(A)^2)
    \end{align*}
    Gathering all estimates together, we conclude.
\end{proof}
\begin{lemma}\label{l:alt-appr-Qcubic}
Let~\eqref{eq:u-v-I} holds, then
    \begin{align*}
        &\abs{\Qcubic_k(\tilde u)-\Qcubic_k(u)}\\
        &\leqslant\const\Big[\abs{u_x}(\R\setminus A) \cdot \abs{V_k(\tilde  u)-V_k(u)}+\abs{u_x}(A)^2+\abs{\tilde u_x}(A)^2\Big](\TV(u)+\TV(\tilde u)).
    \end{align*}
\end{lemma}
\begin{proof}
    Let $\sigma_k,\tilde \sigma_k$ denote the speed maps associated with $u$ and $\tilde  u$ respectively as in~\eqref{eq:abs-sigma-map}.
    Put
    \begin{align*}
        &S_k = (0,\abs{\mu_k}(\R)),\qquad Z_k = (\abs{\mu_k}(A_-), \abs{\mu_k}(A_-\cup A))\\
        &\tilde  S_k = (0,\abs{\tilde \mu_k}(\R)),\qquad \tilde  Z_k = (\abs{\tilde \mu_k}(A_-), \abs{\tilde \mu_k}(A_-\cup A)).
    \end{align*}

    Note that $\tilde  \sigma_k$ can be obtained from $\sigma_k$ by insertion the subinterval $\tilde  Z_k$ in place of $ Z_k$ and changing values there with possible shifting the values on the right of it, see Fig.~\ref{fig:hat_Z}.
    \begin{figure}[h]
        \centering
        \includegraphics{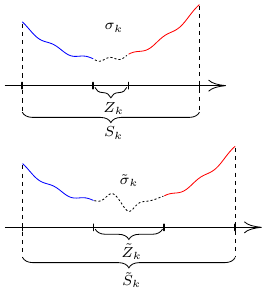}
        \caption{To the proof of Lemma~\ref{l:alt-appr-Qcubic}: the structure of $\sigma_k$ and $\tilde \sigma_k$.}
        \label{fig:hat_Z}
    \end{figure}
        Assume that $\varnothing\subsetneq A\subsetneq\R$, the opposite case is trivial. W.l.o.g., we may assume $\inf A\in\R$. 
        Put $p = u(-\infty)+u_x(A_-)$.
        For every $x\in A$, one has $\abs{u(x)-p}\leqslant\abs{u_x}(A)$, $\abs{\tilde u(x)-p}\leqslant\abs{\tilde u_x}(A)$. Hence, by the construction of $\sigma_k$, one has
    \begin{align}
        \abs{\sigma_k(s) - \lambda_k(p)} \leqslant \const \abs{u_x}(A)\quad (s\in Z_k),\label{eq:dif-sigma-k}\\
        \abs{\tilde\sigma_k(s) - \lambda_k(p)} \leqslant \const \abs{\tilde u_x}(A)\quad (s\in \tilde Z_k).\label{eq:dif-tsigma-k}
    \end{align}
    
    We have 
    \begin{align*}
        \Qcubic_k(u) = \iint\limits_{\smash{S_k\times S_k}}\abs{\sigma_k(s') - \sigma_k(s'')}\, ds'ds'', \hspace{1mm}
    \Qcubic_k(\tilde  u) = \iint\limits_{\smash{\tilde  S_k\times\tilde  S_k}}\abs{\tilde  \sigma_k(s') - \tilde  \sigma_k(s'')}\, ds'ds''.
    \end{align*}
    By the structure of $\sigma_k$, 
    \begin{align*}
        \Qcubic_k(\tilde u)-\Qcubic_k(u) = 
        \iint\limits_{\mathclap{(\tilde  Z_k\times\tilde   S_k)\cup (\tilde  S_k\times \tilde  Z_k)}}\abs{\tilde \sigma_k(s') - \tilde \sigma_k(s'')}\, ds'ds'' -\iint\limits_{\mathclap{(Z_k\times S_k)\cup (S_k\times Z_k)}}\abs{\sigma_k(s') - \sigma_k(s'')}\, ds'ds''.
    \end{align*}
    Let us now show that, in the latter difference, the subtrahend (resp. minuend) is close to the quantity $2|Z_k|\cdot\text{int}$ (resp.~$2|\tilde Z_k|\cdot\text{int}$), where
    \begin{align*}
        \text{int} :=\int\limits_{ S_k\setminus Z_k} \abs{\lambda_k(p)  - \sigma_k(s)} \, ds = \int\limits_{ \tilde S_k\setminus \tilde Z_k} \abs{\lambda_k(p)  - \tilde\sigma_k(s)} \, ds.
    \end{align*}
    Indeed, using $\abs{Z_k}\leqslant \const \abs{u_x}(A)$, we obtain
    \begin{align*}
        &\Big|\text{subtrahend}-2\abs{Z_k}\cdot \text{int}\Big| \leqslant \Bigg|\ \iint\limits_{Z_k\times Z_k}\abs{\sigma_k(s') - \sigma_k(s'')}\, ds'ds''\Bigg| \\
        &+2\iint\limits_{\mathclap{Z_k\times (S_k\setminus Z_k)}}\hspace{5mm} \Bigg|\abs{\sigma_k(s') - \sigma_k(s'')} -\abs{\lambda_k(p)-\sigma_k(s'')}\Bigg|\, ds'ds''\leqslant\const\abs{u_x}(A)^3\\
        &+2\iint\limits_{\mathclap{Z_k\times (S_k\setminus Z_k)}} \abs{\sigma_k(s')- \lambda_k(p)}\, ds'ds''\stackrel{\eqref{eq:dif-sigma-k}}\leqslant \const \abs{u_x}(A)^2\TV(u).
    \end{align*}
    Analogously, using~\eqref{eq:dif-tsigma-k}, we get $
        \Big|\text{minuend}-2|\tilde Z_k|\cdot \text{int}\Big| \leqslant \const \abs{\tilde u_x}(A)^2\TV(\tilde u)$.
    Then
    \begin{align*}
        &\abs{\Qcubic_k(\tilde  u)-\Qcubic_k(u)} \\ &\leqslant2\Abs{\underbrace{|\tilde  Z_k|-|Z_k|}_{V_k(\tilde  u) - V_k(u)}}\cdot \text{int}+\const(\abs{u_x}(A)^2+\abs{\tilde u_x}(A)^2)(\TV(u)+\TV(\tilde u)).
    \end{align*}
    Finally, due to $\text{int}\leqslant\const\abs{u_x}(\R\setminus A)\TV(u)$, we conclude.
\end{proof}
\begin{lemma}\label{l:alt-appr-Qcost}
Let~\eqref{eq:u-v-I} holds, then $$
        \abs{\Qcost(\tilde u)-\Qcost(u)}\leqslant\const (\abs{u_x}(A)^2+\abs{\tilde u_x}(A)^2).$$
\end{lemma}
\begin{proof}
    We have
    \begin{align*}
        \Qcost_k(u) = \sum_{\mathclap{\abs{u_x}(\{x\})>0}}\cost_k(u_-(x),u_+(x)), \quad \Qcost_k(\tilde u) = \sum_{\mathclap{\abs{\tilde u_x}(\{x\})>0}}\cost_k(\tilde u_-(x),\tilde u_+(x))
    \end{align*}
    and, by~\eqref{eq:u-v-I},
    \begin{align*}
        \Qcost_k(\tilde  u)-\Qcost_k(u) =  \sum_{\substack{x\in A\colon\\ \mathclap{\abs{\tilde u_x}(\{x\})>0}}}\cost_k(\tilde u_-(x),\tilde u_+(x)) -  \sum_{\substack{x\in A\colon\\ \mathclap{\abs{u_x}(\{x\})>0}}}\cost_k(u_-(x),u_+(x)).
    \end{align*}
    Using the quadratic bound on $\cost$ (Prop.~\ref{prop:cost_Linfty_est-recall}), we bound the terms in the r.h.s. by $\const\abs{u_x}(A)^2$ and $\const\abs{\tilde u_x}(A)^2$, respectively.
\end{proof}
\begin{proof}[Proof of Prop.~\ref{prop:u-v-I}]
Combine Lemmas~\ref{l:alt-appr-Qtrans}--\ref{l:alt-appr-Qcost} with the inequality 
\begin{align*}
\abs{V^2(\tilde u)-V^2(u)}&=\abs{V(\tilde u)+V(u)}\cdot \abs{V(\tilde u)-V(u)}\\
&\leqslant\const(\TV(u)+\TV(\tilde u)\sum_k\abs{V_k(\tilde u)-V_k(u)}    \qedhere
\end{align*}
\end{proof}

\subsection{Proof of ``\texorpdfstring{$\geqslant$}{>=}'' in Theorem~\ref{prop:sup-prop-crop}}
\label{ss:crop-Upsilon-reduce}

\begin{lemma}\label{lemma:crop-Upsilon-reduce}
    Let $u\in\BV(\R;\Beone)$ and $I\subseteq\R$ be a nonempty interval, then $\Upsilon(\Cr I (u))\leqslant\Upsilon(u)$.
\end{lemma}
\begin{proof}
    There are three possible nontrivial cases: either $\inf I\in\R$, $\sup I=+\infty$, or $\sup I\in\R$, $\inf I = -\infty$, or $I$ is bounded. The first two cases are analogous, while the third one reduces to them through using~\eqref{eq:absorption-crop}. Thus, let us w.l.o.g. assume the first case. Notice that $\R\setminus I$ is now ensured to be an interval.

    Put $v = \Cr {\R\setminus I} (u)$, $w = \Cr I (u)$.
    We claim that 
    \begin{align}\label{eq:superadd-crop}
        \Upsilon(u)\geqslant\Upsilon(v)+\Upsilon(w).
    \end{align}
    Indeed, the $k$-th wave measure of $w$ (resp. $v$) is the restriction of the $k$-th wave measure of $u$ onto $I$ (resp.~$\R\setminus I$). Thus, using $u_x\restr{I} = w_x\restr I$, $u(+\infty) = w(+\infty)$, and $u\restr{\R \setminus I} = v \restr{\R \setminus I}$, $u(-\infty) = v(-\infty)$, we get $V(u) = V(v)+V(w)$ (and hence $V(u)^2 \geqslant V(v)^2+V(w)^2$) and $\Qcost(u) = \Qcost(v)+\Qcost(w)$. Moreover, the inequality~\eqref{eq:superadd-crop} holds true for $\Qtrans$, $\Qcubic$ instead of $\Upsilon$ due to the  cross-terms. This justifies~\eqref{eq:superadd-crop} and concludes the proof.
\end{proof}

\subsection{Proof of ``\texorpdfstring{$\leqslant$}{<=}'' in Theorem~\ref{prop:sup-prop-crop}}
\label{ss:crop-Upsilon-reduce-e}
The proof of ``$\leqslant$'' in Theorem~\ref{prop:sup-prop-crop} follows immediately from the following lemma.

\begin{lemma}\label{lemma:crop-Upsilon-reduce-e}
    Let $u\in\BV(\R;\Beone)$ and $I\subseteq\R$ be a nonempty interval, then 
    \begin{align*}
        \abs{\Upsilon(\Cr I (u))-\Upsilon(u)}\leqslant\const \abs{u_x}(\R\setminus I) \cdot\TV(u)\cdot (1+\TV(u)).
    \end{align*}
\end{lemma}
\begin{proof}
    Again, by  using~\eqref{eq:absorption-crop}, the general case reduces to the case of half-unbounded~$I$. Then $\R\setminus I$ is an interval, and applying Prop.~\ref{prop:u-v-I} to $\tilde u =\Cr I (u)$,  $A = \R\setminus I$, we get
    \begin{align*}
         \abs{\Upsilon(\Cr I (u))-\Upsilon(u)}&\leqslant\const\Big(\TV(u)\cdot \abs{u_x}(A)+\abs{u_x}(A)^2\Big)(1+\TV(u))\\
         &\leqslant\const\abs{u_x}(A)\cdot \TV(u)\cdot(1+\TV(u)).\qedhere
    \end{align*}
\end{proof}

\section{The approximation property w.r.t. Riemann cuts}
\label{s:cut}
We now discuss approximations of a given map $u$ by means of a (finite number of)  \emph{Riemann cuts}. This is the second statement (together with Theorem~\ref{prop:sup-prop-crop}) corresponding to \ref{thm-main-b}, item~\eqref{thm-main:appr}.

We have already informally introduced Riemann cut in Section~\ref{s:on-add-prop}, see Fig.~\ref{fig:riemann_cut}. 
It is nevertheless convenient to adopt a more general definition of Riemann cut.

Let $I \subseteq \R$ be an interval. We denote by $\mcP(I)$ the set of all finite partitions of $I$ into subintervals. On $\mcP(I)$ we define the partial ordering 
$\mcF' \preccurlyeq \mcF$ if for every $I' \in \mcF'$ there is $I$ in $\mcF$ such that $I' \subseteq I$. 

\begin{figure}[h]
    \centering
    \includegraphics{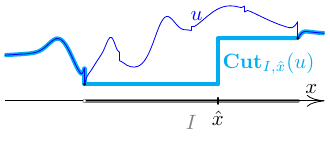}
    \caption{Edge jumps in a simple Riemann cut in the case $\inf I\notin I$, $\sup I\in I$. }
    \label{fig:riemann_cut_jumps}
\end{figure}

\begin{definition}
    \label{def:riemann-cut}
Let  $\mcF= \{ I_i\}_{i=1}^m \in \mcP(I)$. We say that the set $X = \{x_i\}_{i=1}^m$ is \emph{compatible} with $\mcF$  if $x_i \in I_i$. 
In that case, for $u \in \BV(\R; \R^n)$, we define the \emph{Riemann cut} $\Rc\mcF X(u)$ as the map $\hat u$ such that
\begin{equation*}
    \hat u(- \infty) = u(-\infty) \text{ and } \hat u_x = u_x - u_x\restr {I} + \sum_{i=1}^m u_x(I_i) \delta_{x_i}.
\end{equation*}
If  $\mcF = \{I\}$ and $X = \{\hat x\}$, we say that the Riemann cut is \emph{simple} and write $\Rc I{\hat x}(u)$.
\end{definition}
In particular, closed (resp. open) end of the interval in the simple Riemann cut corresponds to the inclusion (resp. not inclusion) of the edge jumps, if there are any, to the resulting $\BV$ function, see Fig.~\ref{fig:riemann_cut_jumps}.

Notice that, for $u\in\BV(\R;\Beone)$, the value of $\Upsilon(\Rc \mcF X(u))$ does not depend  on the particular choice of $X$. In view of that, we will write shortly $\Upsilon(\Rc \mcF \bullet(u))$. By Theorem~\ref{thm:interaction-estimate}, \indfixbeta{if $\beta \ll 1$,}
\begin{align}
    \label{eq:lyap-prop}
    \forall u \in \mcD_\beta\colon\quad \mcF' \preccurlyeq \mcF \ \Longrightarrow \ \Upsilon(\Rc {\mcF'}{\bullet}(u)) \geqslant \Upsilon(\Rc \mcF \bullet (u)).
\end{align}

\bigskip

We can now state the following approximation result.
\begin{theorem}[Approximation w.r.t. $\RcOp$]
\label{thm:sup-prop-cut}
  Let $\beta \ll 1$. Let  $\Theta\subseteq I$ be dense in $I$. Then, for every $u\in\mcD_\beta$,
\begin{align*}
    \Upsilon(u) = \sup \Big\{ \Upsilon \big(\Rc \mcF \bullet(u) \big) \colon \mcF \in \mcP_\Theta(I)
    \Big\},
\end{align*}
where
\begin{align*}
    \mcP_{\Theta}(I) = \big\{\mcF \in \mcP(I) \text{ s.t. } \forall J \in \mcF, \inf J, \sup J \in \Theta \cup \{ \inf I, \sup I\}\big\}.
\end{align*}
\end{theorem}
As an immediate corollary, we get \ref{thm-main-b}, item~\eqref{thm-main:rcut}.

\begin{corollary}
\label{cor:riemann-cut}
Let $\beta\ll 1$ and $u\in \mcD_\beta$. Let $I\subseteq\R$ be an interval and $\hat x\in I$. Then $\Upsilon(\Rc I{\hat x} (u))\leqslant\Upsilon(u)$.

\end{corollary}

\begin{proof}
Apply Theorem~\ref{thm:sup-prop-cut} to the case $\mcF=\{I\}$.
\end{proof}

Note that, in view of~\eqref{eq:lyap-prop}, the finer $\mcF$ is, the greater value $\Upsilon(\Rc\mcF \bullet(u))$ takes.
To prove Theorem~\ref{thm:sup-prop-cut}, we first develop  estimates in spirit of Section~\ref{ss:perturb-ins-interval} by analyzing separately intervals where the total variation is \emph{small} (Section~\ref{ss:direct-est})  and those where \emph{large} jumps are present (Section~\ref{ss:punct-I}). 
We then introduce a class of \emph{$\e$-accurate partitions} (Section~\ref{ss:aux-appr-lemma}) that approximate $\Upsilon$ with a prescribed precision, and then show that for every partition $\mcF\in\mcP_\Theta(I)$ there exists an $\e$-accurate partition $\mcF'\succcurlyeq\mcF$ (Section~\ref{ss:proof-rc}).

\begin{remark}
\label{rmk:sup-prop-cut-bianchini}
    The same strategy, using Prop.~\ref{prop:Qknown-Lyap2} (resp. Remark~\ref{rmk:qknown-lyap})  instead of Theorem~\ref{thm:interaction-estimate}, would show that  the approximation property by Riemann cuts holds as well for the functionals $\Qknown$ (resp. $V+C\Qknown$) introduced in \cite{Bianchini2002, Bianchini2003}:
    \begin{align*}
    \Qknown(u) & = \sup \Big\{ \Qknown \big(\Rc \mcF \bullet(u) \big) \colon \mcF \in \mcP_\Theta(I)
    \Big\}, \\
    (V + C\Qknown)(u) & = \sup \Big\{ (V+C\Qknown) \big(\Rc \mcF \bullet(u) \big) \colon \mcF \in \mcP_\Theta(I) 
    \Big\},
\end{align*}
and therefore also that $\Qknown$ and $V+C\Qknown$ are decreasing across Riemann cuts. 
\end{remark}

\subsection{Variation of \texorpdfstring{$\Upsilon$}{Upsilon} across small Riemann cuts}
\label{ss:direct-est}
Throughout this section, we fix an interval $I\subset \R$, $\hat x\in I$, and a function $u\in\mcD$. Put $\hat u =\Rc{I}{\hat x}(u)$ and denote $\mu_k$ and $\hat\mu_k$ the $k$-th wave measures associated with $u$ and $\hat u$ respectively. By construction, their restrictions on $\R\setminus I$ coincide, and, moreover,  it is straightforward to check that
    \begin{align}
        \hat \mu_k = \mu_k -\mu_k\restr I+ s_k(\hat u_-(\hat x), \hat u_+(\hat x))\delta_{\hat x}.\label{eq:hat-mu-k}
    \end{align}
    This section is devoted to the estimate of the quantity $\abs{\Upsilon(\hat u)-\Upsilon(u)}$ in terms of quantities $\abs{u_x}(I)$, $\mu_k(I)$ and $\abs{\mu_k}(I)$ that are assumed small. More precisely, we aim at proving the following proposition.
\begin{prop}\label{l:alt-appr-Upsilon}
Let $u\in \BV(\R;\Beone)$ and $\hat u =\Rc I{\hat x}$, then
\begin{align*}
    \abs{\Upsilon(\hat u) -\Upsilon(u)}
    \leqslant \const \Big(\abs{u_x}(I)^2+\sum_{k=1}^n\Big[\abs{\mu_k}(I) - \abs{\mu_k(I)}\Big]\Big) \Big(1+\TV(u)^2\Big).
\end{align*}
\end{prop}
The strategy is to estimate the variation of $V_k$ in terms of wave measures (Lemma~\ref{l:alt-appr-Vk}), and then apply Prop.~\ref{prop:u-v-I}.

\begin{lemma}\label{l:alt-appr-Vk} Let $u\in \BV(\R;\Beone)$ and $\hat u =\Rc I{\hat x}$, then
\begin{align*}
    \abs{V_k(\hat u) - V_k(u)}\leqslant  \abs{\mu_k}(I) -\abs{\mu_k(I)}+\const\abs{u_x}(I)^2 .
\end{align*}
\end{lemma}
\begin{proof}
    By~\eqref{eq:hat-mu-k},
    \begin{align}
        V_k(\hat u) - V_k(u) =\abs{\hat \mu_k}(\R)-  \abs{\mu_k}(\R)=  \abs{s_k(\hat u_-(\hat x), \hat u_+(\hat x))}- \abs{\mu_k}(I) .\label{eq:V-k-comp}
    \end{align}
    Then, we have
    \begin{align*}
        \mu_k^\cont(I) = \int_I\scal{\ell_k(u(x))}{u_x^\cont(dx)} = \scal{\ell_k(\hat u_-(\hat x))}{u_x^\cont(I)}+\const \cdot\abs {u_x}(I)^2.
    \end{align*}
    Finally, for every atom $\xi\in I$, by Lemma~\ref{l:length-vs-scalarprod},
    \begin{align*}
        \mu^\atom_k(\{\xi\}) &= \scal{\ell_k(u_-(\xi))}{u_x(\{\xi\})}+\const\abs{u_x(\{\xi\})}^2 \\
        &=\scal{\ell_k\left(\hat u_-(\hat x)\right)}{u_x(\{\xi\})}+\const \abs{u_x}(\{\xi\})\cdot\abs{u_x}(I).
    \end{align*}
    In total, adding $\mu_k^\cont(I)$ and $\mu^\atom_k(\{\xi\}) $ for every atom $\xi$, we obtain
    \begin{align}
        \mu_k(I) = \scal{\ell_k(\hat u_-(\hat x))}{u_x(I)} +\const \abs{u_x}(I)^2.\label{eq:mu_k(I)}
    \end{align}
    Again by Lemma~\ref{l:length-vs-scalarprod},
    \begin{align}
        s_k(\hat u_-(\hat x), \hat u_+(\hat x)) = \scal{\ell_k(\hat u_-(\hat x))}{\smash{\overbrace{\hat u_+(\hat x)-\hat u_-(\hat x)}^{u_x(I)}}}+\const\smash{\overbrace{\abs{\hat u_+(\hat x)-\hat u_-(\hat x)}^2}^{\leqslant\abs{u_x}(I)^2}}.\label{eq:s-k-hat-u}
    \end{align}
    Comparing~\eqref{eq:V-k-comp},~\eqref{eq:mu_k(I)} and~\eqref{eq:s-k-hat-u}, the conclusion follows.
\end{proof}

\begin{proof}[Proof of Prop.~\ref{l:alt-appr-Upsilon}]
 Combine Lemma~\ref{l:alt-appr-Vk} with Prop.~\ref{prop:u-v-I}.
\end{proof}

\subsection{Variation of \texorpdfstring{$\Upsilon$}{Upsilon} across Riemann cuts around big jumps}
\label{ss:punct-I}

\begin{prop}\label{l:alt-appr-Upsilon-2} Let $\beta\ll1$, $u\in \mcD_\beta$ and $\hat u =\Rc I{\hat x}(u)$. Let $z \in I$. Then
    \begin{align*}
    \abs{\Upsilon(\hat u) -\Upsilon(u)}\leqslant \const\Big(\abs{u_x}(J)+\abs{u_x}(J)^2 \Big)\Big(1+\TV(u)^2\Big),\quad J := I\setminus\{z\}.
\end{align*}
\end{prop}
Recall that we threat the jump of~$u$ at~$z$ as \emph{big jump}. That means that this jump dominates the remaining total variation of~$u$ on $I\setminus \{z\} = J$ and~$\hat u$ is obtained from~$u$ by reducing this \emph{big jump} and remaining total variation to a single jump. This motivates the estimate in Prop.~\ref{l:alt-appr-Upsilon-2}  to be in terms of $J$.
\begin{proof}
W.l.o.g. we may assume $\hat x = z$.

Let $ \Jleft = (-\infty, z)\cap J\supseteq\varnothing$, $\Jright = (z, +\infty)\cap J\supseteq\varnothing$, so that $J = \Jleft\sqcup \Jright$. Notice that \begin{align*}
    \hat u = \Rc I z(u) = (\Rc{\Jright\sqcup\{z\}}z\circ\Rc{\Jleft\sqcup\{z\}}{z})(u)
\end{align*}
so it is enough  now to prove the statement in the cases $I = \Jleft\sqcup\{z\},  \Jright\sqcup\{z\}$. 

From now on, we assume then $I = \Jleft\sqcup\{z\}$ (the remaining case being completely analogous). In particular, $J = \Jleft$. If $J = \varnothing$, $u = \hat u$, hence there is nothing to prove, so assume $J \neq \varnothing$.

 Pick $p,p'\in J$, $p<p'$ and define a bounded variation, right continuous function  $\mathring u$ such that $\mathring u (-\infty) = u(-\infty)$ and
    $
        \mathring u_x = \hat u_x + u_x(J)[\delta_{p}-\delta_{p'}]
   $.
   We have $\Rc{J}{p}(u) = \Rc{(p,z]}{z}(\mathring u)$. Indeed,
   \begin{align*}
        &(\Rc{(p,z]}{z}(\mathring u))_x = \mathring u_x-\overbrace{\mathring u_x\restr {(p,z]}}^{\mathclap{u_x(I)\delta_z-u_x(J)\delta_{ p'}\hspace{4mm}}}+\overbrace{\mathring u_x((p,z])}^{\hspace{2mm}u_x(\{z\})}\delta_z 
        =\hat u_x+u_x(J)\delta_p-u_x(J)\delta_z\\
        &=u_x-u_x\restr {I}+u_x(\{z\})\delta_z+u_x(J)\delta_p=u_x-u_x\restr {J}+u_x(J)\delta_p = (\Rc{J}{p}(u))_x.
    \end{align*}
    We now claim that
    \begin{align}\label{l:alt-appr-Upsilon-2:eq0.1}
        \Upsilon(\hat u)\leqslant\Upsilon(\Rc Jp(u))\leqslant\Upsilon(\mathring u).
    \end{align}
    Indeed,
    \begin{align*}
        \Upsilon(\hat u)&\smash{\stackrel{\eqref{eq:lyap-prop}}\leqslant\Upsilon(\overbrace{\Rc {\big\{J, \{z\}\big\}}{\{p,z\}}(u)}^{\Rc J p(u)})}
        =\Upsilon(\Rc{(p,z]}{z}(\mathring u))\\
        &\stackrel{\eqref{eq:lyap-prop}}\leqslant\Upsilon(\Rc{\big\{(p,p'), [p', z),\{z\}\big\}}{\big\{\tfrac{p+p'}2, p',z\big\}}(\mathring u)) =\Upsilon(\mathring u).
    \end{align*}
    Now one has
    \begin{align*}
        \abs{\Upsilon(u)-\Upsilon(\hat u)}\leqslant\abs{\Upsilon( u )-\Upsilon(\Rc Jp(u))}+\abs{\Upsilon(\Rc Jp(u))-\Upsilon(\hat u)}&\\
        \stackrel{\eqref{l:alt-appr-Upsilon-2:eq0.1}}\leqslant\abs{\Upsilon(u )-\Upsilon(\Rc Jp(u))}+&\smash{|\Upsilon(\mathring u)-\Upsilon(
            \stackrel{\substack{\mathclap{\Rc Jp(\mathring u)}\\ \rotatebox{90}=} }{\hat u}
        )|}
    \end{align*}
    Both terms in the r.h.s. can be dominated by $\const\Big(\abs{u_x}(J)+\abs{u_x}(J)^2 \Big)\Big(1+\TV(u)^2\Big)$ using Prop.~\ref{l:alt-appr-Upsilon}.
\end{proof}
    \begin{figure}[h]
        \centering
         \includegraphics{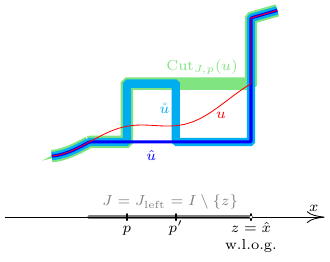}\\
          \caption{To Prop.~\ref{l:alt-appr-Upsilon-2}: the case $I = \Jleft\sqcup\{z\}$.}
        \label{fig:check-u-mr-u}
    \end{figure}
\begin{remark}
\label{remark:shrinking}
    There is a small gap in the proof of Prop.~\ref{l:alt-appr-Upsilon-2}. Namely, $\mathring u$ does not necessary belongs to $\mcD_\beta$ a priori. In particular,~\eqref{eq:lyap-prop} would not be applicable. However, $\TV(\mathring u)\leqslant 3\TV(\hat u)$. So, replacing  $\beta$ in~\eqref{eq:def-D-beta} by \indfixbeta{$\beta'<\beta$ small enough}, we may guarantee that $u\in\mcD_{\beta'}$ implies $v\in\mcD_{\beta} $ for every $v$ with $\TV(v)\leqslant3\TV(u)$ and $v(-\infty) = u(-\infty)$. 
\end{remark}

\subsection{\texorpdfstring{$\e$}{epsilon}-accurate partitions}
\label{ss:aux-appr-lemma} 

We now introduce a class of \emph{$\e$-accurate partitions}  that approximate $\Upsilon$ with a prescribed precision.

\begin{definition}
    Let $u \in \BV(\R; \Beone)$, $\e>0$, $I \subseteq \R$ be an interval. We say that $\mcF \in \mcP(I)$ is \emph{$\e$-accurate} (for $u$ on $I$) if one can write $\mcF = \mcB \sqcup \mcS$ s.t. 
        \begin{enumerate}[(I)]
        \item\label{a2} every $I\in \mcB$ contains exactly one point $z_I$ with $\abs{u_x}(\{z_I\})\geqslant\e$ and, moreover, $\abs{u_x}(J_I)\leqslant \e^2$, where $J_I:=I\setminus \{z_I\}$;
        \item\label{a3} for every $I\in \mcS$, $\abs{u_x}(I)\leqslant\varepsilon$; 
        \item \label{a4} for every $k=1,\dots, n$, $ \sum_{I \in \mcF}  \big( |\mu_k|(I) - |\mu_k(I)| \big) \leqslant \e$.  
    \end{enumerate}
\end{definition}
\noindent We treat $z_I$ in~(\ref{a2}) as the locations of \emph{big jumps} of $u$ while the intervals from $\mcS$ as the ones with \emph{small} variation of $u$ that corresponds to the frameworks of Prop.~\ref{l:alt-appr-Upsilon-2} and Prop.~\ref{l:alt-appr-Upsilon}, respectively. Notice also that $\card \mcB \leqslant \TV(u) \e^{-1}$ and
\begin{align}\label{eq:absorption}
    \mcF\text{ is $\e$-accurate, }\mcF'\succcurlyeq\mcF\quad\Longrightarrow\quad \mcF'\text{ is $\e$-accurate.}
\end{align}

The following lemma shows that $\e$-accurate partitions approximate $\Upsilon$ with prescribed precision.
\begin{lemma}\label{lemma:aux-appr}
    Let $u\in\mcD$, $I \subseteq\R$ be an interval, $\mcF \in \mcP(I)$ be $\e$-accurate. Then
    \begin{align*}
        \abs{\Upsilon( \Rc{\mcF}{\bullet})-\Upsilon(u)}\leqslant\const(\e+\e^3)(1+\TV(u)^3).
    \end{align*}
    
\end{lemma}

\begin{proof}
    Applying Prop.~\ref{l:alt-appr-Upsilon} (resp. Prop.~\ref{l:alt-appr-Upsilon-2}) to the Riemann cuts $\Rc I{\hat x}$, $I\in\mcS$ (resp. $I\in \mcB$), we get
    \begin{align*}
        \abs{\Upsilon(\Rc \mcF\bullet(u))-\Upsilon(u)}\leqslant\const\Big(1+\TV(u)^2\Big)\Big(\sum_{I\in\mcB}\Big[\abs{u_x}(J_I)+\abs{u_x}(J_I)^2\Big]\\
        +\sum_{I\in\mcS}\Big[\abs{u_x}(I)^2+{\sum_{k}\big(\abs{\mu_k}(I) - \abs{\mu_k(I)}\big)\Big]} 
        \Big).
    \end{align*}
    Notice that the statement now follows since
    \begin{itemize}
        \item $\sum_{I\in\mcB}\Big[\abs{u_x}(J_I)+\abs{u_x}(J_I)^2\Big]\leqslant(\e^2+\e^4)\card\mcB\leqslant(\e+\e^3)\TV(u)$ by the assertion~(\ref{a2}); 
        \item $\sum_{I\in\mcS}\abs{u_x}(I)^2\leqslant\e\sum_{I\in\mcS}\abs{u_x}(I)\leqslant\e\TV(u)$ by the assertion~(\ref{a3});
        \item $\sum_{I\in\mcS}\sum_{k}\big(\abs{\mu_k}(I) - \abs{\mu_k(I)}\big)\leqslant n\e$ by the assertion~(\ref{a4}).\qedhere
    \end{itemize}
\end{proof}

The following lemma claims the existence of $\e$-accurate partitions, for every $\e>0$.
\begin{lemma}\label{lemma:FFrpime}
     Let $\e>0$,  $u\in\BV(\R; \Beone)$, $I$ is an interval and $\Theta\subset I$ is dense in $I$. Then there is an $\e$-accurate partition $\mcF\in\mcP_\Theta(I)$.
\end{lemma}
\begin{proof}
    It is straightforward to construct $\tilde \mcF\in\mcP_\Theta(I)$, satisfying the assertions~(\ref{a2})--(\ref{a3}). 

    Recall now that $\sgn\mu_k\in L^1(I,\abs{\mu_k})$ denotes the Radon--Nikodym derivative $d\mu_k/d\abs{\mu_k}$. Notice that $\Theta$ is dense in $I$ and $C_c(I)$ is dense in $ L^1(I,\abs{\mu_k})$. Thus, we can find $\mcG\in\mcP_\Theta(I)$ such that the distance from the function $\sgn\mu_k $ to the subspace $\vspan\{\ind{J}\colon J\in\mcG\}$ in $L^1(I,\abs{\mu_k})$ is less than $\e$ for every $k= 1,\ldots, n$.

    Observe that $\psi_k:=
        \sum_{J\in\mcG}\ind J\sgn(\mu(J))$
    minimizes the $L^1(I,\abs{\mu_k})$-distance between $\sgn \mu_k$ and $\vspan\{\ind{J}\colon J\in\mcG\}$. Therefore, 
    \begin{align*}
        \sum_{J\in\mcG}\big(\abs{\mu_k}(J) - \abs{\mu_k(J)}\big) = {\int\limits_{I} }({\sgn\mu_k - \psi_k})\, d\mu_k\leqslant \norm{\sgn\mu_k-\psi_k}_{L^1(I,\abs{\mu_k})}\leqslant\e,
    \end{align*}
    so that the assertion~\eqref{a4} is satisfied by $\mcG$. Putting $\mcF = \{\tilde J\cap J\colon\tilde J \in\tilde \mcF, J\in\mcG\}$ concludes the proof.
\end{proof}

\subsection{Proofs of  approximation property w.r.t. to Riemann cuts (Theorem~\ref{thm:sup-prop-cut})}
\label{ss:proof-rc}

\begin{proof}[Proof of Theorem~\ref{thm:sup-prop-cut}]
    Let $\mcF\in\mcP_\Theta(I)$ and $\e>0$. By Lemma~\ref{lemma:FFrpime}, construct an $\e$-accurate partition $\mcF_\e \in \mcP_\Theta(I)$, then $\mcF' = \{J\cap J_\e\colon J\in\mcF, J_\e\in\mcF_\e\} \in \mcP_\Theta(I)$ is also $\e$-accurate by~\eqref{eq:absorption}. Then, $\abs{\Upsilon(\Rc{\mcF'}\bullet)-\Upsilon(u)}\leqslant\const(\e+\e^3)(1+\TV(u)^3)$ due to Lemma~\ref{lemma:aux-appr}. By~\eqref{eq:lyap-prop}, $\Upsilon(\Rc\mcF\bullet)\leqslant\Upsilon(\Rc{\mcF'}\bullet)$.
    
    Thus, up to a possible shrinking of $\e$, we have proved that
    \begin{align*}
    \forall \mcF\in\mcP_\Theta(I),\e>0\colon&\exists \mcF'\in\mcP_\Theta(I)\colon \\
    &\Upsilon(\Rc\mcF\bullet)\leqslant\Upsilon(\Rc{\mcF'}\bullet)\in(\Upsilon(u)-\e,\Upsilon(u)+\e).
\end{align*}
    This implies the statement.
\end{proof}

\section{Stability w.r.t. perturbations in state space}
\label{s:cont-wrt-states}

We now state  and prove rigorously \ref{thm-main-b}, item~\eqref{thm-main:stab}. 

\begin{theorem}
\label{thm:cont-wrt-states} 
    Let $\beta \ll 1$. Let $u \in \mcD_\beta$ and let $\Phi\colon \Omega \to \R^n$ be Lipschitz continuous with $\Phi \circ u \in \mcD_\beta$. Then
         \begin{align*}
            \abs{\Upsilon(u)-\Upsilon(\Phi(u))}= 
            \const(1+M^2)M\norm{\Phi-\id}_{W^{1,\infty}},\\
            M := \max\{\TV(u),\TV(\Phi\circ u)\}\leqslant \TV(u)(1+\norm{D\Phi}_{L^\infty}),
        \end{align*}
        where $\id(p) = p$ for all $p \in \R^n$.
\end{theorem}

We first prove an auxiliary Lemma~\ref{lemma:par-trans}. Then we show Prop.~\ref{prop:perturb-space-2}, measuring the change of $\Upsilon$ w.r.t. the perturbation of piecewise constant function. Finally we conclude the proof of Theorem~\ref{thm:cont-wrt-states}.
\begin{lemma}\label{lemma:par-trans}
    Let $u^L, u^R, v^L, v^R \in \Beone$, with 
    \begin{equation*}
        v^L = u^L + p, \qquad v^R = u^R + p
    \end{equation*}
    for some $p \in \R^n$. Then
    \begin{align}
        \abs{s_k(v^L
    , v^R)-s_k(u^L
    , u^R)} \leqslant \const\abs{p}\cdot\abs{u^R-u^L}.\label{lemma:par-trans:eq1}
    \end{align}
\end{lemma}
\begin{proof}
    Denote 
    $\xi = u^R-u^L$. By Prop.~\ref{prop:sdef-quadr-bound} one has 
    \begin{align*}
        s_k(u^L, v^R) &= s_k(u^L,u^R)+s_k(u^R,v^R)+\const \abs {p} \cdot\abs\xi\\
        &= s_k(u^L,v^L)+s_k(v^L,v^R)+\const \abs {p} \cdot\abs\xi,
    \end{align*}
    hence
    \begin{align*}
        &{s_k(v^L, v^R)-s_k(u^L,u^R)}  ={s_k(u^R, v^R)-s_k(u^L, v^L )}+\const \abs {p} \cdot\abs\xi  \\
        &\xlongequal{\text{ Lemma~\ref{l:length-vs-scalarprod}}}
        \scal{\ell_k(u^R)-\ell_k(u^L)}{p }+\const\big(\abs{p}^2+ \abs {p} \cdot\abs\xi\big)=\const\big(\abs{p}^2+ \abs {p} \cdot\abs\xi\big).
    \end{align*}
    Now consider the function  $p\mapsto s_k(u^L+p, u^R+p)$ with open, convex domain consisting of $p$ such that $u^L+p, u^R+p \in \Beone$. The equality we just have proved implies that this function is locally $\const\abs{\xi}$-Lipschitz at point $p=0$. Similarly, we may prove its local $\const\abs{\xi}$-Lipschitz continuity at any other point inside $\Beone$ that implies~\eqref{lemma:par-trans:eq1} due to the convexity of $\Beone$.
\end{proof}
\begin{prop}\label{prop:perturb-space-2}
    Let $\beta \ll 1$. 
    Let $w, \overline w\in\mcD_\beta$ be piecewise constant functions taking subsequent values $w_0,\ldots, w_m$ and $\overline w_0,\ldots, \overline w_m$.
    Let $A=\max_{i=0,\ldots,m}\abs{\overline w_i-w_i}$, $B=\sum_{i=1}^m
    |(\overline w_i\allowbreak-\overline w_{i-1}) - (w_i-w_{i-1})|$, $M = \max\{\TV(w),\TV(\overline w)\}$.
    Then
    \begin{align}
        \abs{\Upsilon(\overline w)-\Upsilon(w)}\leqslant\const (1+M^2)(AM+B).
        \label{lemma:perturb-space-2:eq1}
    \end{align}
\end{prop}
    \begin{proof}    
        We divide the proof in steps. 
        
        \textit{Step 1 (homotopy argument)}. Put $p_i = \overline w_i-w_i$ and $w_i^\tau = w_i + \tau p_i$, $i=0,\ldots,m$, $\tau\in[0,1]$.
        Construct piecewise constant function $w^\tau$ taking subsequent values $w_0^\tau\ldots w_m^\tau$ depending on $\tau\in[0,1]$ (see Fig.~\ref{fig:w-bar-w}(a)).

        First of all, notice that $\TV(w^\tau)\leqslant M$. Thus, in spirit of Remark~\ref{remark:shrinking}, for $\beta'<\beta$ small enough, $w,\overline w\in \mcD_{\beta'}$ implies $w^\tau\in \mcD_{\beta}$ for all $\tau\in[0,1]$.

        Our goal is to show that, for $|\tau - \tilde \tau| \leqslant \e$, \begin{align}
        \label{eq:pert-state-3}
            |\Upsilon(w^\tau)-\Upsilon(w^{\tilde\tau})|\leqslant \const \e(1+M^2)(AM+B)+\boxed{\text{terms of order $\e^2$}}
        \end{align}
        from which~\eqref{lemma:perturb-space-2:eq1} easily follows. 
        
        \smallskip
        \textit{Step 2 (Auxiliary function $v$)}.
        Fix thus $\e\in(0,1)$ and pick $\tau,\tilde \tau\in[0,1]$ with $\abs{\tilde\tau-\tau}\leqslant\e$. 
        In order to compare $\Upsilon(w^\tau)$ with $\Upsilon(w^{\tilde \tau})$, we construct an auxiliary map $v$, as follows. 
        Let $s_k^i := s_k(w^\tau_{i-1},w^\tau_i)$, $i = 1,\ldots, m$.
        For every $i = 1,\ldots,m$, find $v_{i}'$ such that $s_k(w^{\tilde \tau}_{i-1}, v_i') = s_k^i$, $k=1,\ldots, n$ (see Fig.~\ref{fig:w-bar-w}(b)). The piecewise constant function $v$ is constructed by taking subsequent values $w_0^{\tilde\tau},\allowbreak v_1',\allowbreak w_1^{\tilde\tau},\allowbreak v_2',\allowbreak\ldots,\allowbreak w_m^{\tilde\tau}$ (see Fig.~\ref{fig:w-bar-w}(a)).
        \begin{figure}[h]
            \centering
            \begin{minipage}{0.45\textwidth}
                \centering
                \includegraphics{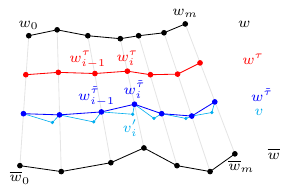}\\
                (a)
            \end{minipage}
            \hfill
            \begin{minipage}{0.45\textwidth}
                \centering
                \includegraphics{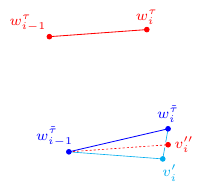}\\
                (b)
            \end{minipage}
            \caption{To the proof of Prop.~\ref{prop:perturb-space-2}. Cyan curve corresponds to $v$.}
            \label{fig:w-bar-w}
        \end{figure}

        Observe that, roughly speaking, jumps in $v$ are of two types: jumps of the form $(w^{\tilde \tau}_{i-1}, v_i')$, which we will treat as \emph{big} jumps (they are translations, in terms of wave strengths, of the jumps of $w^\tau$); and jumps of the form $(v_i', w_i^\tau)$, that we will treat as \emph{small} jumps.

        \smallskip
        \textit{Step 3 (Estimates on small jumps in $v$)}.
        First of all, we show that the sum of all jumps of the second type $\sum_{i=1}^m\abs{v_i'-w^{\tilde \tau}_i}$
        is indeed small (precisely, of order $\e$).
        
        To this aim, let $v_i'' :=w^{\tilde\tau}_{i-1}+({w^\tau_{i}-w^\tau_{i-1}})$ (see Fig.~\ref{fig:w-bar-w}(b)).
        Applying Lemmas~\ref{lemma:3.5}--\ref{lemma:3.6}, we get
        \begin{align*}
            \abs{v'_i-v''_i} 
            & = \const\sum_{k=1}^n \big|\smash{\overbrace{s_k(w^{\tilde \tau}_{i-1}, v_i')}^{s_k(w^\tau_{i-1}, w_i^{\tau})} - \overbrace{s_k(w^{\tilde\tau}_{i-1}, v_i'')}^{\hspace{2.3cm}\mathclap{s_k(w^{\tau}_{i-1} + (\tilde \tau - \tau) p_{i-1}, w^\tau_i + (\tilde \tau - \tau) p_{i-1})}} }\big| \\
            & \xlongequal{\text{Lemma~\ref{lemma:par-trans}}}\const\abs{p_{i-1}}\cdot\e\cdot \abs{{w^\tau_{i}-w^\tau_{i-1}}} = \const A \e \abs{{w^\tau_{i}-w^\tau_{i-1}}}. 
        \end{align*}
       Thus $\sum_{i=1}^m\abs{v'_i-v''_i} = \const A\e M$ and consequently
        \begin{align}
        \label{eq:am+b}
        \sum_{i=1}^m\abs{v_i'-w^{\tilde \tau}_i}\leqslant\sum_{i=1}^n\Big[{\abs{v_i'-v_i''}}+{\abs{v_i''-w^{\tilde \tau}_i}}\Big]=\const \smash{\overbrace{(A M+B)}^{=:C}}\e.
        \end{align}
    Then 
    \begin{align}
        \label{eq:tw+ce}
        \TV(v) \leqslant \TV(w^{\tilde \tau})+\const C\e,
    \end{align}
    and so, for $\beta'<\beta$ and $\e$ small enough, $w,\overline w\in \mcD_{\beta'}$ implies $ v\in\mcD_\beta $ (independently of choice of $\tau, \tilde\tau)$.

      \smallskip
        \textit{Step 4 (Comparison of $\Upsilon(w^\tau)$ with $\Upsilon(v)$)}.
      We now compare $\Upsilon(w^\tau)$ with $\Upsilon(v)$ (whilst later  we will compare $\Upsilon(v)$ with $\Upsilon(w^{\tilde \tau})$). We argue separately for the four functionals $V, \Qtrans, \Qcubic, \Qcost$. 
      We have $V(v)\geqslant V(w^\tau)$ and $\Qtrans(v)\geqslant\Qtrans(w^\tau)$, since $v$ ``contains'' all the jumps of $w^\tau$ (in terms of wave strengths). Next, denote $\sigma_k^i, \tilde\sigma_{k}^i\colon[0,|{s_k^i}|]\to\R $ the $k$-th speed maps associated with Riemann problems $(w^{\tau}_{i-1}, w^\tau_{i})$ and $(w^{\tilde \tau}_{i-1}, v_i')$, respectively.
        We have
        \begin{align*}
        \rule{0pt}{.6cm}\Qcubic_k(v)\geqslant&\sum_{i,j=1}^m\int_0^{|{s_k^i}|}\int_0^{|{s_k^j}|}\abs{\tilde\sigma_k^i(s) - \tilde\sigma_k^j(s')}\,ds' ds\geqslant - 2\sum_{i,j=1}^m|s_k^i|\cdot\|\sigma_k^j-\tilde\sigma_k^j\|_{L^1}\\
        +&\sum_{i,j = 1}^m\int_0^{|{s_k^i}|}\int_0^{|{s_k^j}|}\abs{\sigma_k^i(s) - \sigma_k^j(s')}\,ds' ds,
        \end{align*}
        the latter term being equal to $\Qcubic_k(w^\tau)$. By Lemma~\ref{lemma:3.5}, $\|\sigma_k^j-\tilde\sigma_k^j\|_{L^1}\leqslant \const |s_k^j|\cdot|w^{\tilde\tau}_{j-1}-w^\tau_{j-1}|\leqslant\const A\e|s_k^j|$, and thus
        \begin{align*}
            \Qcubic_k(v)\geqslant\Qcubic_k(w^\tau) -\const\sum_{i,j=1}^m|s_k^i|\cdot A\e|s_k^j| \geqslant\Qcubic_k(w^\tau)-\const\TV(w^\tau)^2A\e.
        \end{align*}
        Finally, denote $f_k^i, \tilde f_k^{i}\colon \I(s_k^i)\to\R$ the reduced fluxes associated with Riemann problems $(w^\tau_{i-1}, w^\tau_{i})$ and $(w^{\tilde\tau}_{i-1}, v_i')$, respectively. We have, using the quadratic bound and stability of $\cost$ (Propositions~\ref{prop:cost_Linfty_est-recall} and~\ref{prop:stability-recall}),        
        \begin{align*}
            &\Qcost_k(v )\leqslant\const\sum_{i=1}^m\abs{v_i'-w^{\tilde\tau}_i}^2+\sum_{i=1}^m\cost_{\tilde  f_k}(0,s_k^i) \stackrel{\eqref{eq:am+b}}{\leqslant}\const C^2\e^2+\sum_{i=1}^m\cost_{ f_k}(0,s_k^i) \\
            &+\const\sum_{i=1}^m|s_k^i|\cdot\norm{\derder {f_k^i}-\derder{\tilde  f_{k}^i}}_{L^1}\hspace{.8cm}\stackrel{\mathclap{\text{Lemma~\ref{lemma:3.5}}}}\leqslant\hspace{.8cm}\const C^2\e^2+\Qcost_k(w^\tau)\\
            &+\const\sum_{i=1}^m|s_k^i|^2\cdot|w^{\tilde\tau}_{i-1}-w^\tau_{i-1}|
            \leqslant\Qcost_k(w^\tau) +\const(C^2\e^2+\smash{\underbrace{\TV(w^\tau)^2}_{\leqslant M^2}} A\e).
        \end{align*}
        Gathering everything together, we obtain
        \begin{align}
        \label{eq:pert-state1}
            \Upsilon(v)\geqslant\Upsilon(w^\tau) - \const(
                           M^2 A\e +C^2\e^2
            )          
        \end{align}

        \smallskip
        \textit{Step 5 (Comparison of $\Upsilon(v)$ with $\Upsilon(w^{\tilde \tau})$)}.
        Let us now compare $\Upsilon(v)$ with $\Upsilon(w^{\tilde \tau})$. Notice that $w^{\tilde\tau}$ is obtained from $v$ by a sequence of Riemann cuts (up to a change of jump locations).  and  so, by Prop.~\ref{l:alt-appr-Upsilon-2},
        \begin{align}
            \Upsilon(v)-\Upsilon(w^{\tilde\tau})&\leqslant\const(1+\TV(v)^2)\sum_{i=1}^m\Big(\abs{v_i'-w^{\tilde \tau}_i}+\abs{v_i'-w^{\tilde\tau}_i}^2\Big) \nonumber \\
            &\stackrel{\mathclap{\eqref{eq:am+b}, \eqref{eq:tw+ce}}}{\leqslant}\hspace{3mm}\const(1+\const[M+C\e]^2)(C\e+C^2\e^2) \nonumber \\
            & \leqslant \const (1 + M^2)C\e + \e^2 \times \boxed{\text{polynomial in $M,C$}}.\label{eq:pert-state-2}
        \end{align}

        \smallskip
        \textit{Step 6 (Comparison of $\Upsilon(w^\tau)$ with $\Upsilon(w^{\tilde \tau})$)}.
        Putting together \eqref{eq:pert-state1},~\eqref{eq:pert-state-2}, we get 
        \begin{align*}     \Upsilon(w^\tau)-\Upsilon(w^{\tilde\tau})\leqslant \const \e(1+M^2)(AM+B)+\boxed{\text{terms of order $\e^2$}}.
         \end{align*}
    Writing an analogous inequality for swapped $w^\tau$ and $w^{\tilde \tau}$, we obtain~\eqref{eq:pert-state-3}, and this concludes the proof. 
    \end{proof}
    
    \begin{proof}[Proof of Theorem~\ref{thm:cont-wrt-states}]
        Notice that the statement holds true for piecewise constant $u$, due to Prop.~\ref{prop:perturb-space-2}: indeed, putting $w = u$, $\overline{w} = \Phi\circ u$, we have $A\leqslant \norm{\Phi-\id}_{L^\infty}$,  $B\leqslant \norm{D\Phi-\id_n}_{L^\infty} M $.

    For an arbitrary $u\in \mcD$, for every $\mcF\in\mcP(\R)$ and $X\subset \R$ compatible with~$\mcF$ we have
    \begin{align}\label{eq:Phi-comm}
        \Phi\circ\Rc \mcF X (u) = \Rc{\mcF} X(\Phi\circ u).
    \end{align}
    Now, by Theorem.~\ref{thm:sup-prop-cut} and~\eqref{eq:lyap-prop}, we may take fine enough partition $\mcF\in\mcP_\R(\R)$ to approximate $\Upsilon(u)$ and $\Upsilon(\Phi\circ u)$ by $\Upsilon(\Rc\mcF X(u))$  and $\Upsilon(\Rc\mcF X(\Phi\circ u))$, respectively. Estimating the difference of the two last quantities by the piecewise constant case and~\eqref{eq:Phi-comm}, we conclude.
    \end{proof}
  
\section{Further properties of \texorpdfstring{$\Upsilon$}{Upsilon}}
\label{s:further-properties}

We now state  and prove rigorously \ref{thm-main-a}, item~\eqref{thm-main:lsc} (Theorem~\ref{thm:lsc}), item~\eqref{thm-main:lyap} (Theorem~\ref{thm:decr-along-semigroup}), and item~\eqref{thm-main:invar} (Theorem~\ref{thm:reparametr-state-curve}). 

\subsection{Lower semicontinuity}

\label{ss:lsc}

\begin{theorem}[Lower semicontinuity in $\Loneloc$]
    \label{thm:lsc}
   Let $\beta \ll 1$. If $(u_m)_m$ and $u$ are in~$\mcD_\beta$ and $u_m \to u$ in $\Loneloc$, then
        \begin{equation*}
            \Upsilon(u) \leqslant \liminf_{m \to \infty} \Upsilon(u_m).
        \end{equation*}
\end{theorem}

\begin{proof}
    Possibly passing to a subsequence, we may assume that $\liminf_m \Upsilon(u_m) = \lim_m \Upsilon(u_m)$ and $u_m\to u$ pointwise on some full measure set. Let $\Theta$ be the intersection of this set with the set of common continuity points of $\{u_m\}_m$ and $u$. Then $\Theta$ is a full measure subset of $\R$, and thus it is dense. 

    Pick an arbitrary interval $I$ with $\inf I,\sup I\in\Theta$, and arbitrary partition $\mcF = \{I_1<\ldots<I_N\}\in\mcP_\Theta(I) $. Let $X$ be a finite set of points compatible with $\mcF$. By construction and continuity of $u$ at the points of $\Theta$, the function $(\Rc \mcF X \circ \Cr I )(u)$ is piecewise constant and takes subsequent values
    \begin{align}\label{eq:crop-cut-u}
        u(\inf I), u(\sup I_1),\ldots, u(\sup I_N),
    \end{align}
    see Fig.~\ref{fig:crop_cut}.
    \begin{figure}[h]
        \centering
        \includegraphics{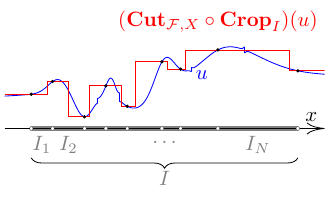}
        \label{fig:crop_cut}
        \caption{The values of $(\Rc \mcF X \circ \Cr I )(u)$.}
    \end{figure}
    Similarly, for every $m$, the function $(\Rc \mcF X \circ \Cr I )(u_m)$ is piecewise constant and takes subsequent values
    \begin{align}\label{eq:crop-cut-um}
        u_m(\inf I), u_m(\sup I_1),\ldots, u_m(\sup I_N).
    \end{align}
    By pointwise convergence of $u_m$ to $u$ on $\Theta$, the values~\eqref{eq:crop-cut-um} converge to the values~\eqref{eq:crop-cut-u}. Hence, applying Theorem~\ref{thm:sup-prop-cut}, Theorem~\ref{prop:sup-prop-crop} and Prop.~\ref{prop:perturb-space-2}, we obtain 
    \begin{align*}
        \Upsilon(u_m)\geqslant\Upsilon\Big((\Rc \mcF X \circ \Cr I )(u_m)\Big)\to \Upsilon\Big((\Rc \mcF X \circ \Cr I )(u)\Big)
    \end{align*}
    so that 
    \begin{align*}
        \liminf_m\Upsilon(u_m)\geqslant\sup_{\substack{I \subseteq \R \\ \inf I,  \sup I \in \Theta}} \sup_{\mcF\in\mcP_\Theta(I)}\Upsilon\Big((\Rc \mcF X \circ \Cr I )(u)\Big).
    \end{align*}
    Again by Theorem~\ref{thm:sup-prop-cut} and Theorem~\ref{prop:sup-prop-crop}, the r.h.s. is equal to $\Upsilon(u)$ that concludes the proof.
\end{proof}

\begin{remark}
\label{rmk:lsc-bianchini}
    The same proof shows that the functionals $\Qknown$, $V+C\Qknown$ of \cite{Bianchini2002, Bianchini2003} are l.s.c., as indeed claimed in \cite{Bianchini2003}. 
\end{remark}

\subsection{Lyapunov property}
\label{ss:decr-along-semigroup}

\begin{theorem}
\label{thm:decr-along-semigroup} 
Let $\beta \ll 1$. Let $\bar u \in \mcD_\beta$ and let $u(t, \cdot) := S_t \bar u$, $t\in [0, + \infty)$, be the vanishing viscosity solution with initial value $\bar u$. Then $S_t \bar u \in \mcD_\beta$ and $$t \mapsto \Upsilon(u(t))$$ is non-increasing. 
\end{theorem}

\begin{proof}
By semigroup property, it is enough  to prove $\Upsilon(\overline u)\geqslant\Upsilon(u(t,\cdot))$, $t>0$.
    We are going to apply the Glimm scheme~\cite{Glimm1965, BressanMarson, Liu1977, BM-Glimm}. \indfixbeta{Assume $\beta \ll 1$.} One can construct  a sequence $(u_m)_m$ of Glimm approximate solutions to~\eqref{eq:syst-cons-law}, with initial data $u_m(0,\cdot) = \overline u_m$.
    We can also assume (using respectively Theorem~\ref{thm:sup-prop-cut} and Remark~\ref{rmk:sup-prop-cut-bianchini}) that   
    $\Upsilon(\overline u_m)\leqslant\Upsilon(\overline u)$ and, similarly, $V(\overline u_m) + C\Qknown(\overline u_m) \leqslant V(\overline u) + C \Qknown (\overline u)$. This implies, by l.s.c. (Theorem~\ref{thm:lsc} and Remark~\ref{rmk:lsc-bianchini} respectively) that $V(\overline u_m) + C\Qknown(\overline u_m) \to V(\overline u) + C \Qknown (\overline u)$ and $\Upsilon(\overline u_m)\to\Upsilon(\overline u)$.  
    
    Observe that each Glimm approximate solution $u_m$ is constructed from the initial datum $\overline u_m$ by iteratively applying Riemann cuts and solving local Riemann problems. Therefore both maps   $t \mapsto \Upsilon(u_m(t,\cdot))$, $t\mapsto V(u_m(t,\cdot)) + C\Qknown(u_m(t,\cdot))$ are non-increasing by Theorem~\ref{thm:interaction-estimate} (Riemann cuts) and Theorem~\ref{thm:const-if-no-interaction} (persistency) (resp. Remark~\ref{rmk:qknown-lyap} and Remark~\ref{rmk:persistencty-bianchini}). In particular, recalling the definition of $\mcD_\beta$ in~\eqref{eq:def-D-beta}, $u_m(t,\cdot) \in \mcD_\beta$, for all $t$ and all $m$.  

    By Helly's Theorem, $u_m \to \tilde u$ (up to not relabeled subsequences). By l.s.c. of $V+C \Qknown$ (Remark~\ref{rmk:lsc-bianchini}) we have $V(\tilde u(t,\cdot)) + C\Qknown(\tilde u(t,\cdot)) \leqslant \beta$ for all~$t$. By tightness of the family of measures $\{\abs{(u_m(t,\cdot))_x}\}_m$ (due to finite speed of propagation), we also get  $\tilde u(t,-\infty) = \overline u(-\infty)$. Overall $\tilde u(t,\cdot) \in \mcD_\beta$, for all~$t$. 

    On the other side, by the analysis in \cite{BM-Glimm} (see also \cite{AM2011}), the limit $\tilde u$ of Glimm scheme must coincide with the vanishing viscosity solution $u(t,\cdot) = S_t \bar u$ to \eqref{eq:syst-cons-law}--\eqref{eq:initial-datum}, for all $t$. Hence $S_t \bar u \in \mcD_\beta$ and, by l.s.c. of $\Upsilon$ (Theorem~\ref{thm:lsc}), we also get $\Upsilon(u(t,\cdot)) \leqslant \Upsilon(\overline u)$. 
\end{proof}
\begin{remark}
    The consistency of Glimm scheme without GNL/LD assumptions was shown in~\cite{BM-Glimm} by means of non-local quadratic functional $\mfQ(t)$ discussed in Section~\ref{ss:intro-locality-time}. Alternatively, one might prove the consistency using our local-in-time functional $\Upsilon(u^\e(t))$ applied to Glimm scheme approximations $u^\e(t)$. No circular reasoning would occur in the proof of Theorem~\ref{thm:decr-along-semigroup}, since~\cite{BM-Glimm} only requires functional to be defined on approximate, piecewise-constant solutions, and thus only Theorems~\ref{thm:interaction-estimate},~\ref{thm:const-if-no-interaction}, that we have already proven. 
\end{remark}

\subsection{Invariance under reparametrization of independent variable}
\label{ss:reparametr-state-curve}

\begin{theorem}
\label{thm:reparametr-state-curve} 
    Let $u\in\BV(\R;\Beone)$ and $\psi\colon \R\to\R$ be strictly increasing. Let $\tilde u\in\BV(\R;\Beone)$ be such that $\tilde u (-\infty) = u(-\infty)$ and $\tilde u_x = \psi_{\sharp}u_x$. Then $\Upsilon(\tilde u) = \Upsilon(u)$.
\end{theorem}
Notice that the pair of conditions $\tilde u(-\infty) = u(-\infty)$ and $\tilde u_x = \psi_\sharp u_x$ is equivalent to $\tilde u(y) = u(\psi^{-1}(y))$ for smooth $u$. We prefer to use the push-forward notation since we assume only $u \in \BV$.

\begin{proof}
    Denote~$\mu_k$ (resp.~$\tilde \mu_k$) the $k$-th wave measure of $u$ (resp. $\tilde u$). Since $\psi$ is injective, we have $(\tilde u _x)^\atom = \psi_\sharp[u_x^\atom]$ and hence $(\tilde u _x)^\cont = \psi_\sharp[u_x^\cont]$. This, by definition of $k$-th wave measure (see Section~\ref{ss:def-V-k-Qknown}), implies $\tilde\mu_k = \psi_\sharp\mu_k$.
    Using the fact that $\psi$ is strictly increasing, we infer $\abs{\tilde \mu_k}=\psi_\sharp\abs{{\mu_k}}$ and hence $V(\tilde u ) = V(u)$. Taking into account the identity $(\psi\otimes\psi)\inv(\udhp) =\udhp$, we get also $\Qtrans(\tilde u)=\Qtrans(u)$.
     
    Next, $\tilde x\in\R$ is a jump point of $\tilde u $ iff $\psi\inv(\{\tilde x\})$ is nonempty (and so $\psi\inv(\{\tilde x\})=\{x\}$ for some $x$) and $x $   is a jump point of $u$. In that case $\tilde u_{\pm}(\tilde x) = u_{\pm}(x)$. This implies $\Qcost(\tilde u) = \Qcost(u)$.

    Finally, let~$\sigma_k$ (resp.~$\tilde\sigma_k$) be the mapping defined in~\eqref{eq:abs-sigma-map} for~$u$ (resp.~$\tilde u$). We claim that $\sigma_k = \tilde\sigma_k$.

    Indeed, fix $\tau\in(0,\abs{\mu_k}(\R)) = (0,\abs{\tilde\mu_k}(\R))$. Let $x$ (resp. $\tilde x$) be a maximal number such that $\abs{\mu_k}((-\infty, x))\leqslant\tau$ (resp. $\abs{\mu_k}((-\infty, \tilde x))\leqslant\tau$). 
    Then $x$ is the maximal number satisfying $\psi(x)\leqslant \tilde x$. This implies $\psi\inv((\psi(x),\tilde x])=\varnothing$. Thus
    $u_\pm(x) = \tilde u_{\pm}(\tilde x)$, and by~\eqref{eq:abs-sigma-map} we get $\sigma_k(\tau) =\tilde\sigma_k(\tau)$. Consequently, $\Qcubic_k(\tilde u) = \Qcubic_k(u)$, and the proof is complete.
\end{proof}

\begin{remark}  If $\beta\ll 1$ and $u,\tilde u\in\mcD_\beta$, an alternative proof is possible. Namely, notice that Theorem~\ref{thm:reparametr-state-curve} trivially holds true for piecewise constant $u$ and $\tilde u$, and hence $\Upsilon(\Rc\mcF\bullet (\tilde u)) = \Upsilon(\Rc{\psi\inv(\mcF)}\bullet (u)) $ for every $\mcF\in\mcP_\R(\R)$. It remains to apply Theorem~\ref{thm:sup-prop-cut}.
\end{remark}

\bibliographystyle{plain}

\bibliography{references}

\end{document}